\documentclass[12pt,reqno]{amsart}
\usepackage{amsmath,amsthm,amssymb,mathrsfs,stmaryrd,color,mathtools}

\usepackage[all]{xy}
\usepackage{url}
\usepackage{tikz-cd}
\usepackage[margin=1in,footskip=.5in]{geometry}

\usepackage[utf8]{inputenc}
\usepackage[T1]{fontenc}

\usepackage{relsize}
\usepackage[bbgreekl]{mathbbol}
\usepackage{amsfonts}

\DeclareSymbolFontAlphabet{\mathbb}{AMSb}
\DeclareSymbolFontAlphabet{\mathbbl}{bbold}
\newcommand{\Prism}{{\mathlarger{\mathbbl{\Delta}}}}
\usepackage{enumitem}
\usepackage[colorlinks=true,hyperindex, linkcolor=magenta, pagebackref=false, citecolor=cyan,pdfpagelabels]{hyperref}
\usepackage[capitalize]{cleveref}

\theoremstyle{plain}
\newtheorem{thm}{Theorem}[section]
\newtheorem{thm2}{Theorem}
\newtheorem{conj}[thm]{Conjecture}
\newtheorem{conv}[thm]{Convention}

\newtheorem{prop}[thm]{Proposition}
\newtheorem{prop2}[thm2]{Proposition}

\newtheorem{cor}[thm]{Corollary}

\newtheorem{lem}[thm]{Lemma}

\theoremstyle{definition}
\newtheorem{defi}[thm]{Definition}
\newtheorem{rmk}[thm]{Remark}
\newtheorem{exa}[thm]{Example}
\newtheorem{const}[thm]{Construction}
\newtheorem{ass}[thm]{Assumption}

\numberwithin{equation}{section}

\newcommand{\bb}[1]{\mathbb{#1}}
\newcommand{\cl}[1]{{\mathcal{#1}}}
\newcommand{\scr}[1]{{\mathscr{#1}}}
\newcommand{\msf}[1]{{\mathsf{#1}}}
\newcommand{\mfr}[1]{{\mathfrak{#1}}}
\newcommand{\mrm}[1]{{\mathrm{#1}}}

\newcommand{\mbf}[1]{\mathbf{#1}}

\newcommand{\ov}[1]{{\overline{#1}}}

\newcommand{\wtd}[1]{{\widetilde{#1}}}

\newcommand{\Qla}{{\overline{\mathbb{Q}}_\ell}}
\newcommand{\Qlax}{{\overline{\mathbb{Q}}_\ell^\times}}

\newcommand{\LG}{{{}^LG}}

\newcommand{\Proj}{{\operatorname{Proj}}}
\newcommand{\Isom}{{\operatorname{Isom}}}
\newcommand{\Mat}{{\operatorname{Mat}}}

\newcommand{\KGL}{{\operatorname{KGL}}}

\newcommand{\Lie}{{\operatorname{Lie}}}
\newcommand{\Img}{{\operatorname{Im}}}
\newcommand{\clos}{{\operatorname{cl}}}

\newcommand{\ints}{{\operatorname{int}}}

\newcommand{\norm}{{\operatorname{norm}}}

\newcommand{\spc}{{\operatorname{sp}}}
\newcommand{\cosp}{{\operatorname{cosp}}}

\newcommand{\Spd}{{\operatorname{Spd}}}
\newcommand{\Spec}{{\operatorname{Spec}}}
\newcommand{\Spa}{{\operatorname{Spa}}}
\newcommand{\Spf}{{\operatorname{Spf}}}
\newcommand{\Perf}{{\operatorname{Perf}}}

\newcommand{\Sh}{{\operatorname{Sh}}}
\newcommand{\univ}{{\operatorname{univ}}}

\newcommand{\et}{{\acute{e}t}}

\newcommand{\crys}{{\operatorname{crys}}}
\newcommand{\perf}{{\operatorname{perf}}}

\newcommand{\bas}{{\operatorname{bas}}}

\newcommand{\Ker}{{\operatorname{Ker}}}
\newcommand{\Cok}{{\operatorname{Cok}}}

\newcommand{\id}{{\operatorname{id}}}
\newcommand{\an}{{\operatorname{an}}}

\newcommand{\colim}{{\operatorname{colim}}}%{\mathop{\operatorname{colim}}}
\newcommand{\pr}{{\operatorname{pr}}}
\newcommand{\pre}{{\operatorname{pre}}}
\newcommand{\ad}{{\operatorname{ad}}}
\newcommand{\red}{{\operatorname{red}}}

\newcommand{\Hom}{{\operatorname{Hom}}}

\newcommand{\Ind}{{\operatorname{Ind}}}
\newcommand{\cInd}{{\operatorname{cInd}}}

\newcommand{\Fr}{{\operatorname{Fr}}}

\newcommand{\diag}{{\operatorname{diag}}}

\newcommand{\FS}{{\operatorname{FS}}}

\newcommand{\LT}{{\operatorname{LT}}}
\newcommand{\CM}{{\operatorname{CM}}}
\newcommand{\Rep}{{\operatorname{Rep}}}
\newcommand{\Irr}{{\operatorname{Irr}}}

\newcommand{\PGL}{{\operatorname{PGL}}}
\newcommand{\GL}{{\operatorname{GL}}}

\newcommand{\std}{{\operatorname{std}}}
\newcommand{\triv}{{\operatorname{triv}}}
\newcommand{\Sht}{{\operatorname{Sht}}}

\newcommand{\sm}{{\operatorname{sm}}}

\newcommand{\Dr}{{\operatorname{Dr}}}
\newcommand{\Ad}{{\operatorname{Ad}}}

\usepackage[citestyle=alphabetic,bibstyle=alphabetic,backend=bibtex, maxalphanames=10, maxnames=10, url=false, doi=false]{biblatex}
\title{On depth-zero integral models of local Shimura varieties}
\author{Yuta Takaya}
\address{Graduate School of Mathematical Sciences, The University of Tokyo, 3-8-1 Komaba, Meguro-ku, Tokyo 153-8914, Japan}
\email{takaya@ms.u-tokyo.ac.jp}

\begin{document}
\begin{abstract}
    We specify explicit affinoids in depth-zero local Shimura varieties whose reductions are parabolic Deligne-Lusztig varieties, and construct explicit Jacquet-Langlands pairs of regular depth-zero supercuspidal representations in the cohomology of local Shimura varieties. Along the way, we develop the theory of integral moduli spaces of depth-zero level structures on local shtukas, especially at non-parahoric levels. As a consequence, we provide a direct construction of Yoshida's generalized semistable models of depth-zero Lubin-Tate spaces that were originally constructed as successive blowups. 
\end{abstract}

\maketitle
\small 
\setcounter{tocdepth}{2}
\tableofcontents

\section*{Introduction}

\subsection{Background}

The cohomology of local Shimura varieties simultaneously admits smooth $p$-adic group actions and Weil group actions and plays an important role in the local Langlands program. For general linear groups $\GL_n$, the local Langlands correspondence established by Harris-Taylor \cite{HT01} is realized in the cohomology of Lubin-Tate spaces. 

A standard approach to compute the cohomology of local Shimura varieties is via globalization to Shimura varieties. However, this globalization method does not apply to all local Shimura varieties and also lacks an explicit description of Galois representations. 
For Lubin-Tate spaces, there is a more direct approach to compute their cohomology via nearby cycles of \textit{special affinoids}, which answers \cite[Question 9]{Harris02_ICM} depending on types under consideration. Though this approach depends on types, it is purely local and provides the \textit{explicit} local Langlands correspondence. It was first carried out by Yoshida \cite{Yos10} for depth-zero supercuspidal types, and then further developed by Boyarchenko-Weinstein \cite{BW16}, Imai-Tsushima \cite{IT20}, \cite{IT21} and Tokimoto \cite{Tok20} for several other types. 

The strategies in the literature depend on explicit defining equations of special affinoids in Lubin-Tate spaces. Since Rapoport-Zink spaces, or general local Shimura varieties, are defined as moduli spaces of $p$-divisible groups or local shtukas, it is quite difficult to compute their explicit defining equations. In this paper, building on recent developments in $p$-adic geometry, we develop a method that avoids these computations in the case of regular depth-zero supercuspidal types, thereby extending \cite{Yos10} to general local Shimura varieties beyond Lubin-Tate spaces. 

We expect that our approach is also applicable to other types, and in our forthcoming work, we will generalize \cite{BW16} to general local Shimura varieties. 

\subsection{Special affinoids at depth zero}

In this paper, we work with general (basic unramified) local Shimura varieties. Let $(G, [b], \mu)$ be a basic unramified local Shimura datum defined over a local field $F / \bb{Q}_p$. We fix a reductive model $\cl{G}$ of $G$ over $O_F$ and consider the first congruence subgroup $\cl{G}(1) = \Ker(\cl{G}(O_F) \to \cl{G}(k))$. Here, $k$ is the residue field of $F$. As observed in the case of Lubin-Tate spaces, our depth-zero special affinoid
\[
    \cl{W}(0) \subset \cl{M}_{G, b, \mu, \cl{G}(1)}
\]  
is a neighborhood of an unramified CM point (though implicitly in this paper). 

Fix a Borel pair $T \subset B \subset G$ over $F$ and let $\breve{F}$ be the completed maximal unramified extension of $F$. In the following, calligraphy fonts denote integral models and the breve accent denotes the base change to $\breve{F}$. First, as a technically important choice, we take a \textit{length-zero} representative (see \Cref{sssec:length_zero}) of $[b]$ so that
\[
    b = \mu(-\pi) w, \quad w \in  N_{\cl{G}}(\cl{T})(O_{\breve{F}}). 
\]
Let $\breve{M} \subset \breve{G}$ be the centralizer of $(w\sigma)^i \mu$ for all $i \geq 0$. We fix $p_w \in \cl{G}(O_{\breve{F}})$ so that $p_w^{-1} \sigma(p_w) = w$, which provides an embedding 
\begin{center}
    \begin{tikzcd}
        G & \breve{M}^{w\sigma} \ar[l, "\Ad(p_w)"', hook'] \ar[r, hook] & G_b. 
    \end{tikzcd}
\end{center}
Let $\mu_w = \Ad(p_w)(\mu)$ and $b_w = \mu_w(-\pi)$. Then, the local Shimura variety attached to $(\breve{M}^{w\sigma}, b_w, \mu_w)$ is zero-dimensional and the image of a base point to $\cl{M}_{\cl{G}, b, \mu}$ is an unramified CM point $x_\CM$ at hyperspecial level. 

The main idea is to construct $\cl{W}(0)$ inside the \textit{tubular neighborhood} containing the unramified CM point, and the key input is the local representability of the $v$-sheaf theoretic integral model $\cl{M}^\ints_{\cl{G}, b, \mu}$ proved by \cite{Bar22} and \cite{Ito25a}: the specialization of $x_\CM$ in $\cl{M}^\ints_{\cl{G}, b, \mu}$ is $[1] \in X_{\mu}(b) \cong (\cl{M}^\ints_{\cl{G}, b, \mu})^\red$ and the tubular neighborhood at $[1]$ is described as
\[
    (\cl{M}^\ints_{\cl{G}, b, \mu})^\wedge_{/[1]} \cong \Spd(R_{\cl{G}, \mu}), \quad
    R_{\cl{G}, \mu} = O_{\breve{F}}\llbracket u_{\alpha} \vert \alpha \in \Phi_{\mu < 0} \rrbracket. 
\] 
Here, $\Phi_{\mu < 0}$ denotes the set of roots $\alpha$ with $\langle \alpha, \mu \rangle = -1$. By the length-zero condition on $b$, the stabilizer of $[1]$ under $G_b(F)$ is a parahoric subgroup $G_{b, \mfr{f}}$ (see \Cref{lem:fparah}) and it acts on $\Spd(R_{\cl{G}, \mu})$. Now, our special affinoid $\cl{W}(0)$ is described as follows in terms of the parameter $u_\alpha$. 

\begin{thm2}\textup{(\Cref{thm:reduction_of_specaff} and \Cref{prop:Gbact})}
    \label{thm2:prop_specialaff}
    There is a $G_{b, \mfr{f}}$-stable rational subdomain 
    \[
        \cl{U}(0) = \{ \lvert u_\alpha \rvert \leq \lvert \pi \rvert^{r_\alpha} \neq 0 \mid \alpha \in \Phi_{\mu < 0} \} \subset \Spf(R_{\cl{G}, \mu})_\eta, \quad
        0 < r_\alpha < 1
    \]
    such that its $\cl{G}(k)$-cover $\cl{W}(0)$ defined by the Cartesian diagram
    \begin{center}
        \begin{tikzcd}
            \cl{W}(0) \ar[r, hook] \ar[d] & \cl{M}_{G, b, \mu, \cl{G}(1)} \ar[d] \\
            \cl{U}(0) \ar[r, hook] & \cl{M}_{\cl{G}, b, \mu}
        \end{tikzcd}
    \end{center}
    has good reduction over a tamely ramified extension $\breve{F}_e / \breve{F}$ and the reduction $\msf{W}(0)$ is realized as a parabolic Deligne-Lusztig variety associated with maximal reductive quotients of $\cl{G}(O_F)$ and $G_{b, \mfr{f}}$. 
\end{thm2}

Here, $e$ and $r_\alpha$ can be explicitly described (see \Cref{lem:tame_ramification} and \Cref{lem:r_alpha_defi}). Nevertheless, it is difficult to give an explicit defining equation of $\cl{W}(0)$ due to the moduli-theoretic nature of the cover $\cl{W}(0) \to \cl{U}(0)$. Instead, we first introduce a parabolic Deligne-Lusztig variety $\msf{W}(0)$ as a $\cl{G}(k)$-cover of the reduction of $\cl{U}(0)$ and lift it to a $\cl{G}(k)$-cover of $\cl{U}(0)$. Then, we compare it to $\cl{W}(0)$ via the moduli interpretation of $\msf{W}(0)$. 

Another difficult part in the proof is to verify the stability of $\cl{U}(0)$ under $G_{b, \mfr{f}}$. For this, we rely on an auxiliary \textit{depth-zero integral model}, which is another important topic of this paper. 

By studying group actions on $\cl{W}(0)$, we may apply the uniform procedure explained in \cite{Mie16} to see that $\cl{W}(0)$ provides the following contribution to the middle \'{e}tale cohomology $H_c^r(G, b, \mu)$ of local Shimura varieties. Let $E$ be the reflex field of $(G, \mu)$. 

\begin{thm2}\textup{(\Cref{prop:computcoh} and \Cref{prop:computsplit})} \label{thm2:cohomology_computation}
    For every depth-zero unramified elliptic regular pair $(S, \theta)$ with $S \subset \breve{M}^{w\sigma}$, there is a finite-dimensional smooth representation $\sigma_{(S, \theta)}$ of $W_E$ such that 
    \[
        \pi_{(S, \theta)} \boxtimes \rho_{(S, \theta)}^\vee \boxtimes \sigma_{(S, \theta)} \subset H^r_c(G, b, \mu)_{\chi \boxtimes \chi^{-1}}
    \]
    for $\chi = \theta\vert_{Z_G(F)}$. Let $d \geq 1$ be the minimum positive integer such that $w^d \in \breve{T}$ and let $F_d$ be the unramified extension of $F$ of degree $d$. When $G$ is split, $\sigma_{(S, \theta)}$ can be taken as 
    \[
        \sigma_{(S, \theta)} = \Ind_{W_{F_d}}^{W_F} \xi_{(S, \theta)}
    \]
    for a certain explicit tame character $\xi_{(S, \theta)} \colon W_{F_d} \to \Qlax$. 
\end{thm2}

Here, $\pi_{(S, \theta)}$ and $\rho_{(S, \theta)}$ are depth-zero regular supercuspidal representations of $G(F)$ and $G_b(F)$, respectively, attached to the pair $(S, \theta)$. Moreover, $H^r_c(G, b, \mu)_{\chi \boxtimes \chi^{-1}}$ denotes the coinvariant with respect to the central character. 

The inertial part of $\xi_{(S, \theta)}$ is given by $\theta$ via local class field theory, and its Frobenius value is a Frobenius eigenvalue of the middle \'{e}tale cohomology of $\msf{W}(0)$. For the latter, we leave the exact computation due to the lack of references on Frobenius eigenvalues of parabolic Deligne-Lusztig varieties. In the case studied in \cite{Yos10}, \Cref{thm2:cohomology_computation} is compatible with the computation in \cite{Mie16}. 

%Here, $\pi_{(S, \theta)}$ and $\rho_{(S, \theta)}$ are depth-zero regular supercuspidal representations 
%\[
%    \pi_{(S, \theta)} = \cInd_{Z_G(F) \cl{G}(O_F)}^{G(F)} \pi_{(S, \theta)}^0, \quad
%    \rho_{(S, \theta)} = \cInd_{Z_G(F) G_{b, \mfr{f}}}^{G(F)} \rho_{(S, \theta)}^0. 
%\]
%Let $\msf{G}^\sigma = \cl{G}(k)$ and $\msf{M}^{w\sigma} = \cl{\breve{M}}^{w\sigma}(k)$. The maximal reductive quotient of $G_{b, \mfr{f}}$ is $\msf{M}^{w\sigma}$, and $\pi_{(S, \theta)}^0$ and $\rho_{(S, \theta)}^0$ are related by the Deligne-Lusztig induction along
%\[
%    \msf{G}^\sigma \times \msf{M}^{w\sigma} \curvearrowright \msf{W}(0). 
%\]
%In particular, our choice of $(S, \theta)$ ensures that one of the maximal reductive quotients appearing in types is a twisted Levi subgroup of the other. Without the length-zero condition on $b$, it is not always the case (add an example of $\GSp_4$). 

Let us mention a relation to \cite{HKW22}. Let 
\[
    R\Gamma(G,b,\mu)[\pi_{(S, \theta)}] = R\Hom_{G(F)}(\colim_{K \subset G(F)} R\Gamma_c(\cl{M}_{G,b,\mu,K,\bb{C}_p}, \ov{\bb{Q}}_\ell[r](\tfrac{r}{2})), \pi_{(S, \theta)})^{\sm}. 
\]
When $p$ and $q$ are sufficiently large so that endoscopic character identities are available for the $L$-packets of $\varphi_{(S, \theta)}$ (see \cite[Section 3.4]{Kal11}), \cite[Theorem 1.0.2]{HKW22} computes the Euler characteristic of $R\Gamma(G,b,\mu)[\pi_{(S, \theta)}]$ up to parabolically induced representations. It implies that $\rho_{(S, \theta)}$ contributes to the Euler characteristic. \Cref{thm2:cohomology_computation} shows that $\rho_{(S, \theta)}$ contributes to $H^0R\Gamma(G, b, \mu)[\pi_{(S, \theta)}]$, which is a particular case of the vanishing conjecture of the cohomology of local Shimura varieties. 

\subsection{Depth-zero integral models}

We mentioned that we need auxiliary depth-zero integral models of tubular neighborhoods for \Cref{thm2:prop_specialaff}. In this paper, we develop a general construction of so-called depth-zero integral models, which classify depth-zero level structures of local shtukas. 

In the classical setting, one usually resorts to Drinfeld level structures, resulting in integral models finite over some parahoric levels, but such models do not have good singularities. Instead, by passing to local shtukas (via the Dieudonn\'{e} theory for $p$-divisible groups), our construction provides depth-zero integral models \textit{non-quasi-finite} over the base. The representability of such moduli spaces is nontrivial since we may not resort to the algebraicity of $p$-divisible groups anymore, but such a difficulty has been already overcome by the author's representability criterion \cite{Tak26_rel}.  

The basic idea is as follows. Let $X$ be a $p$-divisible group of height $n$ over a perfectoid ring $R$. To any homomorphism $\alpha \colon \bb{F}_p^{\oplus n} \to X[p]$, the torsion Dieudonn\'{e} theory (see \cite{ALB23} for perfectoid rings) attaches a homomorphism $\bb{D}(\alpha) \colon (R^\flat)^{\oplus n} \to \bb{D}(X) / p \bb{D}(X)$. If we set $\cl{P}_X$ to be the local $\GL_n$-shtuka associated to $X$, then $\bb{D}(\alpha)$ is identified with a section of $\cl{P}_X\vert_{R^\flat} \times^{\GL_n} \Mat_n$. By replacing $\GL_n \subset \Mat_n$ with an arbitrary \textit{equivariant compactification} $\cl{G}_k \subset \msf{\ov{G}}$, we may think of a section of $\cl{P}\vert_{\Spa(R^\flat, R^{\flat+})} \times^{\cl{G}_k} \msf{\ov{G}}$ as a \textit{depth-zero level structure} on a local $\cl{G}$-shtuka $\cl{P}$ over a perfectoid space $\Spa(R, R^+)$. 

In fact, a more general setting is allowed in our construction. Let $\cl{G}$ be an arbitrary parahoric group scheme and let $\cl{G}^+$ be an arbitrary dilatation of $\cl{G}$ such that $\msf{G}^\sigma = \cl{G}(O_F) / \cl{G}^+(O_F)$ is reductive (see \Cref{ssec:level_structure_in_char0}). We fix an equivariant compactification $\msf{G}^\sigma \subset \msf{\ov{G}}$. 

\begin{prop2} \textup{(\Cref{prop:replevel})}
    \label{prop2:replevel}
    Let $\mfr{X}$ be a reduced excellent $p$-adic formal $O_F$-scheme admitting a good cover and let $\cl{P}$ be a local $\cl{G}$-shtuka over $\mfr{X}^\diamond _{/O_F}$. The $v$-closure of the embedding
    \begin{equation} \label{eq:proper_embedding}
        \Sht_{\cl{G}^+}\times_{\Sht_{\cl{G}}, \cl{P}} \mfr{X}^{\diamond}_\eta \subset \lbrack \msf{\ov{G}}^{\diamond} /(L_W^+\cl{G})^\diamondsuit \rbrack \times_{\lbrack \ast/(L_W^+\cl{G})^\diamondsuit \rbrack, \cl{P}} \mfr{X}^\diamond
    \end{equation}
    is representable by a unique proper $p$-adic formal scheme $\mfr{Y}$ over $\mfr{X}$ admitting a maximal good cover with $\mfr{Y}_\eta$ finite \'{e}tale over $\mfr{X}_\eta$. Moreover, the $\msf{G}^\sigma$-action on $\mfr{Y}_\eta$ over $\mfr{X}_\eta$ extends to $\mfr{Y}$ over $\mfr{X}$. 
\end{prop2}

The terminology `good cover' is introduced in \cite{Tak26_rel} and the construction of \eqref{eq:proper_embedding} is explained in \Cref{ssec:moduli_level_structures}. In \Cref{sec:level_structure_pdivgrp}, we study the relation of this construction to classical constructions: when $\cl{G} = \GL_n / \bb{Z}_p$ and $\cl{P}$ comes from a $p$-divisible group of height $n$, $\mfr{Y}$ (essentially) recovers the flat closure of the moduli of depth-zero Drinfeld level structures by using a compactification
\[
    \msf{\ov{G}} = \bb{P}(\Mat_n \oplus \triv) \supset \Mat_n \supset \GL_n = \msf{G}. 
\]
When $\cl{G} = \cl{I} / \bb{Z}_p$ is the Iwahori group scheme of $\GL_n$ and $\cl{P}$ comes from a full flag of $p$-divisible groups of height $n$, $\mfr{Y}$ (essentially) recovers the flat closure of the moduli of $\Gamma_1(p)$-level structures as introduced in \cite{HR12} by using a compactification
\[
    \msf{\ov{I}} = (\bb{P}^1)^n \supset \bb{A}^n \supset \bb{G}_m^n = \msf{I}. 
\]

Now, we apply this construction to the universal deformation over $R_{\cl{G}, \mu}$. By using a \textit{toroidal compactification} $\msf{\ov{G}}$, we get the following integral model of the tubular neighborhood at $[1]$. Let $\mfr{W}(0)$ be the smooth affine formal model of $\cl{W}(0) \otimes \breve{F}_e$. 

\begin{prop2}\textup{(\Cref{thm:role_of_integralmodel})}
    \label{prop2:characterization_intmodel}
    There is a proper formal scheme $\mfr{X}_b$ over $R_{\cl{G}, \mu}$ with a $\cl{G}(k)$-action such that we have a Cartesian diagram 
    \begin{center}
        \begin{tikzcd}
            \mfr{X}_{b, \eta} \ar[r, hook] \ar[d] & \cl{M}_{G, b, \mu, \cl{G}(1)} \ar[d] \\
            \Spf(R_{\cl{G}, \mu})_\eta \ar[r, hook] & \cl{M}_{\cl{G}, b, \mu}
        \end{tikzcd}
    \end{center}
    and $\cl{W}(0) \cong (\mfr{X}_b\vert_{U_b})_\eta$ for some affine open subset $U_b \subset \lvert \mfr{X}_b \rvert$. Moreover, the natural map
    $
        \mfr{W}(0) \to \mfr{X}_b\vert_{U_b}
    $
    is finite and surjective. 
\end{prop2}

In particular, $\mfr{X}_b$ admits a $p$-adic open formal subscheme whose underlying space is $r$-dimensional, though $R_{\cl{G}, \mu}$ is far from $p$-adic and $\lvert \Spf(R_{\cl{G}, \mu}) \rvert$ is a singleton. Here, we find another interesting feature: there is a natural closed embedding $\mfr{X}_{b, \red}^\perf \subset \msf{\ov{G}}^\perf$ and the image of the composition
\[
    \msf{W}(0)^\perf \to \mfr{X}_{b, \red}^\perf \subset \msf{\ov{G}}^\perf
\]
recovers a natural embedding of a parabolic Deligne-Lusztig variety $\msf{W}(0)$ into a flag variety. 

As in \cite{Tak26_rel}, we may apply \Cref{prop2:replevel} to canonical integral models $\{\scr{S}_{\msf{K}}\}_{\msf{K}^p}$ of Shimura varieties to get depth-zero integral models. Let $(\mbf{G}, \mbf{X})$ be a Shimura datum satisfying the axiom (SV5) and let $G = \mbf{G}_{\bb{Q}_p}$. Let $\msf{K}_p = \cl{G}(\bb{Z}_p)$ and let $\msf{K}=\msf{K}_p\msf{K}^p\subset \mbf{G}(\bb{A}_f)$ be a compact open subgroup for a sufficiently small compact open subgroup $\msf{K}^p\subset \mbf{G}(\bb{A}_f^p)$. Let $\mbf{E}$ be the reflex field of $(\mbf{G},\mbf{X})$. Let $v$ be a place of $\mbf{E}$ over $p$ and set $E=\mbf{E}_v$. 

Now, suppose that $\msf{K}_p$ is contained in a hyperspecial subgroup $\msf{K}'_p$ and there exists a system of canonical integral models $\{\scr{S}_{\msf{K}'}\}_{\msf{K}^p}$ of Shimura varieties in the sense of \cite{PR24}. Here, $\msf{K}' = \msf{K}'_p \msf{K}^p$ and the latter assumption is verified in many contexts (see \Cref{rmk:PR_canonical_models}). Then, by \cite[Theorem 6.32]{Tak26_rel}, there exists a system of canonical integral models $\{\scr{S}_{\msf{K}}\}_{\msf{K}^p}$ of Shimura varieties as algebraic spaces. Let $\cl{P}_\msf{K}$ be the local $\cl{G}$-shtuka over $\scr{S}_\msf{K}$. Then, we get the following depth-zero integral models $\scr{S}_{\msf{K}^+}$. Here, $\msf{K}_p^+ = \cl{G}^+(\bb{Z}_p)$ and $\msf{K}^+ = \msf{K}^+_p \msf{K}^p$. 

\begin{prop2}\textup{(\Cref{prop:depth_zero_integral_models})}
    In the above situation, for every $\msf{K}^p$, there is a flat separated algebraic space $\scr{S}_{\msf{K}^+}^{\msf{\ov{G}}}$ of finite type over $O_E$ with $(\scr{S}_{\msf{K}^+}^{\msf{\ov{G}}})_E \cong \mrm{Sh}_{\msf{K}^+}(\mbf{G},\mbf{X})_E$ such that 
    \begin{enumerate}
        \item The $\msf{G}^\sigma$-action on $\mrm{Sh}_{\msf{K}^+}(\mbf{G},\mbf{X})_E$ over $\mrm{Sh}_{\msf{K}}(\mbf{G},\mbf{X})_E$ extends to $\scr{S}_{\msf{K}^+}^{\msf{\ov{G}}}$ over $\scr{S}_{\msf{K}}$, 
        \item $\scr{S}_{\msf{K}^+}^{\msf{\ov{G}}}$ is proper and surjective over $\scr{S}_{\msf{K}}$, and
        \item $(\scr{S}_{\msf{K}^+}^{\msf{\ov{G}}})^\wedge$ is the flat moduli space of $\cl{G}^+$-level structures of type $\msf{\ov{G}}$ on $\cl{P}_\msf{K}$. 
    \end{enumerate}
    Moreover, the system $\{\scr{S}_{\msf{K}^+}^{\msf{\ov{G}}}\}_{\msf{K}^p}$ admits finite \'{e}tale prime-to-$p$ Hecke correspondences and satisfies
    \[
        (\varprojlim_{\msf{K}^p} \Sh_{\msf{K}^+}(\mbf{G},\mbf{X}))(R[\tfrac{1}{p}]) = (\varprojlim_{\msf{K}^p} \scr{S}_{\msf{K}^+})(R)
    \]
    for every discrete valuation ring $R$ of mixed characteristic over $O_E$. 
\end{prop2}

In a recent preprint, Pappas-Rapoport \cite{PR26} proposed axioms of $\Gamma_1(p)$-level integral models \cite[Conjecture 1.6.1]{PR26}, where $\cl{G}$ is Iwahori and $\msf{G}^\sigma = \cl{G}(O_F) / \cl{G}^+(O_F)$ is a torus. Our construction satisfies conditions (i) and (iii) loc. cit. and recovers classical $\Gamma_1(p)$-level integral models of unramified unitary Shimura varieties up to absolute weak normalization (see \Cref{exa:Gamma_1(p)_unitary}). It would be interesting to clarify the relation to the construction in \cite[Theorem 1.6.2]{PR26}. 

In classical settings, $\scr{S}_{\msf{K}^+}$ is always finite over $\scr{S}_{\msf{K}}$, but one should not expect such finiteness in our construction, at least outside $\Gamma_1(p)$-levels. In fact, we expect that Yoshida's successive blowup of Harris-Taylor integral models (see \cite[Definition 4.5]{Yos10}), which has mild singularities \cite[Remark 4.9]{Yos10}, should be recovered by our construction (see \Cref{exa:Harris_Taylor}). In this paper, we verify a local variant of this expectation for Lubin-Tate spaces.

Let $\LT_n$ be the Lubin-Tate space associated to $\cl{G} = \GL_n / \bb{Z}_p$ and let $\cl{G}^+$ be its first congruence subgroup scheme. As a toroidal compactification of $\GL_n$, we use the one constructed in \cite{Kau00}, denoted by $\KGL_n$. 

\begin{thm2}\textup{(\Cref{thm:comparison})}
    Let $\LT_n^{\KGL_n}$ be the flat moduli space of $\cl{G}^+$-level structures of type $\KGL_n$ on the universal $p$-divisible group over $\LT_n$. Then, the normalization of $\LT_n^{\KGL_n}$ is isomorphic to Yoshida's generalized semistable model $\LT_n^{1, n-1}$ of the depth-zero Lubin-Tate space. 
\end{thm2}

Though $\LT_n$ is not $p$-adic, $\LT_n^{\KGL_n}$ can be naturally defined (see \Cref{defi:LTn_KGLn}). As a direct computation of $\LT_n^{\KGL_n}$ is difficult, we instead compute its specialization map (see \Cref{sssec:spKGLn}). Then, we compare it with the one for $\LT_n^{1, n-1}$, which we compute in \Cref{sssec:spcLT}, and invoke \cite[Corollary 6.6]{Lou17} to get the claim. 

\addtocontents{toc}{\protect\setcounter{tocdepth}{1}}

\subsection*{The structure of the paper}
In \Cref{sec:introlocsht}, we review the theory of local shtukas and introduce depth-zero level structures in characteristic $0$. In \Cref{sec:locShim}, we explain length-zero conditions and review an explicit description of tubular neighborhoods. In \Cref{sec:DL}, we specify special affinoids at depth zero and show that their reductions are parabolic Deligne-Lusztig varieties. In \Cref{sec:dep0_integral_model}, we introduce depth-zero integral models and provide the characterization of special affinoids in terms of integral models. In \Cref{sec:group_action_on_special_affinoids}, we use the characterization to show the stability of special affinoids under several group actions. In \Cref{sec:nearby}, we review the nearby cycle method to special affinoids relating the cohomology of reductions and generic fibers. In \Cref{sec:dep0reg}, we construct explicit Jacquet-Langlands pairs of regular depth-zero supercuspidal representations in the cohomology of local Shimura varieties. In \Cref{sec:level_structure_pdivgrp}, we apply our construction to $p$-divisible groups and compare it with classical constructions in terms of Drinfeld level structures. In \Cref{sec:comparison_Yoshida}, we show that our construction recovers Yoshida's generalized semistable models after normalization. 

\subsection*{Acknowledgements}
I would like to thank my advisor Yoichi Mieda for his constant support and encouragement. I would also like to thank Naoki Imai, Kai-Wen Lan and Kazuki Tokimoto for helpful discussions. This work was supported by the WINGS-FMSP program at the Graduate School of Mathematical Sciences, the University of Tokyo and JSPS KAKENHI Grant number JP24KJ0865.

\subsection*{Notation}
Fix a prime number $p$. All rings are assumed to be commutative. The reduced underlying scheme of a formal scheme $\mfr{X}$ is denoted by $\mfr{X}_\red$. The perfection of a scheme $X$ over $\bb{F}_p$ is denoted by $X^\mrm{perf}$. For a scheme $X$ over a ring $A$, its base change to an $A$-algebra $B$ is denoted by $X_B$ or $X\otimes_A B$. 

The center of an algebraic group $G$ is denoted by $Z_G$. The normalizer of a subgroup $H$ in $G$ is denoted by $N_G(H)$. For a cocharacter $\lambda \colon \bb{G}_m \to G$ over a ring $A$, $\lambda(a)$ is also denoted by $a^\lambda$ for $a\in A^\times$.
Let $G^{\lambda=0} \subset G$ denote the centralizer of $\lambda$ and let $G^{\lambda \geq 0} \subset G$ denote the parabolic subgroup
\[
	G^{\lambda \geq 0} = \{ g\in G \mid \lim_{t \to 0} \lambda(t) g \lambda(t)^{-1} \hspace{3pt} \mrm{exists}\}. 
\]
Let $G^{\lambda > 0}$ denote the unipotent radical of $G^{\lambda \geq 0}$. For a torus $T$, $X^*(T)$ (resp.\ $X_*(T)$) denotes the character (resp.\ cocharacter) group of $T$ and the natural pairing between $X^*(T)$ and $X_*(T)$ is denoted by $\langle - , - \rangle$. 

Let $F$ be a non-archimedean local field over $\bb{Q}_p$ and let $O_F$ be the ring of integers of $F$ with residue field $k$. Let $\pi$ be a uniformizer of $F$ and let $q$ be the number of elements of $k$. For a perfect $k$-algebra $R$, we set $W_{O_F}(R)=W(R)\otimes_{W(k)} O_F$. The Frobenius on $W_{O_F}(R)$ relative to $k$ is denoted by $\sigma$. The Weil group of $F$ is denoted by $W_F$. The inertia group of $F$ is denoted by $I_F$. The dual group of a connected reductive group $G$ over $F$ is denoted by $\widehat{G}$, which is defined over $\ov{\bb{Q}}_\ell$ unless otherwise stated. The $L$-group of $G$ is denoted by $\LG = \widehat{G} \rtimes W_F$. 

Let $\breve{F}$ be the completed maximal unramified extension of $F$ and let $\ov{k}$ be the residue field of $\breve{F}$. In this paper, italic, calligraphy, sans-serif types denote generic fibers over $F$, integral models over $O_F$, and special fibers over $\ov{k}$, and breve accents denote the base change to $O_{\breve{F}}$ or $\breve{F}$. For a $\ov{k}$-variety $\msf{X}$ with a Frobenius automorphism $\sigma$ relative to $k$, $\msf{X}^{\sigma}$ denotes the set of $\sigma$-fixed points of $\msf{X}(\ov{k})$. Depending on the context, $\msf{X}^{\sigma}$ also denotes the associated $k$-rational form of $\msf{X}$. For a $\ov{k}$-scheme $\msf{X}$ and a $\ov{k}$-algebra $A$, the Frobenius pullback $\msf{X}(A) \to (\sigma^*\msf{X})(A)$ is simply denoted by $\sigma$. When $\msf{X}$ is defined over $k$, $\sigma^*\msf{X}$ is identified with $\msf{X}$. 

Let $R$ be a perfectoid ring over $O_F$. Its tilt is denoted by $R^\flat$ and the Frobenius on $W_{O_F}(R^\flat)$ is denoted by $\varphi_R$. Let
\[
    \theta_R\colon W_{O_F}(R^\flat)\to R
\]
be the canonical surjection and let $\xi_R$ denote a generator of $\Ker(\theta_R)$. The image of $\xi_R$ under $W_{O_F}(R^\flat) \to R^\flat$ is denoted by $\xi_{R,0}$. When we just say that $R$ is a perfectoid ring, we apply the above notation for $O_F = \bb{Z}_p$. 
Let $\Perf$ be the category of perfectoid spaces over $k$ and let $\ast = \Spd(k)$. When we work with the category $\Perf_{\ov{k}}$ of perfectoid spaces over $\ov{k}$, we set $\ast = \Spd(\ov{k})$ by abuse of notation. 
The $v$-sheaf represented by a $p$-adic ring $R$ is denoted by $\Spd(R)$. For a $v$-sheaf $X$ over $\Spd(\bb{Z}_p)$, let $X_s = X \times_{\Spd(\bb{Z}_p)} \Spd(\bb{F}_p)$ and $X_\eta = X \times_{\Spd(\bb{Z}_p)} \Spd(\bb{Q}_p, \bb{Z}_p)$. 

\addtocontents{toc}{\protect\setcounter{tocdepth}{2}}

\section{Preliminaries on local shtukas} \label{sec:introlocsht}

\subsection{Generalities} \label{ssec:introlocsht}

In this section, we set up some notation and recall the definition of local shtukas introduced in \cite{SW20}. 

For each $S \in \Perf$, there is an analytic adic space $\cl{Y}_S$. When $S = \Spa(R,R^+)$ with a pseudo-uniformizer $\varpi \in R^+$, we have
\[
    \cl{Y}_S = \Spa(W_{O_F}(R^+)) - V([\varpi]). 
\]
For each integer $n \geq 1$, let 
\[
    \cl{Y}_{S,[0,n]} = \{ \lvert \pi^n \rvert \leq \lvert [\varpi] \rvert \neq 0 \} \subset \cl{Y}_S
\]
be the rational domain. By the proof of \cite[Proposition II.1.1]{SW20}, we have
\begin{equation}
    \cl{Y}_{S, [0,n]} = \Spa(B_{S,[0,n]},B^+_{S,[0,n]}), \quad B_{S, [0, n]}^\wedge \cong W_{O_F}(R). \label{eq:YSn}
\end{equation}
Here, the wedge denotes the $p$-adic completion. Let $G$ be a connected reductive group over $F$ and let $\cl{G}$ be a smooth affine model of $G$ with connected fibers over $O_F$. 

\begin{defi}(\cite[Section 23.1]{SW20})
    For each $S \in \Perf$, a local $\cl{G}$-shtuka over $S$ is a tuple 
    \[
        ((S^\sharp,\iota), \cl{P}, \varphi_{\cl{P}})
    \]
    where $(S^\sharp, \iota)\in \Spd(O_F)(S)$ is an untilt of $S$ over $O_F$, $\cl{P}$ is a $\cl{G}$-torsor over $\cl{Y}_S$ and 
    \[
        \varphi_{\cl{P}}\colon \varphi^*\cl{P}\vert_{\cl{Y}_S\backslash S^\sharp} \cong \cl{P}\vert_{\cl{Y}_S\backslash S^\sharp}
    \]
    is an isomorphism of $\cl{G}$-torsors meromorphic along $S^\sharp$. The $v$-stack of local $\cl{G}$-shtukas on $\Perf$ is denoted by $\Sht_{\cl{G}}$. 
\end{defi}

\begin{defi}
    Let $X$ be a $v$-stack on $\Perf$. A local $\cl{G}$-shtuka over $X$ is a map
    \[
        X \to \Sht_{\cl{G}}
    \]
    of $v$-stacks. When $X$ is a $v$-stack over $\Spd(O_F)$, a local $\cl{G}$-shtuka over $X_{/O_F}$ is a map $X \to \Sht_{\cl{G}}$ over $\Spd(O_F)$. 
\end{defi}

The restriction along \eqref{eq:YSn} provides the following map. 

\begin{lem}\label{prop:mapLW}
    Let $(L_W^+\cl{G})^\diamondsuit$ be the $v$-sheaf on $\Perf$ sending $\Spa(R,R^+)$ to $\cl{G}(W_{O_F}(R))$. There is a natural map 
    \[
        \Sht_{\cl{G}}\to [\ast/(L_W^+\cl{G})^\diamondsuit]. 
    \]
\end{lem}
\begin{proof}
    Since $\cl{G}$ is affine, $(L_W^+\cl{G})^\diamondsuit$ is a small $v$-sheaf on $\Perf$ by \cite[Theorem 8.7]{Sch17}. The quotient $v$-stack $[\ast/(L_W^+\cl{G})^\diamondsuit]$ classifies the groupoid of $\cl{G}$-torsors on $\Spec(W_{O_F}(R))$ over an affinoid perfectoid space $S=\Spa(R,R^+)$ over $k$ since $\cl{G}$-torsors on $W_{O_F}(R)$ satisfy $v$-descent by \cite[Corollary 17.1.9]{SW20}. As in \cite[Section 20.3]{SW20}, we obtain a map $\Sht_{\cl{G}}\to [\ast/(L_W^+\cl{G})^\diamondsuit]$ by pulling back $\cl{G}$-torsors on $\cl{Y}_{S,[0,n]}$ to $W_{O_F}(R)$. 
\end{proof}

%Let $R$ be a perfectoid ring over $O_F$. The Frobenius map on $W_{O_F}(R^\flat)$ is denoted by $\varphi_R$. By \cite[Lemma 2.4.2]{Ito25b}, the kernel of the surjection $\theta_R\colon W_{O_F}(R^\flat)\to R$ is principal. By an abuse of notation, a generator of $\Ker(\theta_R)$ is also denoted by $\xi_R$. The image of $\xi_R$ under the projection $W_{O_F}(R^\flat) \to R^\flat$ is denoted by $\xi_{R,0}$. 

\begin{defi}
    A $\cl{G}$-Breuil-Kisin-Fargues ($\cl{G}$-BKF) module over a perfectoid ring $R$ over $O_F$ is a pair $(P,\varphi_P)$ such that $P$ is a $\cl{G}$-torsor on $W_{O_F}(R^\flat)$ and $\varphi_P$ is an isomorphism
    \[
        \varphi_P \colon (\varphi_R^*P)[1/\xi_R] \cong P[1/\xi_R]. 
    \]
\end{defi}

\begin{prop}\textup{(\cite{Gut23})}\label{prop:GBKF}
    For a perfectoid ring $R$ over $O_F$, the category of local $\cl{G}$-shtukas over $\Spd(R)_{/O_F}$ is equivalent to the category of $\cl{G}$-BKF modules over $R$. 
\end{prop}

The above result is also proved in \cite{GIZ25} when $R$ is in characteristic $p$. For each perfectoid Huber pair $(S,S^+)$ over $R$ and a $\cl{G}$-BKF module $(P, \varphi_P)$, the associated local $\cl{G}$-shtuka over $(S,S^+)$ is
$
    (P\vert_{\cl{Y}_S}, \varphi_{P}\vert_{\cl{Y}_S}). 
$
When $P$ is a trivial $\cl{G}$-torsor, a $\cl{G}$-BKF module is given by
\[
    (\cl{G} \otimes W_{O_F}(R^\flat), g\sigma), \quad
    g \in \cl{G}(W_{O_F}(R^\flat)[1 / \xi_R]). 
\]

\subsection{Depth-zero level structures in characteristic $0$} \label{ssec:level_structure_in_char0}

In this section, we introduce depth-zero level structures of local $\cl{G}$-shtukas on perfectoid spaces over $F$. Here, let $G$ be a connected reductive group over $F$ and let $\cl{G}$ be a parahoric group scheme of $G$. 

Let $\msf{G}^\sigma$ be a reductive quotient of $\cl{G}_{k}$ such that $\Ker(\cl{G}_k \to \msf{G}^\sigma)$ is smooth and connected, and let $\msf{G} = (\msf{G}^\sigma)_{\ov{k}}$. Let $\cl{G}^+$ be the dilatation of $\cl{G}$ along $\Ker(\cl{G}_k \to \msf{G}^\sigma)$ (see \cite[Appendix A]{KP23}). Then, $\cl{G}^+$ is a smooth affine group scheme with connected fibers such that
\[
    \cl{G}^+(O_F) = \Ker(\cl{G}(O_F) \to \msf{G}^\sigma)
\]  
(\cite[Lemma A.5.7, Lemma A.5.10]{KP23}). We set $\cl{G}(1) = \cl{G}^+(O_F)$. The natural map
\[
    \Sht_{\cl{G}^+} \to \Sht_{\cl{G}}
\]
is a finite \'{e}tale $\msf{G}^\sigma$-torsor over $\Spd(F)$ by \cite[Proposition 22.6.1]{SW20}. 

\begin{rmk} \label{rmk:choice_G+}
    Aside from the general construction in \Cref{ssec:moduli_level_structures}, we will assume that $\cl{G}$ is reductive and $\msf{G}^\sigma = \cl{G}_k$. Note that $\msf{G}^\sigma$ also denotes the set of $k$-valued points $\msf{G}^\sigma(k)$ depending on the context. % when there is no confusion.
\end{rmk}

Due to this convention, we will use the following terminology omitting $\cl{G}^+$. 

\begin{defi}
    Let $X$ be a $v$-stack over $\Spd(F)$ and let $\cl{P}$ be a local $\cl{G}$-shtuka over $X$. A depth-zero level structure on $\cl{P}$ is a commutative diagram
    \begin{center}
        \begin{tikzcd}
            X \ar[r, "\cl{P}^+"] \ar[rd, "\cl{P}"'] & \Sht_{\cl{G}^+} \ar[d] \\
            & \Sht_{\cl{G}}. 
        \end{tikzcd}
    \end{center}
    In other words, it consists of a local $\cl{G}^+$-shtuka $\cl{P}^+$ over $X$ and an identification $\cl{P}^+\times^{\cl{G}^+} \cl{G} \cong \cl{P}$. 
\end{defi}

Since $\Sht_{\cl{G}^+, \eta} \to \Sht_{\cl{G}, \eta}$ is a $\msf{G}^\sigma$-torsor, $\msf{G}^\sigma$ acts on the set of depth-zero level structures. This $\msf{G}^\sigma$-set can be described explicitly. 

\begin{lem} \label{lem:Gpsubtors}
    Let $A$ be a $p$-torsion free and $p$-complete $O_F$-algebra. For every $\cl{G}$-torsor $P$ over $A$, there is a bijection between the set of $\cl{G}^+$-subtorsors of $P$ and the set of sections of $P_{A/\pi A} \times^{\cl{G}_k} \msf{G}^\sigma$. 
\end{lem}
\begin{proof}
    Suppose that we have a section of $P_{A/\pi A} \times^{\cl{G}_k} \msf{G}^\sigma$. Since there are no nontrivial automorphisms of $\cl{G}^+$-subtorsors of $P$, we may localize $A$ $p$-completely \'{e}tale locally to get a $\cl{G}^+$-subtorsor of $P$. Thus, we may assume that $P$ is trivial and the section of $P_{A/\pi A} \times^{\cl{G}_k} \msf{G}^\sigma$ can be lifted to $x\in P(A)$. 

    The $\cl{G}^+$-subtorsor generated by $x$ is independent of the choice of $x$ because any other choice $x'\in P(A)$ can be written as $x' = x g^+$ with
    \[
        g^+ \in \Ker(\cl{G}(A) \to \msf{G}^\sigma(A / \pi A)) = \cl{G}^+(A). 
    \]
    Thus, we may attach a $\cl{G}^+$-subtorsor $\cl{G}^+\cdot x \to \cl{G}$ to each liftable section of $P_{A/\pi A} \times^{\cl{G}_k} \msf{G}^\sigma$, and by descent, we get a desired map in one direction. For the converse direction, let $P^+ \to P$ be a $\cl{G}^+$-subtorsor. By localizing $A$ $p$-completely \'{e}tale locally, we may assume that $P^+$ is trivial. Then, every section of $P^+$ maps to the same element of $(P \times^{\cl{G}_k} \msf{G}^\sigma)(A / \pi A)$. By descent, we get the converse map. 
    %Suppose that we have a section of $P_{A/\pi A}$. Since $A$ is $p$-complete and $P$ is smooth and affine, it lifts to an element $x\in P(A)$. Then, the $\cl{G}^+$-subtorsor generated by $x$ is independent of the choice of $x$ because any other choice $x'\in P(A)$ can be written as $x' = x g^+$ with
    %\[
    %    g^+ \in \Ker(\cl{G}(A) \to \cl{G}(A / \pi A)) = \cl{G}^+(A). 
    %\]
    %Thus, we get a map in one direction. For the converse direction, let $P^+ \subset P$ be a $\cl{G}^+$-subtorsor. By localizing $A$ $p$-completely \'{e}tale locally, we may assume that $P^+$ can be trivialized. Then, every section of $P^+$ maps to the same element in $P(A / \pi A)$. By descent, we get the converse map. 
\end{proof}

\begin{lem} \label{lem:finetdesc}
    Let $(R,R^+)$ be a perfectoid Huber pair over $F$ and let $\cl{P}$ be a local $\cl{G}$-shtuka over $(R,R^+)$. Let $P$ be the $\msf{G}^\sigma$-torsor over $R^\flat$ with a Frobenius action $\varphi_P$ obtained as 
    \[
        P = \cl{P}\vert_{\Spa(R^\flat, R^{\flat+})} \times^{\cl{G}_k} \msf{G}^\sigma, \quad
        \varphi_P = \varphi_{\cl{P}}\vert_{\Spa(R^\flat, R^{\flat +})} \colon \varphi^*P \cong P, 
    \]
    i.e. the restriction along the locus $\{p = 0\}$. The set of depth-zero level structures on $\cl{P}$ is bijective to
    \[
        \{ 
            p \in P(R^\flat) , \quad \varphi_P(p) = p
        \}
    \]
    and the right action of $\msf{G}^\sigma(k)$ on each $p$ is given by $p \cdot g = pg$ for $g \in \msf{G}^\sigma(k)$. 
\end{lem}
\begin{proof}
    The set of $\cl{G}^+$-subtorsors of $\cl{P}$ is equivalent to the set of $\cl{G}^+$-subtorsors of
    \[
        (\cl{P}\vert_{\cl{Y}_{S, [0,n]}}) \otimes W_{O_F}(R^\flat)
    \]
    by \eqref{eq:YSn} and the Beauville-Laszlo lemma as $\cl{G}^+_F = \cl{G}_F$. By \Cref{lem:Gpsubtors}, it is bijective to $P(R^\flat)$. 

    A depth-zero level structure on $\cl{P}$ is equivalent to a $\cl{G}^+$-subtorsor $\cl{P}^+$ of $\cl{P}$ stable under $\varphi_{\cl{P}}$. The action of $\varphi_{\cl{P}}$ on the set of $\cl{G}^+$-subtorsors is given by $\varphi_P$ via the above bijection. Thus, we get the desired bijection. 

    For each $g \in \msf{G}^\sigma(k)$ with a lift $\wtd{g} \in \cl{G}(O_F)$, $g$ sends a $\cl{G}^+$-subtorsor $\cl{P}^+ \to \cl{P}$ to 
    \[
        \cl{P}^+\cdot \wtd{g} \to \cl{P}. 
    \]   
    Through the above bijection, this action is given as in the statement. 
    %equals the stated action on the set of $\varphi_P$-equivariant sections. 
\end{proof}

\section{Local Shimura varieties and universal deformations} \label{sec:locShim}

\subsection{Local Shimura data} \label{ssec:LSD}

In this section, we explain a choice of a length-zero representative $b \in [b]$ for each basic local Shimura datum $(G, [b], \mu)$. From now on, $G$ is unramified and $\cl{G}$ is a reductive model. Fix a Borel pair $T\subset B \subset G$ over $F$. Note that $T$ uniquely extends to an $O_F$-torus $\cl{T} \subset \cl{G}$. 

In this paper, italic, calligraphy, sans-serif types denote generic fibers over $F$, integral models over $O_F$, and special fibers over $\ov{k}$. Moreover, breve accents denote the base change to $O_{\breve{F}}$ or $\breve{F}$. For example, 
\[
    \breve{G} = G_{\breve{F}} , \quad \breve{\cl{G}} = \cl{G}_{O_{\breve{F}}}, \quad \msf{G} = \cl{G}_{\ov{k}}. 
\]

\subsubsection{Length-zero representatives} \label{sssec:length_zero}

Let $\cl{A}(\breve{G},\breve{T})$ be the reduced Bruhat-Tits apartment and let
\[
    o \in \mfr{a} \subset \cl{A}(\breve{G},\breve{T})
\]
be the hyperspecial point and the base alcove associated to $\cl{G}$ and $B$. Fix an identification
\[
    X_*(T/Z_G) \otimes \bb{R} \cong \cl{A}(\breve{G},\breve{T}), \quad 
    0 \mapsto o. 
\]
Then, $\mfr{a}$ lies in the codominant Weyl chamber with respect to $B$. 

Let $\breve{\cl{I}}$ be the Iwahori subgroup associated to $\mfr{a}$. Let
\[
    \wtd{W}=N_G(T)(\breve{F})/T(\breve{F}) \cap \breve{\cl{I}}
\]
be the Iwahori-Weyl group. There is a length function
$
    \ell \colon \wtd{W} \to \bb{Z}
$
and the set $\Omega \subset \wtd{W}$ of length-zero elements equals the stabilizer of $\mfr{a}$. Let $B(G)_\bas$ be the set of basic $\sigma$-conjugacy classes in $G(\breve{F})$. As in \cite{GHKR10} and \cite{He14}, there is a surjection
\[
    \Omega \twoheadrightarrow B(G)_\bas
\]
and each $[b]\in \Omega$ has a decent representative 
\[
    b = \mu(-\pi)w ,\quad w \in N_{\cl{G}}(\cl{T})(O_{\breve{F}})
\]
such that $\mu$ is dominant and minuscule, and $(w\sigma)^N = \sigma^N$ for some $N \geq 1$. Let $G_b$ be the extended pure inner form of $G$ associated to $b$. We fix an identification
\[
    G_b(F) = \{ g \in G(\breve{F}) \mid b\sigma(g)b^{-1} = g \} \subset G(\breve{F}). 
\]
\begin{rmk}
    The association between $[b] \in \Omega$ and minuscule dominant cocharacters $\mu \in X_*(\breve{T})$ is bijective. This is because for such $\mu$, we have 
    $
        \mu(-\pi)\cdot o \in \mfr{a}, 
    $
    so there is a unique $w$ (up to $\cl{T}(O_{\breve{F}})$) such that $w\cdot \mfr{a} = \mu(-\pi)^{-1}\cdot \mfr{a}$. Then, $[b]=\mu(-\pi)w \in \Omega$ is attached. 
\end{rmk}

As $[b] \in B(G, \mu)$, $(G, [b], \mu)$ is a basic local Shimura datum. From now on, we fix $(G, [b], \mu)$ and a length-zero representative $b \in [b]$. We adopt the opposite sign convention for $\mu$ compared to \cite{SW20}. 

\subsubsection{Local Shimura varieties}

For each compact open subgroup $K \subset G(F)$, a local Shimura variety $\cl{M}_{G, b, \mu, K}$ is defined as a rigid analytic variety over $\breve{F}$. 
When $K = \cl{G}(O_F)$, it is also denoted by $\cl{M}_{\cl{G}, b, \mu}$ and the associated moduli extends to perfectoid spaces over $O_{\breve{F}}$. 

\begin{defi}\textup{(\cite[Definition 25.1.1]{SW20})} \label{defi:v_sheaf_integral_model}
    The $v$-sheaf theoretic integral model $\cl{M}_{\cl{G},b,\mu}^\ints$ of $\cl{M}_{\cl{G}, b, \mu}$ is a $v$-sheaf sending a perfectoid space $S$ over $O_{\breve{F}}$ to the set of pairs $(\cl{P}, \iota)$ where 
    \begin{enumerate}
        \item $\cl{P}$ is a local $\cl{G}$-shtuka whose Frobenius $\varphi_{\cl{P}}$ is bounded by $-\mu$, and
        \item $\iota \colon \cl{P}\vert_{\cl{Y}_{[r, \infty)}} \cong \cl{E}^b\vert_{\cl{Y}_{[r, \infty)}}$ (for sufficiently large $r$) is a quasi-isogeny from $\cl{P}$ to $\cl{E}^b$. 
    \end{enumerate}
\end{defi}

\begin{rmk}
    In most cases of abelian type, it is known to be representable by a formal scheme over $O_{\breve{F}}$ (see \cite{PR26_local}), but much less is known in general, especially in exceptional cases. In a recent ongoing work of Lee and Madapusi \cite{Si24}, the representability at hyperspecial levels is announced in full generality. 
\end{rmk}

In general, $\cl{M}_{\cl{G}, b, \mu}^{\ints}$ is known to be a prekimberlite and the underlying space $(\cl{M}_{\cl{G}, b, \mu}^{\ints})_\red$ is equal to the affine Deligne-Lusztig variety $X_\mu(b)$ (see \cite[{}2.61, 2.63]{Gle22a}). Then, for each closed point $x \in X_\mu(b)$, we can define a formal completion $(\cl{M}^\ints_{\cl{G},b,\mu})^\wedge_{/x}$ by \cite[Definition 4.18]{Gle24}: it is a $v$-subsheaf of $\cl{M}^\ints_{\cl{G},b,\mu}$ consisting of geometric points whose rank $1$ generalizations map to $x$ under the specialization map. When $\cl{M}^\ints_{\cl{G}, b, \mu}$ is representable by a formal scheme $\cl{S}$, $(\cl{M}^\ints_{\cl{G},b,\mu})^\wedge_{/x}$ is represented by the formal completion of $\cl{S}$ at $x$. 

In particular, the representability of $(\cl{M}^\ints_{\cl{G},b,\mu})^\wedge_{/x}$ is milder than the whole representability, and it is in fact proved when $\cl{G}$ is reductive by \cite{Bar22} and \cite{Ito25a}. 

\begin{thm}\textup{(\cite[Theorem 4.4.2, Theorem 5.3.5]{Ito25b})} \label{thm:univdef}
    There exist a prismatic $(\cl{G}, \mu)$-display $\mfr{Q}^\univ$ over $(R_{\cl{G}, \mu})_{\Prism, O_F}$ and an isomorphism of $v$-sheaves 
    \[
        (\cl{M}^\ints_{\cl{G},b,\mu})^\wedge_{/x} \cong \Spd(R_{\cl{G}, \mu})
    \]
    such that the universal local $\cl{G}$-shtuka over $(\cl{M}^\ints_{\cl{G},b,\mu})^\wedge_{/x}$ is sent to the one induced by $\mfr{Q}^\univ$. 
\end{thm}

In the next section, we will briefly review the theory of prismatic displays to explain the full power of this theorem. Especially, we study a deformation at the distinguished closed point $[1]$ in 
\[
    (\cl{M}_{\cl{G}, b, \mu}^\ints)^\red \cong X_\mu(b) = \{ g \in G(\breve{F}) / \cl{G}(O_{\breve{F}}) \mid g^{-1} b \sigma(g) \in \cl{G}(O_{\breve{F}}) \mu(\pi) \cl{G}(O_{\breve{F}}) \}. 
\]

\subsection{Universal deformation at $[1]$} \label{ssec:univdef}

In this section, we explicitly describe the universal deformation $\mfr{Q}^\univ$ at $[1] \in X_\mu(b)$. Let $\Phi \subset X^*(\breve{T})$ be the set of roots of $\breve{G}$. For each $\alpha\in \Phi$, fix an identification
\[
    i_{\alpha} \colon \bb{A}^1 \cong \cl{U}_{\alpha} \subset \breve{\cl{G}}
\]
for the root group of $\alpha$ and let $u_\alpha$ be its coordinate. The special fiber of $\cl{U}_\alpha$ is denoted by $\msf{U}_\alpha$. Let $\Phi_{\mu<0} \subset \Phi$ be the set of roots $\alpha$ with $\langle \mu, \alpha \rangle = -1$ and let
\[
    \cl{U}_{\mu < 0} = \prod_{\alpha \in \Phi_{\mu<0}} \cl{U}_\alpha. 
\]
Since $\mu$ is minuscule, $\cl{U}_{\mu<0}$ is abelian. Then, we have
$
    R_{\cl{G},\mu} = O_{\breve{F}}\llbracket u_\alpha \mid \alpha \in \Phi_{\mu<0} \rrbracket. 
$

To provide an explicit description of $\mfr{Q}^\univ$, we adopt the equivalent formulation of prismatic $(\cl{G},\mu)$-displays given in \cite[Section 1.5.2]{IKY25}. 

\begin{defi}
    A prismatic $(\cl{G},\mu)$-display over an $O_F$-prism $(A,I)$ is a pair $(P, \varphi_{P})$ where $P$ is a $\cl{G}$-torsor over $A$ and
    \[
        \varphi_{P} \colon \varphi_A^*P[1/I] \cong P[1/I]
    \]
    is an isomorphism \'{e}tale locally described as 
    \[
        p \mapsto g_1 \mu(d) g_2 \varphi_A(p)
    \]
    for some $g_1,g_2 \in \cl{G}(A)$ and $I=(d)$.
    A prismatic $(\cl{G}, \mu)$-display over an $O_F$-algebra $R$ is a collection of prismatic $(\cl{G}, \mu)$-displays $P_{(A, I)}$ for each $(A, I) \in (R)_{\Prism, O_F}$ equipped with a functorial isomorphism
    \[
        P_{(A, I)} \otimes_A B \cong P_{(B, J)}
    \]
    for every map $(A, I) \to (B, J)$ in $(R)_{\Prism, O_F}$. 
\end{defi}

\begin{rmk}
    A prismatic $(\cl{G}, \mu)$-display over a perfect $O_F$-prism $(A, I)$ is equal to a $\cl{G}$-BKF module over a perfectoid $O_F$-algebra $A / I$ of type $\mu$.
\end{rmk}

Now, we may describe the universal deformation $\mfr{Q}^\univ$ over $R_{\cl{G},\mu}$ at $[1] \in X_\mu(b)$. Let
\[
    (\mfr{S},\cl{E}) = (O_{\breve{F}}\llbracket t, u_\alpha \mid \alpha \in \Phi_{\mu<0} \rrbracket, \pi-t) ,\quad \mfr{S}/\cl{E} \cong R_{\cl{G},\mu}
\]
be an (oriented) Breuil-Kisin $O_F$-prism associated with $R_{\cl{G},\mu}$.

\begin{prop}\textup{(\cite[Theorem 6.1.3]{Ito25b})}
    The restriction functor from the groupoid of prismatic $(\cl{G},\mu)$-displays over $R_{\cl{G},\mu}$ to the one over $(\mfr{S},\cl{E})$ is an equivalence. 
\end{prop}

\begin{prop}\textup{(\cite[Example 2.4.6, Theorem 4.4.2]{Ito25a})}
    The prismatic $(\cl{G}, \mu)$-display $\mfr{Q}^\univ$ over $R_{\cl{G}, \mu}$ such that
    \[
        \mfr{Q}^\univ(\mfr{S},\cl{E}) = \bigl(\cl{G} \otimes \mfr{S}, \mu(\cl{E}) \prod_{\alpha \in \Phi_{\mu<0}} i_\alpha(-u_\alpha) \mu(-1)w \sigma \bigr)
    \]
    is the universal deformation at a closed point $[1] \in X_\mu(b)$. 
\end{prop}
\begin{proof}
    It follows from the description in the proof of \cite[Theorem 4.4.2]{Ito25a}. Here, $\mu(-1)$ appears due to our choice $b = \mu(-\pi)w$. 
\end{proof}

\begin{cor} \label{lem:expBKF}
    Let $R$ be a perfectoid $R_{\cl{G},\mu}$-algebra. Any choice of elements $\pi^\flat, u_\alpha^\flat \in R^\flat$ such that $(\pi^\flat)^\sharp = \pi$ and $(u_\alpha^\flat)^\sharp=u_\alpha$ yields a map
    \[
        (\mfr{S},\cl{E})\to (W_{O_F}(R^\flat), \Ker(\theta_R)), \quad
        u_\alpha \mapsto [u_\alpha^\flat], \quad
        t \mapsto [\pi^\flat] 
    \]
    of $O_{F}$-prisms. Then, $\mfr{Q}^\univ(W_{O_F}(R^\flat), \Ker(\theta_R))$ is isomorphic to the $\cl{G}$-BKF module
    \[
        \bigl(\cl{G}\otimes W_{O_F}(R^\flat), \mu([\pi^\flat]-\pi) \prod_{\alpha \in \Phi_{\mu<0}} i_\alpha([u_\alpha^\flat])w\sigma\bigr). 
    \]
\end{cor}
\begin{proof}
    Let $\mfr{S}_n = O_{\breve{F}}\llbracket t^{1/p^n}, u_\alpha^{1/p^n} \mid \alpha \in \Phi_{\mu<0} \rrbracket$ for $n \geq 0$ and let $\mfr{S}_\infty = \colim_{n>0}^\wedge \mfr{S}_n$ be the $p$-completed colimit. Then, $(\mfr{S}_\infty, \cl{E})$ is the perfection of $(\mfr{S},\cl{E})$. For the first claim, we will construct $\mfr{S}_\infty/\cl{E} \to R$, which provides $(\mfr{S}_\infty, \cl{E}) \to (W_{O_F}(R^\flat), \Ker(\theta_R))$ by \cite[Theorem 3.10]{BS22}. We have $\mfr{S}_n/\cl{E} = R_{\cl{G},\mu}[\pi^{1/p^n}, u_\alpha^{1/p^n}]$, and the choice of $\pi^\flat$ and $u_\alpha^\flat$ induces 
    \[
        \mfr{S}_n/\cl{E} \to R, \quad \pi^{1/p^n} \mapsto (\pi^{\flat,1/p^n})^\sharp, \quad u_\alpha^{1/p^n} \mapsto (u_\alpha^{\flat,1/p^n})^\sharp. 
    \]
    These maps are compatible with transition maps for all $n \geq 1$ and induce $\mfr{S}_\infty/\cl{E} \to R$, so we get the first claim. The second claim follows from the description of $\mfr{Q}^\univ(\mfr{S},\cl{E})$. Note that since $\mu$ is minuscule, we have
    \[
        \mu(\pi-[\pi^\flat]) \prod_{\alpha \in \Phi_{\mu<0}} i_\alpha(-[u_\alpha^\flat])\mu(-1)w = \mu([\pi^\flat]-\pi) \prod_{\alpha \in \Phi_{\mu<0}} i_\alpha([u_\alpha^\flat])w.
    \]
\end{proof}

\begin{rmk}
    By Andr\'{e}'s flatness lemma (see \cite[Theorem 7.12]{BS22}), every perfectoid $R_{\cl{G},\mu}$-algebra $R$ admits a $p$-completely faithfully flat map $R\to S$ of perfectoid rings such that we can take $\pi^\flat$ and $u_\alpha^\flat$ in $S^\flat$. 
\end{rmk}

\begin{rmk} \label{rmk:changechoice}
    Let $(R,R^+)$ be a perfectoid Huber pair with a continuous homomorphism $R_{\cl{G},\mu} \to R^+$. For any choices $\pi^{\flat}, \pi^{\flat'}, u_\alpha^\flat,  u_\alpha^{\flat'} \in R^{\flat+}$ with $(\pi^{\flat})^\sharp = (\pi^{\flat'})^\sharp = \pi$ and $(u_\alpha^{\flat})^\sharp = (u_\alpha^{\flat'})^\sharp=u_\alpha$, the functoriality of $\mfr{Q}^\univ$ provides an element $\iota \in \cl{G}(W_{O_F}(R^{\flat+}))$ such that
    \begin{equation}
        \iota \mu([\pi^\flat]-\pi) \prod_{\alpha \in \Phi_{\mu<0}} i_\alpha([u_\alpha^\flat])w \sigma(\iota)^{-1} = \mu([\pi^{\flat'}]-\pi) \prod_{\alpha \in \Phi_{\mu<0}} i_\alpha([u_\alpha^{\flat'}])w. \label{eq:choice}
    \end{equation}
    Let $\varpi \in R^{\flat+}$ be a pseudo-uniformizer that divides $\pi^{\flat}$, $\pi^{\flat'}$, $u_\alpha^\flat$ and $u_\alpha^{\flat'}$. The compositions
    \[
        (\mfr{S},\cl{E}) \to (W_{O_F}(R^{\flat +}), \Ker(\theta_{R^+})) \to  (W_{O_F}(R^{\flat+})/[\varpi], \Ker(\theta_{R^+})/[\varpi])
    \]
    for two choices $(\pi^{\flat}, u_{\alpha}^\flat)$ and $(\pi^{\flat'}, u_{\alpha}^{\flat'})$ are equal, so $\iota$ is trivial modulo $[\varpi]$. By the standard argument as in \Cref{lem:uniqext}, $\iota$ can be characterized as a unique element satisfying \eqref{eq:choice} that is trivial modulo $[\varpi]$ for some pseudo-uniformizer $\varpi \in R^{\flat+}$. 
\end{rmk}

Then, we can provide an explicit description of depth-zero level structures on $\mfr{Q}^\univ$ in characteristic $0$. Here, we choose $\msf{G}^\sigma = \cl{G}_k$ as in \Cref{rmk:choice_G+}. 

\begin{prop} \label{cor:explevel}
    Let $x\colon \Spa(R,R^+) \to \Spf(R_{\cl{G},\mu})_\eta$ be a map from an affinoid perfectoid space over $\breve{F}$. Let $\pi^\flat, u_\alpha^\flat \in R^{\flat+}$ be elements such that $(\pi^\flat)^\sharp = \pi$ and $(u_\alpha^\flat)^\sharp=u_\alpha$. A depth-zero level structure on $\mfr{Q}^\univ$ at $x$ is given by an element $g\in \msf{G}(R^\flat)$ such that 
    \begin{equation}
        g\sigma(g)^{-1} = \mu(\pi^\flat) \prod_{\alpha \in \Phi_{\mu<0}} i_\alpha(u_\alpha^\flat)w. \label{eq:depth0sol}
    \end{equation}
    The right action of $g_0 \in \cl{G}(k)$ is given by $g \cdot g_0 = gg_0$.  
    %When $(R,R^+)$ is an algebraically closed perfectoid field, the specialization of $x$ in $\mfr{X}_b^{\msf{\ov{G}}}$ is identified with the specialization of $g$ in $\msf{\ov{G}}$ via the inclusion $\mfr{X}_{b,\red}^{\msf{\ov{G}},\perf} \hookrightarrow \ov{\msf{G}}^\perf$. 
\end{prop}
\begin{proof}
    The claim follows from the explicit description in \Cref{lem:expBKF} and \Cref{lem:finetdesc}. %The second claim follows from \Cref{rmk:computespc}. 
\end{proof}
\begin{rmk} \label{rmk:changelevel}
    For any different choices $\pi^{\flat'}, u_{\alpha}^{\flat'} \in R^{\flat+}$ as in \Cref{rmk:changechoice}, the isomorphism \eqref{eq:choice} sends a depth-zero level structure $g$ to $\iota g$. %Since $\iota$ is trivial modulo $[\varpi]$ for some pseudo-uniformizer $\varpi \in R^{\flat+}$, the specialization of $g$ is equal to that of $\iota g$ in $\msf{\ov{G}}$. 
\end{rmk}

The following lemma is useful in constructing a solution to \eqref{eq:depth0sol}. 
%A special affinoid is constructed from a certain simple recipe of the above equation. An important lemma is the following. 

\begin{lem} \label{lem:uniqext}
    Let $(R,R^+)$ be a perfectoid Huber pair in characteristic $p$ and let $\varpi \in R^+$ be a pseudo-uniformizer. Let $g\in \msf{G}(R^{+})$. For every element $h_0\in \msf{G}(R^+/\varpi)$ such that $h_0 \sigma(h_0)^{-1} \equiv g$, there is a unique lift $h \in \msf{G}(R^+)$ of $h_0$ such that $h \sigma(h)^{-1} = g$. 
\end{lem}
\begin{proof}
    We construct a sequence of elements $h_n \in \msf{G}(R^+/\varpi^{n+1})$ such that $h_n\sigma(h_n)^{-1} \equiv g$ and $h_n \equiv h_{n+1}$ for $n \geq 0$. Suppose that we have $h_n$ for some $n\geq 0$. Since $\msf{G}$ is smooth, we can take a lift $h'_{n+1} \in \msf{G}(R^+/\varpi^{n+2})$ of $h_n$. Then, $d=h'_{n+1}\sigma(h'_{n+1})^{-1}g^{-1}$ is trivial modulo $\varpi^{n+1}$. Since $\sigma(d)$ is trivial modulo $\varpi^{n+2}$ as $q \geq 2$, $h_{n+1} = d^{-1} h'_{n+1}$ is a unique lift satisfying $h_{n+1}\sigma(h_{n+1})^{-1} \equiv g$. Then, the limit $h=\lim_{n\geq 0} h_n$ satisfies $h\sigma(h)^{-1} = g$ and it is the unique lift of $h_0$ satisfying that condition. 
\end{proof}

\subsubsection{The inner action}

Recall that $G_b(F)$ acts on $X_\mu(b)$ so that $j \cdot [g] = [jg]$ for $j \in G_b(F)$ and $[g] \in X_\mu(b)$. In particular, the stabilizer of $[1]\in X_\mu(b)$ is $G_b(F) \cap \cl{G}(O_{\breve{F}})$, and it acts on $\Spf(R_{\cl{G},\mu})$ by \Cref{thm:univdef}. Here, we provide an explicit description of this inner action. 

\begin{prop} \label{prop:innactuniv}
    Let $(R,R^+)$ be a perfectoid Huber pair with a continuous homomorphism $R_{\cl{G},\mu} \to R^+$. Let $\pi^\flat, u_\alpha^\flat \in R^\flat$ be elements such that $(\pi^\flat)^\sharp = \pi$ and $(u_\alpha^\flat)^\sharp=u_\alpha$. For each $j\in G_b(F) \cap \cl{G}(O_{\breve{F}})$, the composition $R_{\cl{G},\mu} \xrightarrow{j} R_{\cl{G},\mu} \to R^+$ corresponds to a deformation
    \[
        \bigl(\cl{G}\otimes W_{O_F}(R^{\flat+}), j\mu([\pi^\flat]-\pi) \prod_{\alpha \in \Phi_{\mu<0}} i_\alpha([u_\alpha^\flat])w\sigma(j)^{-1} \sigma\bigr). 
    \]
    of the prismatic $(\cl{G}, \mu)$-display $(\breve{\cl{G}}, b\sigma)$. 
\end{prop}
\begin{proof}
    Let $\cl{P} = (\cl{G}\otimes W_{O_F}(R^{\flat +}), \mu([\pi^\flat]-\pi) \prod_{\alpha \in \Phi_{\mu<0}} i_\alpha([u_\alpha^\flat])w\sigma)$ be a deformation of $(\breve{\cl{G}}, b\sigma)$ over $R^+$. The quasi-isogeny of local $\cl{G}$-shtukas between $\cl{P}$ and $(\breve{\cl{G}}, b\sigma)$ parametrized by $\Spd(R^+) \to \Spd(R_{\cl{G},\mu}) \to \cl{M}^\ints_{\cl{G},b,\mu}$ is given by the identification
    \[
        (\cl{G}\otimes W_{O_F}(R^{\flat+})/[\varpi], \mu([\pi^\flat]-\pi) \prod_{\alpha \in \Phi_{\mu<0}} i_\alpha([u_\alpha^\flat])w\sigma) = (\cl{G}\otimes W_{O_F}(R^{\flat+})/[\varpi], b\sigma),
    \]
    where $\varpi$ is taken as in \Cref{rmk:changechoice}, via \cite[Proposition 5.2.6]{Ito25a} (see also \cite{Gle22a}). The action of $j$ sends the quasi-isogeny to the one given by the left translation
    \[
        (\cl{G}\otimes W_{O_F}(R^{\flat+})/[\varpi], \mu([\pi^\flat]-\pi) \prod_{\alpha \in \Phi_{\mu<0}} i_\alpha([u_\alpha^\flat])w\sigma) \xrightarrow{j} (\cl{G}\otimes W_{O_F}(R^{\flat+})/[\varpi], b\sigma). 
    \]
    Since $j\in \cl{G}(O_{\breve{F}})$, this left translation lifts to an isomorphism
    \[
        \cl{P} \xrightarrow{j} (\cl{G}\otimes W_{O_F}(R^{\flat+}), j\mu([\pi^\flat]-\pi) \prod_{\alpha \in \Phi_{\mu<0}} i_\alpha([u_\alpha^\flat])w\sigma(j)^{-1} \sigma). 
    \]
    Thus, the claim follows since the quasi-isogeny is sent to the identification
    \[
        (\cl{G}\otimes W_{O_F}(R^{\flat+})/[\varpi], j\mu([\pi^\flat]-\pi) \prod_{\alpha \in \Phi_{\mu<0}} i_\alpha([u_\alpha^\flat])w\sigma(j)^{-1}\sigma) = (\cl{G}\otimes W_{O_F}(R^{\flat+})/[\varpi], b\sigma)
    \]
    through the above isomorphism. 
\end{proof}

\section{Special affinoids at depth zero} \label{sec:DL}

%In this section, we first specify the parabolic Deligne-Lusztig variety $\msf{W}(0)$ and then introduce the special affinoid $\cl{W}(0)$. The main goal is to prove \Cref{thm2:prop_specialaff}.  

\subsection{Parabolic Deligne-Lusztig varieties}

In this section, we introduce a parabolic Deligne-Lusztig variety $\msf{W}(0)$ that will be shown later to be the reduction of the special affinoid $\cl{W}(0)$. An important role is played by the following rational cocharacter $\lambda$. 

\begin{lem} \label{lem:lambdadef}
    There is a unique element $\lambda \in X_*(\msf{T})_{\bb{Q}}$ such that $\lambda - qw\sigma \lambda = \mu$.
\end{lem}
\begin{proof}
    There is $N\geq 1$ such that $(w\sigma)^N$ acts trivially on $X_*(\msf{T})$. Since 
    \[
        (q^N-1)\lambda = (qw\sigma)^N \lambda - \lambda = -\sum_{0 \leq k < N} (qw\sigma)^k \mu, 
    \]
    we have $\lambda = -\tfrac{1}{q^N-1}\sum_{0 \leq k < N} (qw\sigma)^k \mu$. It satisfies the equation $\lambda - qw\sigma \lambda = \mu$.
\end{proof}

\begin{lem} \label{lem:lambdapst}
    For every $\alpha \in \Phi$, $\langle \alpha, \lambda \rangle >0$ if and only if there is a positive integer $k$ such that 
    \[
        \langle \alpha, (w\sigma)^{-k}\mu \rangle < 0 , \quad
        \langle \alpha, (w\sigma)^{-i}\mu \rangle=0 \quad (0< i < k). 
    \]
    In particular, $\langle \alpha, \lambda \rangle =0$ if and only if $\langle \alpha, (w\sigma)^{-i}\mu \rangle=0$ for every $i>0$. 
\end{lem}
\begin{proof}
    From the description of $\lambda$, it is easy to see that $\langle \alpha, \lambda \rangle =0$ if $\langle \alpha, (w\sigma)^{-i}\mu \rangle=0$ for every $i>0$. Let $k$ be the minimum positive integer such that $\langle \alpha, (w\sigma)^{-k}\mu \rangle \neq 0$. Since $\mu$ is minuscule, $\lvert  \langle \alpha, (qw\sigma)^{-i} \mu \rangle \rvert \leq q^{-i}$ for every $i>0$. Thus, $\lvert \langle \alpha, \sum_{i>k} (qw\sigma)^{-i}\mu \rangle \rvert \leq \tfrac{q^{-k}}{q-1}$ and the equality holds if and only if $\langle \alpha, (w\sigma)^{-i} \mu \rangle$ has the same sign for every $i>0$. Since $\langle \alpha, (w\sigma)^{-k}\mu \rangle=\pm 1$ and we have $\lambda = - \sum_{i>0} (qw\sigma)^{-i}\mu$, the sign of $\langle \alpha, \lambda \rangle$ is opposite to that of $\langle \alpha, (w\sigma)^{-k}\mu \rangle$. 
\end{proof}

\begin{prop} \label{lem:lambdadom}
    The rational cocharacter $\lambda$ is dominant. 
\end{prop}
\begin{proof}
    Let $\phi_\alpha$ be the supremum of $\alpha$ on $\mfr{a}$ 
    for each $\alpha \in \Phi$. We have 
    \[
        \phi_\alpha = 0 \quad ( \alpha > 0) ,\quad 
        \phi_\alpha = 1 \quad ( \alpha < 0). 
    \]
    Since $b\sigma \mfr{a}=\mfr{a}$, we have 
    \[
        \phi_\alpha = \phi_{(w\sigma)^{-1}\alpha} - \langle \alpha, \mu \rangle. 
    \]
    For each $\alpha\in \Phi$ with $\langle \alpha, \lambda\rangle \neq 0$, let $k>0$ be the minimum positive integer such that $\langle \alpha, (w\sigma)^{-k} \mu \rangle \neq 0$. Then, we have $\phi_\alpha = \phi_{(w\sigma)^{k}\alpha} + \langle \alpha, (w\sigma)^{-k}\mu \rangle$, so 
    \begin{equation}
        \text{
        $\langle \alpha, (w\sigma)^{-k}\mu \rangle=1$ if and only if $\alpha<0$ and $(w\sigma)^{k}\alpha>0$.
        } 
        \label{eq:lambdapst}
    \end{equation}
    Thus, if $\langle \alpha, \lambda \rangle >0$, then $\alpha>0$ by \Cref{lem:lambdapst}, so $\lambda$ is dominant.
\end{proof}

We will see later that the denominator of $\lambda$ measures the tame ramification of the associated Galois representation. 

\begin{lem} \label{lem:tame_ramification}
    The minimal positive integer $e \geq 1$ such that $e\lambda \in X_*(\msf{T})$ is coprime to $p$. 
\end{lem}
\begin{proof}
    The claim follows since $e$ divides $q^{N}-1$ for some $N\geq 1$. 
\end{proof}
%We will introduce several closed subgroups of $\msf{G}$ that are needed in our discussion.
\begin{defi}
    Let
    \[
        \msf{P} = \msf{G}^{e\lambda \geq 0} ,\quad 
        \ov{\msf{P}} = \msf{G}^{e\lambda \leq 0} ,\quad 
        \msf{M} = \msf{G}^{e\lambda = 0} ,\quad
        \msf{N} = \msf{G}^{e\lambda > 0} ,\quad
        \ov{\msf{N}} = \msf{G}^{e\lambda < 0}
    \]
    and let
    \[
        \msf{U}_{\mu>0} = \msf{G}^{\mu > 0} ,\quad
        \msf{U}_{\mu<0} = \msf{G}^{\mu < 0}. 
    \]
    Moreover, let
    \[
        \Ad(w\sigma) \colon \msf{G} \to \msf{G} ,\quad g \mapsto w\sigma(g)w^{-1}
    \]
    be a Frobenius-twisted endomorphism. 
\end{defi}

\begin{prop} \label{lem:MNprop}
    $\msf{M}$ is stable under $\Ad(w\sigma)$ and $\msf{N} \cap \Ad(w\sigma)(\msf{\ov{N}}) = \msf{U}_{\mu>0}$. 
\end{prop}
\begin{proof}
    The first claim follows from \Cref{lem:lambdapst}. For the second claim, $\msf{N} \cap \Ad(w\sigma)(\msf{\ov{N}}) \subset \msf{U}_{\mu>0}$ since $\langle (w\sigma)^{-1}\alpha, \lambda \rangle <0$ and $\langle \alpha, \lambda \rangle >0$ together imply $\langle \alpha, \mu \rangle =1$. On the other hand, $\langle \alpha, \mu \rangle =1$ implies  $\langle (w\sigma)^{-1}\alpha,  (w\sigma)^{-1}\mu \rangle=1$, so $(w\sigma)^{-1}\alpha<0$ and $\alpha > 0$ by \eqref{eq:lambdapst}. 
    Then, $\alpha$ is a root inside $\msf{N} \cap \Ad(w\sigma)(\msf{\ov{N}})$ by \Cref{lem:lambdadom} since $\langle \alpha, \lambda \rangle \neq 0$. 
\end{proof}

\begin{defi}
    Let
    \[
        \msf{W}(0) =  \{g \in \msf{G}/\msf{N}\mid g^{-1} \sigma(g) \in \msf{N} w \sigma(\msf{N}) \}
    \]
    be the parabolic Deligne-Lusztig variety attached to $(\msf{P}, w)$. It is equipped with an action of $\msf{G}^\sigma \times \msf{M}^{w\sigma}$ such that $(g_0,m) \cdot g = g_0 g m^{-1}$. 
\end{defi}

Here, we identify $w$ with its reduction in $\msf{G}(\ov{k})$. Now, we will give another description as a Lang-torsor over $\msf{U}_{\mu<0}$. First, we have
\begin{equation}
    \msf{W}(0) = \{g \in \msf{G} \mid g^{-1}\sigma(g)w^{-1} \in \Ad(w\sigma)(\msf{N}) \}/(\msf{N} \cap \Ad(w\sigma)(\msf{N})). 
    \label{eq:firstrewrite}
\end{equation} 
By \Cref{lem:MNprop}, $\msf{\ov{N}}\cap \Ad(w\sigma)(\msf{N})$ is normal in $\Ad(w\sigma)(\msf{N})$. Let 
\begin{equation} \label{eq:definition_piw}
    \pi_w \colon \Ad(w\sigma)(\msf{N}) \to \msf{N} \cap \Ad(w\sigma)(\msf{N})
\end{equation}
be the quotient by $\msf{\ov{N}} \cap \Ad(w\sigma)(\msf{N})$ and let 
\[
    \phi_w \colon \msf{N} \cap \Ad(w\sigma)(\msf{N}) \to \msf{N} \cap \Ad(w\sigma)(\msf{N}) ,\quad
    h \mapsto \pi_w(\Ad(w\sigma)(h))
\]
be a Frobenius-twisted endomorphism. 
\begin{lem} \label{lem:phiwnilp}
    There is a positive integer $N>0$ such that $\phi_w^N$ is trivial. 
\end{lem}
\begin{proof}
    Let $\Phi_0 \subset \Phi$ be the set of roots in $\msf{N} \cap \Ad(w\sigma)(\msf{N})$. Let $h\in  \msf{N} \cap \Ad(w\sigma)(\msf{N})$ be an element with a decomposition $h=\prod_{\alpha \in J} h_\alpha$ for $h_\alpha \in \msf{U}_\alpha$ and an ordered subset $J\subset \Phi_0$. Then, we have
    \[
        \phi_w(h)=\prod_{\alpha \in J \cap (w\sigma)^{-1}\Phi_0} \Ad(w\sigma)(h_\alpha). 
    \]
    Take $N>0$ so that $(w\sigma)^N$ acts trivially on $X^*(\msf{T})$ and let $J_N = J\cap \bigcap_{1\leq i \leq N} (w\sigma)^{-i}\Phi_0$. By the above computation, it follows that $\phi_w^N(h) = \prod_{\alpha \in J_N} \Ad(w\sigma)^N(h_\alpha)$. However, since the pairing of each element of $\Phi_0$ with $\mu$ is zero by \Cref{lem:MNprop}, it follows from \Cref{lem:lambdapst} that the pairing of each element of $J_N$ with $\lambda$ is zero and $J_N$ is empty as $J_N \subset \Phi_0$. In particular, $\phi_w^N$ is trivial. 
\end{proof}

\begin{cor} \label{lem:phiwtors}
    The map
    \[
        \msf{N} \cap \Ad(w\sigma)(\msf{N}) \to \msf{N} \cap \Ad(w\sigma)(\msf{N}) ,\quad
        h \mapsto h \phi_w(h)^{-1}
    \]
    is an isomorphism. 
\end{cor}
\begin{proof}
    As the differential of $\Ad(w\sigma)$ vanishes at every point, the proof of Lang's theorem (see \cite[Theorem 16.3]{Bor91}) works to see that the given map is a finite \'{e}tale torsor under the equalizer of $\phi_w$ and the identity map. Thus, the claim follows from \Cref{lem:phiwnilp}. 
\end{proof}

\begin{prop} \label{cor:W0desc}
    There is a $\msf{G}^\sigma$-equivariant isomorphism
    \[
        \msf{W}(0) \cong \{ g \in \msf{G} \mid g^{-1} \sigma(g)w^{-1} \in \msf{U}_{\mu<0}\}
    \]
    to a Lang torsor over $\msf{U}_{\mu < 0}$. In particular, $\msf{W}(0) \to \msf{U}_{\mu<0}$ is a finite \'{e}tale torsor under $\msf{G}^\sigma$, and $\msf{W}(0)$ is a smooth affine $\ov{k}$-variety of dimension $r$. 
\end{prop}
\begin{proof}
    By \Cref{lem:phiwtors}, for every $g\in \msf{G}$ with $g^{-1}\sigma(g) w^{-1} \in \Ad(w\sigma)(\msf{N})$, there is a unique $h \in \msf{N}\cap \Ad(w\sigma)(\msf{N})$ such that
    \[
        \pi_w(g^{-1}\sigma(g)w^{-1}) = h\phi_w(h)^{-1} \Leftrightarrow \pi_w((gh)^{-1}\sigma(gh)w^{-1})=1. 
    \]
    %It is equivalent to $\pi_w((gh)^{-1}\sigma(gh)w^{-1})=1$. 
    By \eqref{eq:firstrewrite}, every element of $\msf{W}(0)$ has a unique representative $g\in \msf{G}$ such that $g^{-1}\sigma(g) w^{-1} \in \msf{\ov{N}}\cap \Ad(w\sigma)(\msf{N})$. Thus, the claim follows from \Cref{lem:MNprop}. 
\end{proof}

\subsection{Construction of special affinoids}

In this section, we specify the radii $r_\alpha$ for each $\alpha \in \Phi_{\mu < 0}$ and define the special affinoid $\cl{W}(0)$. 

\begin{lem} \label{lem:r_alpha_defi}
    Let $r_\alpha = \langle \alpha, qw\sigma\lambda \rangle$ for $\alpha \in \Phi$. For every $\alpha \in \Phi_{\mu < 0}$, $r_\alpha \in \bb{Z}[\tfrac{1}{e}]$ and $r_\alpha > 0$. 
\end{lem}
\begin{proof}
    It follows from the definition of $e$ and \Cref{lem:MNprop}. 
\end{proof}

We record a more precise property for our future reference. 

\begin{lem} \label{lem:vareps}
    For each $\alpha \in \Phi_{\mu<0}$, there is a unique positive integer $n_\alpha \geq 1$ such that 
    \[
        \langle (w\sigma)^{-n_\alpha}\alpha, \mu \rangle = 1, \quad 
        \langle (w\sigma)^{- i}\alpha, \mu \rangle = 0 \quad (0 < i < n_\alpha). 
    \]
    Let $\beta_\alpha = - (w\sigma)^{-n_\alpha}\alpha$. The map $\alpha \mapsto \beta_\alpha$ is a permutation of $\Phi_{\mu < 0}$ and we have $r_\alpha = q^{n_\alpha}(1-r_{\beta_\alpha})$ for every $\alpha \in \Phi_{\mu<0}$. In particular, $1 - q^{-1} <r_\alpha<1$ for every $\alpha \in \Phi_{\mu<0}$. 
\end{lem}
\begin{proof}
    Since $\msf{U}_{\mu < 0} \subset \msf{\ov{N}}$ by \Cref{lem:MNprop}, $n_\alpha$ is well-defined by \Cref{lem:lambdapst} and it is easy to see that the map $\alpha \mapsto \beta_\alpha$ is injective. It is enough to prove the claims for $r_\alpha$. Since $\lambda - qw\sigma \lambda = \mu$, 
    \[
        r_{(w\sigma)^{-1}\alpha} = \langle (w\sigma)^{-1}\alpha, \lambda - \mu \rangle = q^{-1}r_\alpha - \langle (w\sigma)^{-1}\alpha, \mu \rangle. 
    \]
    Since $\langle (w\sigma)^{-i}\alpha, \mu \rangle = 0$ for $0 < i < n_\alpha$, we have
    \[
        q^{-n_\alpha} r_\alpha = r_{-\beta_\alpha} + \langle -\beta_\alpha, \mu \rangle = 1 -r_{\beta_\alpha}. 
    \]
    For the last claim, take $\gamma \in \Phi_{\mu<0}$ so that $\beta_\gamma = \alpha$. Since $r_\gamma > 0$ and $r_\alpha = 1 - q^{-n_\gamma} r_\gamma$, we have $r_\alpha < 1$. Then, we also have $r_\gamma < 1$, so $r_\alpha > 1 - q^{-1}$. 
\end{proof}

\begin{defi}
    Let
    \[
        \cl{U}(0) = \{ \lvert u_\alpha \rvert \leq \lvert \pi \rvert^{r_\alpha} \neq 0 \mid \alpha \in \Phi_{\mu < 0} \} \subset \Spf(R_{\cl{G}, \mu})_\eta
    \]
    be a rational subdomain and let $\cl{W}(0)$ be the $\msf{G}^\sigma$-cover of $\cl{U}(0)$ defined by the Cartesian diagram
    \begin{center}
        \begin{tikzcd}
            \cl{W}(0) \ar[r, hook] \ar[d] & \cl{M}_{G, b, \mu, \cl{G}(1)} \ar[d] \\
            \cl{U}(0) \ar[r, hook] & \cl{M}_{\cl{G}, b, \mu}. 
        \end{tikzcd}
    \end{center}
\end{defi}

\subsection{Reduction of special affinoids} \label{ssec:specaff}

In this section, we identify the reduction of $\cl{W}(0)$ with $\msf{W}(0)$. For this, we proceed by constructing a smooth affine model over a suitable tame extension. 

\begin{defi}
    Let $\breve{F}_e = \breve{F}[\pi^{1/e}]$ be the tame extension of $\breve{F}$ of degree $e$ and let $O_{\breve{F}_e} = O_{\breve{F}}[\pi^{1/e}]$ be the ring of integers of $\breve{F}_e$. Let
    \[
        R_{\cl{G},\mu, \lambda} = O_{\breve{F}_e}\langle \pi^{-r_\alpha}u_\alpha \mid \alpha \in \Phi_{\mu<0}\rangle
    \]
    be a $p$-adic $R_{\cl{G}, \mu}$-algebra. Here, $\pi^{d}$ denotes $(\pi^{1/e})^{ed}$ for every $d \in \bb{Z}[\tfrac{1}{e}]$. 
\end{defi}

In particular, $\Spf(R_{\cl{G}, \mu, \lambda})$ is a smooth affine model of $\cl{U}(0)$ over $O_{\breve{F}_e}$. We construct a smooth affine model of $\cl{W}(0)$ as a finite \'{e}tale cover of $\Spf(R_{\cl{G}, \mu, \lambda})$. 

\begin{const}
    Let $u_{\alpha, 1} = \pi^{-r_\alpha} u_\alpha \in R_{\cl{G}, \mu, \lambda}$ and fix an identification
    \[
        \Spec(R_{\cl{G}, \mu, \lambda}/\pi^{1/e}) \cong \msf{U}_{\mu < 0} ,\quad
        u_{\alpha, 1} \mapsto u_\alpha. 
    \]
    By the topological invariance of finite \'{e}tale sites, the finite \'{e}tale $\msf{G}^\sigma$-cover
    \[
        \msf{W}(0) \to \msf{U}_{\mu < 0} \cong \Spec(R_{\cl{G}, \mu, \lambda}/\pi^{1/e})
    \]
    uniquely lifts to a finite \'{e}tale $\msf{G}^\sigma$-cover
    \[
        \mfr{W}(0) \to \Spf(R_{\cl{G},\mu, \lambda}). 
    \]
    In particular, $\mfr{W}(0)$ is a smooth affine formal scheme over $O_{\breve{F}_e}$. 
\end{const}

\begin{thm} \label{thm:reduction_of_specaff}
    There is a $\msf{G}^\sigma$-equivariant isomorphism 
    \[
        \mfr{W}(0)_\eta \cong \cl{W}(0) \otimes \breve{F}_e
    \]
    over $\Spf(R_{\cl{G},\mu, \lambda})_\eta \cong \cl{U}(0) \otimes \breve{F}_e$. 
\end{thm}
\begin{proof}
    Since both sides are finite \'{e}tale $\msf{G}^\sigma$-covers of $\Spf(R_{\cl{G},\mu, \lambda})_\eta$, it is enough to construct a $\msf{G}^\sigma$-equivariant map
    \[
        \mfr{W}(0)_\eta^\diamond \to \cl{W}(0)^\diamond
    \]
    over $\cl{U}(0)^\diamond$ by \cite[Lemma 15.6]{Sch17}. 

    Let $(R,R^+)$ be a perfectoid Huber pair over $\mfr{W}(0)_\eta$. After passing to a $v$-cover (e.g. a product of geometric points), we may take $\pi^\flat$ and $u_\alpha^\flat$ as in \Cref{lem:expBKF}. Then, the evaluation $\mfr{Q}^\univ(R^+)$ is given by the $\cl{G}$-BKF module
    \[
        P_{R^+}=(\cl{G}\otimes W_{O_F}(R^{\flat+}), \mu([\pi^\flat]-\pi)\prod_{\alpha\in \Phi_{\mu<0}} i_\alpha([u_\alpha^\flat])w\sigma). 
    \]
    Since $e$ is coprime to $p$, $\pi^{1/e} \in R^+$ determines a compatible choice of $\pi^{\flat,1/e} \in R^{\flat+}$. We will construct a solution to \eqref{eq:depth0sol} using the map $\Spf(R^+) \to \mfr{W}(0)$. Since $\mu=\lambda - qw\sigma\lambda$, we have
    \[
        \mu(\pi^\flat)\prod_{\alpha\in \Phi_{\mu<0}} i_\alpha(u_\alpha^\flat)w
        = (\pi^{\flat,1/e})^{e\lambda} \cdot \prod_{\alpha\in \Phi_{\mu<0}} i_\alpha(\pi^{\flat,-r_\alpha} u_\alpha^\flat)w \cdot \sigma((\pi^{\flat,1/e})^{-e\lambda}). 
    \]
    Here, we use
    \[
        (\pi^{\flat,1/e})^{e\lambda} \cdot \Ad(w\sigma)((\pi^{\flat,1/e})^{-e\lambda}) = (\pi^{\flat,1/e})^{e\lambda} \cdot (\pi^{\flat,q/e})^{-ew\sigma \lambda} = (\pi^{\flat,1/e})^{e(\lambda - qw\sigma\lambda)} = \mu(\pi^\flat). 
    \]
    Now, $\pi^{\flat,-r_\alpha}u_\alpha^\flat \in R^{\flat+}$ since $\pi^{-r_\alpha}u_\alpha \in R^+$. By \Cref{cor:W0desc}, $\Spec(R^+ / \pi^{1/e}) \to \msf{W}(0)$ corresponds to an element $g_0 \in \msf{G}(R^+/\pi^{1/e})$ such that
    \[
        g_0^{-1} \sigma(g_0) \equiv \prod_{\alpha\in \Phi_{\mu<0}} i_\alpha(u_{\alpha, 1})w  \pmod{\pi^{1/e}}
        = 
        \prod_{\alpha\in \Phi_{\mu<0}} i_\alpha(\pi^{\flat,-r_\alpha} u_\alpha^\flat)w \pmod{\pi^{\flat, 1/e}}. 
    \]
    Then, by \Cref{lem:uniqext}, $g_0$ can be uniquely lifted to an element $g\in \msf{G}(R^{\flat +})$ such that 
    \[
        g^{-1}\sigma(g) = \prod_{\alpha\in \Phi_{\mu<0}} i_\alpha(\pi^{\flat,-r_\alpha} u_\alpha^\flat)w. 
    \]
    Thus, we get a solution $(\pi^{\flat,1/e})^{e\lambda}g^{-1}$ to \eqref{eq:depth0sol}. We will show that this depth-zero level structure is independent of the choice of $\pi^\flat$ and $u_\alpha^\flat$. 

    Let $\pi^{\flat'}$ (resp.\ $u^{\flat'}_\alpha$) be another choice of $\pi^\flat$ (resp.\ $u_\alpha^\flat$). Let
    \[
        P'_{R^+}=(\cl{G}\otimes W_{O_F}(R^{\flat+}), \mu([\pi^{\flat'}]-\pi)\prod_{\alpha\in \Phi_{\mu<0}} i_\alpha([u_\alpha^{\flat'}])w\sigma). 
    \]
    As in \Cref{rmk:changechoice}, there is a unique $\iota \in \cl{G}(W_{O_F}(R^{\flat+}))$ representing an isomorphism of $P_{R^+}$ and $P'_{R^+}$ as a deformation of $(\cl{G},b\sigma)$. In particular, $\iota$ is trivial modulo $[\pi^{\flat,1/p^m}]$ for some $m\geq 0$. Let $g'\in \msf{G}(R^{\flat+})$ be the lift of $g_0$ such that
    \[
        g'^{-1}\sigma(g') = \prod_{\alpha\in \Phi_{\mu<0}} i_\alpha(\pi^{\flat',-r_\alpha} u_\alpha^{\flat'})w. 
    \]
    Then, we have two depth-zero level structures on $P'_{R^+}$: 
    \[
        (\ov{\iota} (\pi^{\flat,1/e})^{e\lambda}g^{-1})\sigma(\ov{\iota} (\pi^{\flat,1/e})^{e\lambda}g^{-1})^{-1} = ((\pi^{\flat',1/e})^{e\lambda}g'^{-1})\sigma((\pi^{\flat',1/e})^{e\lambda}g'^{-1})^{-1}. 
    \]
    Here, $\ov{\iota}$ denotes the reduction of $\iota$ in $\msf{G}(R^{\flat+})$. By restricting to a closed and open subset of $\Spa(R,R^+)$, we can take $h\in \msf{G}^\sigma$ so that
    \[
        \ov{\iota} (\pi^{\flat,1/e})^{e\lambda}g^{-1} = (\pi^{\flat',1/e})^{e\lambda}g'^{-1}h.
    \]
    Since $\ov{\iota}$ is trivial modulo $\pi^{\flat, 1/p^m}$ and $\pi^{\flat, 1/e} \cdot \pi^{\flat', -1/e}$ is trivial modulo $\pi^{\flat, 1 - 1/e}$, 
    \[
        g'^{-1}hg = (\pi^{\flat',1/e})^{-e\lambda}\ov{\iota} (\pi^{\flat,1/e})^{e\lambda} \in \ov{p} \cdot \msf{N}(R^\flat)
    \]
    for some $\ov{p} \in \ov{\msf{P}}(R^{\flat +})$ that is trivial modulo $[\pi^{\flat, 1/p^m}]$. Thus, 
    \[
        [g_0] \equiv [hg_0] \in (\msf{G}/\msf{N})(R^{\flat+}/\pi^{\flat, 1/p^m})
    \]
    for sufficiently large $m$. Since $\msf{W}(0) \subset \msf{G}/\msf{N}$ and the $\msf{G}^\sigma$-action on $\msf{W}(0)$ is free, we have $h=1$. Thus, $\iota$ sends $(\pi^{\flat,1/e})^{e\lambda}g^{-1}$ to $(\pi^{\flat',1/e})^{e\lambda}g'^{-1}$, so we get the independence of the choice of $\pi^\flat$ and $u_\alpha^\flat$. This construction provides a map
    \[
        \mfr{W}(0)_\eta^\diamond \to \cl{W}(0)^\diamond, \quad
        g_0 \mapsto (\pi^{\flat,1/e})^{e\lambda}g^{-1}. 
    \]
    The right action of $h \in \msf{G}^\sigma$ on $\mfr{W}(0)_\eta$ sends $g_0$ to $h^{-1} g_0$, and it sends the depth-zero level structure $(\pi^{\flat,1/e})^{e\lambda}g^{-1}$ to $(\pi^{\flat,1/e})^{e\lambda}g^{-1} h$. Since the unique lift of $h^{-1} g_0$ is $h^{-1} g$, the above map is $\msf{G}^\sigma$-equivariant. 
\end{proof}

\section{Characterization in terms of integral models} \label{sec:dep0_integral_model}

To study the inner action on the special affinoid $\cl{W}(0)$, we need another characterization in terms of the \textit{specialization map} of a suitable integral model $\mfr{X}_b$. First, we provide a general construction of integral models. The main goal is to prove \Cref{prop2:characterization_intmodel}. 

\subsection{Equivariant compactifications} \label{ssec:equivcpt}

In this section, we review basic facts on equivariant compactifications of connected reductive groups. Here, let $\msf{G}$ be a connected reductive group over $\ov{k}$. 

\begin{defi}
    An equivariant compactification of $\msf{G}$ is a proper normal $\ov{k}$-variety $\msf{\ov{G}}$ with a $\msf{G}\times \msf{G}$-action and a $\msf{G}\times \msf{G}$-equivariant open immersion $\msf{G} \subset \msf{\ov{G}}$. 
\end{defi}

Here, $\msf{G}$ is equipped with a $\msf{G}\times \msf{G}$-action such that $(g_1,g_2)\cdot g = g_1gg_2^{-1}$ for $g,g_1,g_2\in \msf{G}$. 

\begin{prop} \label{prop:twistable}
    Let $X$ be a $\ov{k}$-scheme and let $P$ be a $\msf{G}$-torsor over $X$. The fppf sheaf $P \times^{\msf{G}}\msf{\ov{G}}$ is representable by a proper $X$-scheme. 
\end{prop}
\begin{proof}
    By Sumihiro's theorem (see \cite[Theorem 1, Lemma 8]{Sum74}), there is a $\msf{G} \times \msf{G}$-stable open covering $\{U_i\}_{i\in I}$ of $\msf{\ov{G}}$ such that $U_i$ admits a $\msf{G} \times \msf{G}$-equivariant ample line bundle. Then, $P \times^{\msf{G}} U_i$ is a scheme by \cite[Theorem 7]{BLR90}, where a descent datum is obtained by fixing an \'{e}tale local trivialization of $P$. By patching them together for all $i\in I$, we get a scheme representing $P \times^{\msf{G}}\msf{\ov{G}}$. 
\end{proof}

The following examples are important in application. Fix a Borel pair $\msf{T}\subset \msf{B}\subset \msf{G}$. 

\begin{enumerate}
    \item The wonderful compactification is an equivariant compactification for an adjoint group (see \cite{DCP83} and \cite{Str87}). Suppose that $\msf{G}$ is adjoint and let $\msf{\ov{G}}$ be the wonderful compactification. The set of $\msf{G}\times \msf{G}$-orbits of $\msf{\ov{G}}$ is bijective to the set of parabolic subgroups $\msf{B}\subset \msf{P}\subset \msf{G}$. Let $\msf{\ov{G}}_{\msf{P}} \subset \ov{\msf{G}}$ be the $\msf{G}\times \msf{G}$-orbit associated to $\msf{P}$. Let $\ov{\msf{P}}$ be the opposite of $\msf{P}$ and let $\msf{M}^\ad$ be the adjoint group of the Levi of $\msf{P}$. Then, 
    \[
        \msf{\ov{G}}_{\msf{P}} \cong (\msf{G}\times \msf{G})\times^{\msf{P}\times \ov{\msf{P}}} \msf{M}^\ad. 
    \]

    \item Toroidal compactifications are generalizations of the wonderful compactification to arbitrary connected reductive groups (see \cite[Section 6.2]{BK05}). Let $W$ be the Weyl group of $\msf{T}$ in $\msf{G}$ and let $\Sigma$ be a $W$-stable complete fan in the cocharacter group $X_*(\msf{T})$. For each $\Sigma$, we have a toroidal compactification $\ov{\msf{G}}_\Sigma$ filling in the following diagram: 
    \begin{center}
        \begin{tikzcd}
            \msf{T} \ar[r, hook] \ar[d, hook] & \msf{\ov{T}}_{\Sigma} \ar[d, hook] \\
            \msf{G} \ar[r, hook] & \msf{\ov{G}}_\Sigma. 
        \end{tikzcd}
    \end{center}
    Here, $\ov{\msf{T}}_\Sigma$ is the toric variety associated to $\Sigma$. When $\msf{G}$ is adjoint and $\Sigma$ is spanned by coroots of $\msf{G}$, $\ov{\msf{G}}_\Sigma$ recovers the wonderful compactification of $\msf{G}$. There is a natural map
    \[
        \msf{\ov{G}}_\Sigma \to \ov{\msf{G}}^\ad
    \]
    to the wonderful compactification of the adjoint group of $\msf{G}$. 

    Let $\ov{\msf{T}}^+_\Sigma \subset \ov{\msf{T}}_\Sigma$ be the open toric subvariety corresponding to the dominant chamber. Let $\msf{U}$ (resp.\ $\ov{\msf{U}}$) be the unipotent radical of $\msf{B}$ (resp.\ $\ov{\msf{B}}$). The map 
    \[
        \msf{U} \times \ov{\msf{T}}^+_\Sigma \times \ov{\msf{U}} \to \msf{\ov{G}}_\Sigma ,\quad
        (n,t,\ov{n}) \mapsto (n,\ov{n}) \cdot t
    \]
    is an open immersion. The translations of this open subset under $\msf{G}\times \msf{G}$ cover $\msf{\ov{G}}_\Sigma$. %(see \cite[6.2.3 Proposition (i)]{BK05}). 
\end{enumerate}

As a simple application, one can prove the existence of a Cartan decomposition using the toroidal compactification. 

\begin{lem} \label{lem:Cartandec}
    Let $V$ be a valuation ring over $\ov{k}$ with a fraction field $K$. Every $g\in \msf{G}(K)$ can be written as
    \[
        g=g_1tg_2 ,\quad g_1,g_2 \in \msf{G}(V) ,\quad t\in \msf{T}(K). 
    \]
\end{lem}

We say that a decomposition $g=g_1tg_2$ is a Cartan decomposition of $g$. 

\begin{proof}
    Let $\ov{\msf{G}}=\ov{\msf{G}}_\Sigma$ be a toroidal compactification associated to a $W$-stable complete fan $\Sigma$. By the valuative criterion, we can take an extension $\wtd{g}\in \msf{\ov{G}}(V)$. The translations of $\msf{U} \times \ov{\msf{T}}_{\Sigma}^+ \times \ov{\msf{U}} \subset \msf{\ov{G}}$ under $\msf{G}(\ov{k}) \times \msf{G}(\ov{k})$ cover $\msf{\ov{G}}$, so there is $h_1,h_2 \in \msf{G}(\ov{k})$ such that $h_1\wtd{g}h_2 \in \msf{U} \times \ov{\msf{T}}_{\Sigma}^+ \times \ov{\msf{U}}$. Thus, we can take $u_1\in \msf{U}(V)$, $u_2 \in \ov{\msf{U}}(V)$ and $t \in \ov{\msf{T}}_{\Sigma}^+(V)$ such that $g=(h_1^{-1} u_1)t(u_2h_2^{-1})$. 
\end{proof}

\subsection{Review of relative representability}

In this section, we review the author's representability criterion \cite{Tak26_rel}. A key ingredient is a Sen-theoretic datum on the base formal scheme, called good covers. First, we review the definition. For simplicity, we restrict to the $p$-adic setting. 

Let $R$ be a reduced excellent $p$-adically complete ring. Let $R_\bullet =  (R=R_0 \to R_1 \to \cdots)$ be a sequence of reduced finite $R$-algebras and let $\Gamma_\bullet = (\{1\} = \Gamma_0 \twoheadleftarrow \Gamma_1 \twoheadleftarrow \cdots)$ be a sequence of finite groups with surjective transition maps that acts on $R_\bullet$ in the sense of \cite[Definition 5.1]{Tak26_rel}. In the following, we set $U_n = \Spa(R_n)_\eta$ to be the generic fiber. 

\begin{defi} \textup{(\cite[Section 5.2]{Tak26_rel})}
    We say that $(R_\bullet, \Gamma_\bullet)$ is a good cover of $R$ if it satisfies the following conditions.  
    \begin{itemize}[align=right, leftmargin=2cm]
        \item[(quot)] $\Spd(R_{n+1}) \to \Spd(R_n)$ is a geometric quotient by $\Gamma_{n+1}^n$ and $\cl{O}(U_n) \to \cl{O}(U_{n+1})$ is injective and satisfies the condition in \cite[Lemma 3.13]{Tak26_rel} for $\Gamma_{n+1}^n$ for every $n\geq 0$. 
        \item[(perfd)] $\lim_{n \geq 0} \Spd(R_n)$ is represented by a perfectoid $R$-algebra $R_\infty$.  
        \item[(density)] the image of $\colim_{n < \infty} \cl{O}(U_n) \to \cl{O}(U_\infty)$ is dense, where $U_\infty = \Spa(R_\infty)_\eta$.
        \item[(unifsec)] we have a continuous $\cl{O}(U_n)$-linear section $s_{n,U_n} \colon \cl{O}(U_\infty) \to \cl{O}(U_n)$ for each $n \geq 0$ such that there is an open ideal $\cl{J} \subset \cl{O}^+(U_\infty)$ mapping to $\cl{O}^+(U_n)$ under $s_{n,U_n}$ for every $n$. 
    \end{itemize}
    We say that a reduced excellent $p$-adic formal scheme $\mfr{X}$ admits a good cover if $\mfr{X}$ is Zariski locally isomorphic to the formal spectrum of a reduced excellent $p$-adic ring admitting a good cover. 
\end{defi}

\begin{rmk}
    By the local nature of good covers (see \cite[Lemma 5.11, Lemma 5.13]{Tak26_rel}), the conditions of good covers are simplified in the $p$-adic setting. The class of $p$-adic rings admitting good covers is explored in \cite[Section 5.5]{Tak26_rel}. It includes $p$-adically smooth algebras over $O_F$ and their completion (see \cite[Lemma 5.7, Proposition 5.9]{Tak26_rel}).   
\end{rmk}

The representability criterion to be applied is the following. 

\begin{prop} \textup{(\cite[Proposition 5.45]{Tak26_rel})}
    \label{prop:finetrep}
    Let $\mfr{X}$ be a reduced excellent $p$-adic formal scheme admitting a good cover. Let $Y\to \mfr{X}^\diamond$ be a quasiseparated map of small $v$-sheaves. Let $Z \subset Y_\eta$ be a closed subsheaf that is finite \'{e}tale over $\mfr{X}_\eta^\diamond$. Suppose that the following condition holds. 
    \begin{quote}
        For every perfectoid ring $R$ with $\Spf(R)\to \mfr{X}$, there is a perfectly proper $R^\flat$-scheme $Y_{R^\flat}$ and a finite \'{e}tale $R^\flat[1/\xi_{R,0}]$-scheme $Z_{R^\flat[1/\xi_{R,0}]}$ with a morphism $Z_{R^\flat[1/\xi_{R,0}]} \to  Y_{R^\flat}$ such that we have the following commutative diagram
        \begin{center}
            \begin{tikzcd}
                Z \times_{\mfr{X}_\eta^\diamond} \Spd(R)^\an \ar[r] \ar[d,"\cong"] & Y \times_{\mfr{X}^\diamond} \Spd(R) \ar[d,"\cong"] \\
                Z_{R^\flat[1/\xi_{R,0}]}^{\diamondsuit}\times_{\Spec(R^\flat[1/\xi_{R,0}])^\diamondsuit} \Spd(R^\flat)^\an \ar[r] & Y_{R^\flat}^\diamond\times_{\Spec(R^\flat)^\diamond} \Spd(R^\flat). 
            \end{tikzcd}
        \end{center}
        that is compatible with tilting equivalence $\Spd(R) \cong \Spd(R^\flat)$. 
    \end{quote}
    Then, $Z^\clos$ is representable by a unique reduced excellent $p$-torsion free $p$-adic formal scheme $\mfr{Y}$ admitting a maximal good cover that is proper over $\mfr{X}$ with $\mfr{Y}_\eta$ finite \'{e}tale over $\mfr{X}_\eta$. 
\end{prop}

Here, the notion of maximal good covers is introduced in \cite[Section 5.4]{Tak26_rel}. The existence of such a cover is helpful due to the following full faithfulness of the $v$-sheafification functor. 

\begin{lem} \textup{(\cite[Lemma 5.34]{Tak26_rel})}
    \label{lem:uniqmapmaxlgood}
    Let $\mfr{X}$ be a reduced excellent $p$-torsion free $p$-adic formal scheme admitting a maximal good cover. Let $\mfr{Y}$ be a Noetherian formal scheme over $\bb{Z}_p$ and let $f\colon \mfr{X}^\diamond \to \mfr{Y}^\diamond$ be a map over $\Spd(\bb{Z}_p)$ such that $f \vert_{\mfr{X}^\diamond_\eta}$ is represented by a morphism $g\colon \mfr{X}_\eta \to \mfr{Y}^\ad$ of adic spaces. Then, $f$ is represented by a unique morphism $\mfr{X}\to \mfr{Y}$ and its restriction to $\mfr{X}_\eta$ is $g$. 
\end{lem}

\subsection{Depth-zero integral models} \label{ssec:moduli_level_structures}

In this section, we introduce integral models of the universal depth-zero level structures constructed as \textit{$v$-sheaf theoretic} flat closures. Here, we do not assume that $G$ is unramified, and let $\cl{G}$ be an arbitrary parahoric subgroup scheme over $O_F$. We employ the same notation as in \Cref{ssec:level_structure_in_char0}. 

By \Cref{prop:mapLW}, we have the following commutative diagram. 
\begin{center}
    \begin{tikzcd}
        \Sht_{\cl{G}^+} \ar[r] \ar[d] & \lbrack \ast/(L_W^+\cl{G^+})^\diamondsuit \rbrack \ar[d]\\
        \Sht_\cl{G} \ar[r] & \lbrack \ast/(L_W^+\cl{G})^\diamondsuit \rbrack . 
    \end{tikzcd}
\end{center}

Since $\cl{G}^+ \otimes F \cong \cl{G} \otimes F$, the natural map
\[
    \Sht_{\cl{G}^+} \hookrightarrow \Sht_{\cl{G}} \times_{[\ast/(L_W^+\cl{G})^\diamondsuit]} [\ast/(L_W^+\cl{G}^+)^\diamondsuit]
\]
is a monomorphism over $\Sht_{\cl{G}}$ by the Beauville-Laszlo lemma (cf.\ \cite[Proposition 6.11]{Tak26_rel}). The difference from the parahoric case in \cite[Section 6.2]{Tak26_rel} is that $\lbrack \ast/(L_W^+\cl{G^+})^\diamondsuit \rbrack \to  \lbrack \ast/(L_W^+\cl{G})^\diamondsuit \rbrack$ is not proper. 

\begin{lem} \label{lem:quotG}
    There is an isomorphism $(L_W^+\cl{G})^\diamondsuit/(L_W^+\cl{G}^+)^\diamondsuit \cong (\msf{G}^\sigma)^\diamondsuit$ of $v$-sheaves on $\Perf$. 
\end{lem}
\begin{proof}
    For each $\Spa(R,R^+) \in \Perf$, the reduction map
    \[
        \cl{G}(W_{O_F}(R)) \to \cl{G}(R) \to \msf{G}^\sigma(R)
    \]
    induces $(L_W^+\cl{G})^\diamondsuit \to \cl{G}_k^\diamondsuit \to (\msf{G}^\sigma)^\diamondsuit$. Since $\cl{G}$ is smooth, $(L_W^+\cl{G})^\diamondsuit \to \cl{G}_k^\diamondsuit$ is surjective. Since $\Ker(\cl{G}_k \to \msf{G}^\sigma)$ is smooth, $\cl{G}_k^\diamondsuit \to (\msf{G}^\sigma)^\diamondsuit$ is surjective in the \'{e}tale topology. Since $\cl{G}^+(W_{O_F}(R)) = \Ker(\cl{G}(W_{O_F}(R))\to \msf{G}^\sigma(R))$, we get the claim.
\end{proof}

We choose an equivariant compactification $\msf{G} \subset \msf{\ov{G}}$ and work over $\Perf_{\ov{k}}$. When $\msf{\ov{G}}$ is defined over $k$, we may work over $\Perf$ and replace $O_{\breve{F}}$ by $O_F$ for the subsequent arguments. 
Then, we have
\[
    [\ast/(L_W^+\cl{G}^+)^\diamondsuit] \cong [\msf{G}^\diamondsuit / (L_W^+\cl{G})^\diamondsuit] \subset [\msf{\ov{G}}^\diamondsuit / (L_W^+\cl{G})^\diamondsuit] \cong [\msf{\ov{G}}^\diamond / (L_W^+\cl{G})^\diamondsuit]. 
\]
This can be seen as a compactification of $[\ast/(L_W^+\cl{G}^+)^\diamondsuit]$ over $[\ast/(L_W^+\cl{G})^\diamondsuit]$. Here, we take the right action of $(L_W^+\cl{G})^\diamondsuit$ on $\msf{G}^\diamondsuit \subset \msf{\ov{G}}^\diamond$ as follows. 

\begin{conv} \label{conv:rightactLW}
    We put the right action of  $(L_W^+\cl{G})^\diamondsuit$ on $\msf{\ov{G}}^\diamond$ so that
    \[
        x \cdot h = \ov{h}^{-1}x ,\quad x \in \msf{\ov{G}}(R) ,\quad h \in \cl{G}(W_{O_F}(R))
    \]
    for each $\Spa(R,R^+) \in \Perf$. Here, $\ov{h} \in \msf{G}(R)$ is the reduction of $h$. 
\end{conv}

\begin{lem} \label{lem:desc_compactification}
    The $v$-stack $[\msf{\ov{G}}^\diamond/(L_W^+\cl{G})^\diamondsuit]$ is the $v$-sheafification of the functor sending $\Spa(R, R^+) \in \Perf_{\ov{k}}$ to the set of pairs $(P, s)$ where $P$ is a $\cl{G}$-torsor over $W_{O_F}(R)$ and $s \in (P\vert_{R} \times^{\cl{G}_k} \msf{\ov{G}})(R)$. The open embedding
    \[
        [\ast/(L_W^+\cl{G}^+)^\diamondsuit] \hookrightarrow [\msf{\ov{G}}^\diamond/(L_W^+\cl{G})^\diamondsuit]
    \]
    is the locus where $s$ comes from the trivialization of $P\vert_R \times^{\cl{G}_k} \msf{G}$. 
\end{lem}
\begin{proof}
    First, $[\ast / (L_W^+\cl{G})^\diamondsuit]$ sends $\Spa(R, R^+) \in \Perf$ to the groupoid of $\cl{G}$-torsors $P$ over $W_{O_F}(R)$. The fiber of $[\msf{\ov{G}}^\diamond/(L_W^+\cl{G})^\diamondsuit] \to [\ast / (L_W^+\cl{G})^\diamondsuit]$ at $P$ is the $v$-sheafification of the functor classifying a $\cl{G}_k$-equivariant map
    \[
        f \colon P \vert_R \to \msf{\ov{G}}. 
    \]
    By \Cref{conv:rightactLW}, $f$ is identified with a section $s \in (P \vert_R \times^{\cl{G}_k} \ov{\msf{G}})(R)$. Moreover, $f$ factors through $\msf{G}$ if and only if $s \in (P\vert_R \times^{\cl{G}_k} \msf{G})(R)$, so we get the second claim. 
\end{proof}

\begin{prop} \label{prop:repatinfty_depth_zero}
    Let $R$ be a perfectoid $O_{\breve{F}}$-algebra and let $\cl{P}$ be a local $\cl{G}$-shtuka over $\Spd(R)_{/O_F}$. There is a perfectly proper $R^\flat$-scheme $Y$ and a finite \'{e}tale $R^\flat[1/\xi_{R,0}]$-scheme $Z$ with a morphism $Z \to  Y$ such that we have the following commutative diagram
    \begin{center}
        \begin{tikzcd}
            \Sht_{\cl{G}^+} \times_{\Sht_{\cl{G}}, \cl{P}} \Spd(R)^\an \ar[r] \ar[d,"\cong"] & \lbrack \msf{\ov{G}}^{\diamond} /(L_W^+\cl{G})^\diamondsuit \rbrack \times_{\lbrack \ast/(L_W^+\cl{G})^\diamondsuit \rbrack, \cl{P}} \Spd(R) \ar[d,"\cong"] \\
            Z^{\diamondsuit}\times_{\Spec(R^\flat[1/\xi_{R,0}])^\diamondsuit} \Spd(R^\flat)^\an \ar[r] & Y^\diamond\times_{\Spec(R^\flat)^\diamond} \Spd(R^\flat). 
        \end{tikzcd}
    \end{center}
    that is compatible with tilting equivalence $\Spd(R) \cong \Spd(R^\flat)$. 
\end{prop}
\begin{proof}
    By \Cref{prop:GBKF}, $\cl{P}$ is given by a $\cl{G}$-BKF module $P$ over $R$. Let
    \[
        Y = P\vert_{R^\flat} \times^{\cl{G}_k} \ov{\msf{G}}
    \]
    be a proper $R^\flat$-scheme (see \Cref{prop:twistable}). By \Cref{lem:desc_compactification}, we have
    \begin{equation} \label{eq:write_with_Y}
        \lbrack \msf{\ov{G}}^{\diamond} /(L_W^+\cl{G})^\diamondsuit \rbrack \times_{\lbrack \ast/(L_W^+\cl{G})^\diamondsuit \rbrack, \cl{P}} \Spd(R) \cong Y^\diamondsuit \cong Y^\diamond. 
    \end{equation}
    Thus, the claim for $Y$ follows. Since $\Sht_{\cl{G}^+,\eta} \to \Sht_{\cl{G}, \eta}$ is finite \'{e}tale, there is a finite \'{e}tale homomorphism $(R^\flat[1/\xi_{R,0}], R^{\flat, +}) \to (S,S^+)$ such that
    \[
        \Spd(S,S^+) \cong \Sht_{\cl{G}^+} \times_{\Sht_{\cl{G}}, \cl{P}} \Spd(R)^\an
    \]
    under tilting equivalence. Let $Z=\Spec(S)$. By construction, we have a depth-zero level structure $\cl{P}^+$ on $\cl{P}$ over $\Spa(S,S^+)$. By \Cref{lem:finetdesc}, it corresponds to an $S$-valued point of $P\vert_{R^\flat} \times^{\cl{G}_k} \msf{G} \subset Y$. Thus, we have a morphism $Z\to Y$ with the given condition. 
\end{proof}

Now, we can apply the representability criterion developed in \cite{Tak26_rel}. 

\begin{prop} \label{prop:replevel}
    Let $\mfr{X}$ be a reduced excellent $p$-adic formal $O_{\breve{F}}$-scheme admitting a good cover and let $\cl{P}$ be a local $\cl{G}$-shtuka over $\mfr{X}^\diamond_{/O_F}$. The $v$-closure of the embedding
    \[
        \Sht_{\cl{G}^+}\times_{\Sht_{\cl{G}}, \cl{P}} \mfr{X}^{\diamond}_\eta \subset \lbrack \msf{\ov{G}}^{\diamond} /(L_W^+\cl{G})^\diamondsuit \rbrack \times_{\lbrack \ast/(L_W^+\cl{G})^\diamondsuit \rbrack, \cl{P}} \mfr{X}^\diamond
    \]
    is representable by a unique proper $p$-adic formal scheme $\mfr{Y}$ over $\mfr{X}$ admitting a maximal good cover with $\mfr{Y}_\eta$ finite \'{e}tale over $\mfr{X}_\eta$. Moreover, the $\msf{G}^\sigma$-action on $\mfr{Y}_\eta$ over $\mfr{X}_\eta$ extends to $\mfr{Y}$ over $\mfr{X}$. 
\end{prop}

We say that $\mfr{Y}$ is the flat moduli space of depth-zero level structures of type $\ov{\msf{G}}$ on $\cl{P}$. 

\begin{proof}
    Let $Y=\lbrack \msf{\ov{G}}^{\diamond} /(L_W^+\cl{G})^\diamondsuit \rbrack \times_{\lbrack \ast/(L_W^+\cl{G})^\diamondsuit \rbrack, \cl{P}} \mfr{X}^\diamond$ and $Z=\Sht_{\cl{G}^+} \times_{\Sht_{\cl{G}}, \cl{P}_\eta} \mfr{X}^\diamond_\eta$. Since $\ov{\msf{G}}$ is separated, $Y$ is separated over $\mfr{X}^\diamond$ by \cite[Proposition 10.11]{Sch17}. Moreover, $Z$ is finite \'{e}tale over $\mfr{X}_\eta^\diamond$, so $Z\hookrightarrow Y_\eta$ is closed. By \Cref{prop:repatinfty_depth_zero}, we may apply \Cref{prop:finetrep} to $Z \hookrightarrow Y$ to see the representability of the $v$-closure of $Z$ in $Y$.

    The right action $x\cdot g = xg$ of $\msf{G}^\sigma$ on $\msf{\ov{G}}$ induces $\msf{G}^\sigma$-actions on $Y$ and $Z$. It is compatible with the natural action on $Z$ by \Cref{lem:finetdesc}. In particular, the $\msf{G}^\sigma$-action on $Z$ extends to the $v$-closure of $Z$ in $Y$. It uniquely provides the action on $\mfr{Y}$ by \Cref{lem:uniqmapmaxlgood}. 
\end{proof}

\begin{rmk} \label{rmk:defined_O_F}
    If $\msf{\ov{G}}$ can be defined over $k$, the above arguments work over $\Perf$ and $\mfr{X}$ can be taken as a reduced excellent $p$-adic formal $O_{F}$-scheme admitting a good cover. 
\end{rmk}

\begin{prop} \label{prop:finite_replevel}
    In the setting of \Cref{prop:replevel}, suppose that there is a $\msf{G}$-equivariant affine open subset $\msf{U} \subset \msf{\ov{G}}$ such that the natural map
    \[
        \mfr{Y}_s^\diamond \subset \lbrack \msf{\ov{G}}^{\diamond} /(L_W^+\cl{G})^\diamondsuit \rbrack \times_{\lbrack \ast/(L_W^+\cl{G})^\diamondsuit \rbrack, \cl{P}} \mfr{X}_s^\diamond
    \]
    factors through $\lbrack \msf{U}^{\diamondsuit} /(L_W^+\cl{G})^\diamondsuit \rbrack \times_{\lbrack \ast/(L_W^+\cl{G})^\diamondsuit \rbrack, \cl{P}} \mfr{X}_s^\diamond$. Then, $\mfr{Y}$ is finite over $\mfr{X}$. 
\end{prop}
\begin{proof}
    We may assume that $\mfr{X} = \Spf(R)$ is affine. Let $A = (R_\red)^\perf$ and let $P$ be a $\cl{G}$-BKF module over $A$ associated to $\cl{P}\vert_{\Spd(A)}$ via \Cref{prop:GBKF}. Let $Y = P\vert_A \times^{\cl{G}_k} \msf{\ov{G}}$. By \cite[Lemma 3.31]{Gle24} and \eqref{eq:write_with_Y}, $\mfr{Y}_s^\perf$ is a closed subscheme of $Y^\perf$.
    
    Let $V = P\vert_A \times^{\cl{G}_k} \msf{U}$. By assumption, $V$ is affine and $(\mfr{Y}_s^\perf)^\diamond \subset Y^\diamond$ factors through $V^\diamondsuit$. We will show $\mfr{Y}_s^\perf \subset V^\perf$. Let $x \in \mfr{Y}_s^\perf$ and let $\kappa$ be the residue field of $x$. Then, the point of $Y^\diamond$ corresponding to the composition
    \[
        \Spa(\kappa((t^{1/p^\infty})), \kappa\llbracket t^{1/p^\infty} \rrbracket) \to \Spec(\kappa) \to \mfr{Y}_s^\perf \to Y
    \]
    lies in $V^\diamondsuit$. It follows that $x \in V$ since $V$ is affine. Then $\mfr{Y}_s^\perf$ is a closed subscheme of $V$, so $\mfr{Y}_s$ is affine. Since $\mfr{Y}_s$ is proper over $\mfr{X}_s$, $\mfr{Y}_s$ is finite over $R / pR$, so $\mfr{Y}$ is finite over $\mfr{X}$. 
\end{proof}

\begin{rmk}
    In \Cref{prop:Drlevel} and \Cref{prop:Gamma_1(p)-level}, we interpret moduli spaces of depth-zero Drinfeld level structures in the framework of \Cref{prop:replevel}. The finiteness of such moduli spaces can be explained by \Cref{prop:finite_replevel}. % by the choice $\msf{U} = \Mat_n \subset \bb{P}(\Mat_n \oplus \triv) = \ov{\msf{G}}$. 
\end{rmk}

Now, we can introduce an integral model of the moduli of depth-zero level structures on $\mfr{Q}^\univ$. Here, we assume that $\cl{G}$ is reductive and set $\msf{G} = \cl{G}_{\ov{k}}$. 

\begin{const} \label{const:Xb}
    Let $\ov{\msf{G}}$ be an equivariant compactification of $\msf{G}$. Since $R_{\cl{G},\mu}$ admits a good cover as a $p$-adic ring by \cite[Proposition 5.5]{Tak26_rel} and \cite[Lemma 5.7]{Tak26_rel}, 
    we can take the flat moduli space
    $
        \mfr{X}_b^{\msf{\ov{G}},\pre}
    $
    of depth-zero level structures of type $\ov{\msf{G}}$ on $\mfr{Q}^\univ$. Then, let
    \[
        \mfr{X}_b^{\msf{\ov{G}}} = (\mfr{X}_b^{\msf{\ov{G}},\pre})^\wedge \to \Spf(R_{\cl{G}, \mu})
    \]
    be the completion along the maximal ideal of $R_{\cl{G}, \mu}$. By \Cref{prop:replevel}, $\mfr{X}_b^{\msf{\ov{G}}}$ is equipped with a $\msf{G}^\sigma$-action over $\Spf(R_{\cl{G}, \mu})$ and we have
    \[
        \mfr{X}_{b, \eta}^{\msf{\ov{G}}} \cong \Spf(R_{\cl{G}, \mu})_\eta \times_{\Sht_{\cl{G}}} \Sht_{\cl{G}^+} ,\quad
        (\mfr{X}_b^{\msf{\ov{G}}})^\perf_\red \subset \msf{\ov{G}}^\perf
    \]
    by \cite[Proposition 18.3.1]{SW20}. The latter is a closed immersion since $(\mfr{X}_b^{\msf{\ov{G}}})^\perf_\red$ is perfectly proper over $\ov{k}$. When $\ov{\msf{G}}$ is clear from the context, $\mfr{X}_b^{\msf{\ov{G}}}$ (resp.\ $\mfr{X}_b^{\msf{\ov{G}},\pre}$) is simply denoted by $\mfr{X}_b$ (resp.\ $\mfr{X}_b^\pre$) by abuse of notation. 
\end{const}

\subsection{Relation to special affinoids}

Fix an equivariant compactification $\msf{\ov{G}}$ and the associated integral model $\mfr{X}_b$. First, we show that \Cref{thm:reduction_of_specaff} can be extended to $\mfr{W}(0)$. 

\begin{prop} \label{prop:mapWX}
    There is a unique $\msf{G}^\sigma$-equivariant map
    \[
        \mfr{W}(0) \to \mfr{X}_b
    \]
    that restricts to the map in \Cref{thm:reduction_of_specaff} on the generic fiber. 
\end{prop}
\begin{proof}
    The uniqueness and the $\msf{G}^\sigma$-equivariance are automatic since $\mfr{W}(0)$ is smooth and $\mfr{X}_b$ is separated. In particular, it is enough to construct a map $\mfr{W}(0) \to \mfr{X}_b^{\pre}$. By \cite[Proposition 5.9]{Tak26_rel} and \Cref{lem:uniqmapmaxlgood}, it is enough to construct a map
    \[
        \mfr{W}(0)^\diamond \to \lbrack \msf{\ov{G}}^{\diamond} /(L_W^+\cl{G})^\diamondsuit \rbrack \times_{\lbrack \ast/(L_W^+\cl{G})^\diamondsuit \rbrack, \cl{P}} \Spd(R_{\cl{G}, \mu})
    \]
    extending the map in \Cref{thm:reduction_of_specaff}. 

    As previously, let $(R,R^+)$ be a perfectoid Huber pair over $\mfr{W}(0)$ and choose $\pi^\flat, u_\alpha^\flat \in R^{\flat +}$. Then, the associated $\cl{G}$-BKF module is given by
    \[
        P_{R^+}=(\cl{G}\otimes W_{O_F}(R^{\flat+}), \mu([\pi^\flat]-\pi)\prod_{\alpha\in \Phi_{\mu<0}} i_\alpha([u_\alpha^\flat])w\sigma). 
    \]
    Let $u_{\alpha, 1}^\flat = \pi^{\flat, -r_\alpha} u_{\alpha}^\flat \in R^{\flat +}$. The map $\Spf(R^+/\pi^{1/e}) \to \msf{W}(0)$ provides $g_0\in \msf{G}(R^+/\pi^{1/e})$ with
    \[
        g_0^{-1}\sigma(g_0) \equiv \prod_{\alpha\in \Phi_{\mu<0}} i_\alpha(u_{\alpha,1}^\flat)w \pmod{\pi^{\flat, 1/e}}. 
    \]
    By \Cref{lem:uniqext}, $g_0$ can be uniquely lifted to $g\in \msf{G}(R^{\flat +})$ such that
    \[
        g^{-1}\sigma(g) = \prod_{\alpha\in \Phi_{\mu<0}} i_\alpha(u_{\alpha,1}^\flat)w. 
    \]
    Recall the description in \Cref{lem:desc_compactification}. We will show that the map
    \[
        (\cl{G}\otimes W_{O_F}(R^{\flat}), (\pi^{\flat,1/e})^{e\lambda} g^{-1}) \colon 
        \Spa(R, R^+) \to \lbrack \msf{\ov{G}}^{\diamond} /(L_W^+\cl{G})^\diamondsuit \rbrack \times_{\lbrack \ast/(L_W^+\cl{G})^\diamondsuit \rbrack, \cl{P}} \Spd(R_{\cl{G}, \mu})
    \]
    is independent of the choice of $\pi^\flat$ and $u_{\alpha,1}^\flat$. Here, $e\lambda \in X_*(\msf{T})$ extends to
    \[
        e\lambda \colon \bb{A}^1 \to \msf{\ov{G}}
    \]
    due to the properness of $\msf{\ov{G}}$ and $(\pi^{\flat,1/e})^{e\lambda} \in \msf{\ov{G}}$ denotes the image of $\pi^{\flat, 1/e}$. When $(R, R^+)$ is in characteristic $0$, this recovers the construction of the map in \Cref{thm:reduction_of_specaff}, so it is enough to check the well-definedness of this map. 

    As $\lbrack \msf{\ov{G}}^{\diamond} /(L_W^+\cl{G})^\diamondsuit \rbrack \times_{\lbrack \ast/(L_W^+\cl{G})^\diamondsuit \rbrack, \cl{P}} \Spd(R_{\cl{G}, \mu})$ is separated, an equality of two maps into it can be checked at each geometric point. Since the characteristic $0$ case is already treated in \Cref{thm:reduction_of_specaff}, we may assume that $(R, R^+)$ is a geometric point in characteristic $p$. In this case, the choice of $\pi^\flat$ and $u_{\alpha,1}^\flat$ is unique, so the claim is obvious. 
\end{proof}

%The map of underlying spaces is given as follows. 

\begin{cor} \label{cor:explicit_underlying}
    Let $x_\lambda = (e\lambda)(0) \in \msf{\ov{G}}$ be the limit point of $e\lambda$. The action of $\msf{\ov{N}}\times \msf{N}$ on $x_{\lambda}$ is trivial, and the underlying map of $\mfr{W}(0) \to \mfr{X}_b$ is given by the composition (up to perfection)
    \[
        \msf{W}(0) \hookrightarrow \msf{G} / \msf{N} \to \msf{\ov{G}} ,\quad
        g \mapsto x_{\lambda} g^{-1}. 
    \]
\end{cor}
\begin{proof}
    The first claim is a general fact. Let us recall the proof here. Let $n\in \msf{N}$ be a closed point. The map $\bb{G}_m \to \msf{N}$ sending $t$ to $\Ad(t^{e\lambda})n$ extends to $\bb{A}^1$ and sends $0$ to $1$. Thus, the extension of the map $t\mapsto t^{e \lambda}n$ sends $0$ to $x_{\lambda}$. Thus, we have $x_{\lambda} \cdot n = x_{\lambda}$. By the same argument, $\ov{n} \cdot x_\lambda = x_\lambda$ for every $\ov{n} \in \ov{\msf{N}}$. 

    By the proof of \Cref{prop:mapWX}, each geometric point 
    \[
        \Spa(C, C^+) \to \msf{W}(0)
    \]
    corresponding to $g_0 \in \msf{G}(C^+)$ is sent to $x_\lambda \cdot g_0^{-1}$. Thus, the underlying map is as given in the statement up to perfection. 
\end{proof}

The main role of the integral model $\mfr{X}_b$ is the following characterization of $\cl{W}(0)$. 

\begin{thm} \label{thm:role_of_integralmodel}
    When $\msf{\ov{G}}$ is a toroidal compactification, the image of 
    \[
        \msf{W}(0) \to \mfr{X}_{b, \red}
    \]
    is open. Let $U_b$ be the open image. Then, $U_b$ is affine and 
    \[
        \mfr{W}(0) \to \mfr{X}_b\vert_{U_b}
    \]
    is finite. In particular, $\cl{W}(0) = (\mfr{X}_b\vert_{U_b})_\eta$. 
\end{thm}

As an immediate corollary, we recover Yoshida's construction of special affinoids. 

\begin{cor} 
    Let $\mfr{X}_{b,e}^\norm$ be the normalization of $\mfr{X}_b \otimes O_{\breve{F}_e}$. When $\msf{\ov{G}}$ is toroidal, we have an open immersion
    \[
        \mfr{W}(0) \subset \mfr{X}_{b, e}^\norm. 
    \]
\end{cor}

\subsection{Proof of \Cref{thm:role_of_integralmodel}}

In this section, we assume that $\msf{\ov{G}}$ is toroidal and prove \Cref{thm:role_of_integralmodel}. This assumption is imposed to control the stabilizer of $x_\lambda$ as follows. 

\begin{lem} \label{lem:stab_xlambda}
    The stabilizer of $x_{\lambda} \in \msf{\ov{G}}$ under $\msf{G} \times \msf{G}$ is contained in $\msf{\ov{P}}\times_{\msf{M}^\ad} \msf{P}$. 
\end{lem}
\begin{proof}
    There is a natural map $\ov{\msf{G}} \to \ov{\msf{G}}^\ad$ to the wonderful compactification of $\msf{G}^\ad$. For the wonderful compactification, the claim is proved in \cite[Proposition 2.2]{Str87}. %, so the claim follows. 
\end{proof}

Since $\mfr{X}_b$ is constructed as an abstract flat closure, it is hard to directly study the geometry of $\mfr{X}_b$. Instead, we study the \textit{specialization map} of $\mfr{X}_b$. In particular, we will deduce all the claims from the following property of the specialization map. 

\begin{prop} \label{prop:specialization_Xb}
    Let $x \in \mfr{X}_{b, \eta}^\pre$ be a geometric point. Then, $x \in \cl{W}(0)$ if and only if the specialization $\spc(x)$ of $x$ lies in the image of $\msf{W}(0) \to \mfr{X}_{b, \red}^\pre$. 
\end{prop}

%For simplicity, we will pretend morphisms defined only up to perfection to be usual morphisms between $\ov{k}$-schemes. This will cause no problems in our argument. 

\begin{proof}[Proof of \Cref{thm:role_of_integralmodel} from \Cref{prop:specialization_Xb}]
    Recall that $\mfr{X}_b^\pre$ is a $p$-adic $p$-torsion free formal scheme that is proper over $R_{\cl{G}, \mu}$ (as a $p$-adic ring). In particular, the specialization map is surjective and closed. Since $\cl{W}(0) \subset \mfr{X}_{b, \eta}^\pre$ is open, 
    \[
        Z = \spc(\mfr{X}_{b, \eta}^\pre - \cl{W}(0)) \subset \lvert \mfr{X}_{b, \red}^\pre \rvert
    \]
    is closed. By \Cref{prop:specialization_Xb}, $\msf{W}(0)$ maps surjectively into the complement of $Z$, so the image $U_b$ of $\msf{W}(0) \to \mfr{X}_{b, \red}$ is open in $\lvert \mfr{X}_{b, \red}^\pre \rvert$. Up to perfection, we have the following diagram. 
    \begin{center}
        \begin{tikzcd}
            \msf{W}(0) \ar[r, two heads] \ar[d, two heads] & U_b \ar[r, hook] & \msf{\ov{G}} \ar[d] \\
            \msf{W}(0) / \msf{Z}_{\msf{M}^{w\sigma}} \ar[r, hook] & \msf{G} / \msf{Z}_{\msf{M}} \msf{N} \ar[r, hook] & \msf{\ov{G}}^\ad. 
        \end{tikzcd}
    \end{center}
    It follows that up to perfection, the finite \'{e}tale quotient $\msf{W}(0) \to \msf{W}(0) / \msf{Z}_{\msf{M}^{w\sigma}}$ factors through $U_b$. In particular, $U_b$ is finite over $\msf{W}(0) / \msf{Z}_{\msf{M}^{w\sigma}}$ up to perfection, so $U_b$ is affine and $\msf{W}(0) \to U_b$ is finite. Since $U_b \subset \mfr{X}_{b, \red}$, we have $\mfr{X}_b\vert_{U_b} = \mfr{X}_b^\pre\vert_{U_b}$ and $\mfr{X}_{b}\vert_{U_b}$ is $p$-adic. Thus, $\mfr{W}(0) \to \mfr{X}_b\vert_{U_b}$ is finite. 
\end{proof}

\begin{rmk} \label{rmk:Xwopen}
    When the stabilizer of $x_{\lambda}$ equals $\msf{\ov{P}}\times_{\msf{M}^\ad} \msf{P}$, we have (up to perfection)
    \[
        \msf{W}(0)/\msf{Z}_\msf{M}^{w\sigma} = U_b \subset \mfr{X}_{b,\red}. 
    \]
    In particular, if $x_{\lambda}$ lies in a closed $\msf{G}\times \msf{G}$-orbit, $U_b$ equals the Deligne-Lusztig variety
    \[
        \msf{X}(w)=\{g\in \msf{G}/\msf{B} \mid g^{-1}\sigma(g) \in \msf{B} w \msf{B} \}
    \]
    up to perfection. 
\end{rmk}

Now, we will prove \Cref{prop:specialization_Xb}. Take a geometric point
\[
    x \colon \Spa(C, C^+) \to \mfr{X}_{b, \eta}^\pre
\]
and suppose that the specialization $\spc(x)$ lies in the image of $\msf{W}(0) \to \mfr{X}_{b, \red}^\pre$. Choose $\pi^\flat$ and $u_{\alpha}^\flat$ as in \Cref{lem:expBKF}. The associated $\cl{G}$-BKF module is
\[
    P_{x}=(\cl{G}\otimes W_{O_F}(C^{\flat+}), \mu([\pi^\flat]-\pi)\prod_{\alpha\in \Phi_{\mu<0}} i_\alpha([u_\alpha^\flat])w\sigma)
\]
and the depth-zero level structure of $x$ corresponds to $g\in \msf{G}(C^\flat)$ such that 
\[
    g\sigma(g)^{-1} = \mu(\pi^\flat) \prod_{\alpha \in \Phi_{\mu<0}} i_\alpha(u_\alpha^\flat)w. 
\]
Then, $\spc(x) \in \msf{\ov{G}}$ is given by the specialization of $g \in \msf{G}(C^\flat) \subset \msf{\ov{G}}(C^\flat)$. 

To compute this specialization, take a Cartan decomposition
\[
    g=h_1th_2 ,\quad h_1, h_2 \in \msf{G}(C^{\flat +}) ,\quad t \in \msf{T}(C^\flat)
\]
as in \Cref{lem:Cartandec}. Let $\kappa$ be the residue field of $C^+$. Let $\ov{h}_1, \ov{h}_2 \in \msf{G}(\kappa)$ be the reductions and let $\ov{t} \in \msf{\ov{T}}(\kappa)$ be the specializations in the associated toric variety $\ov{\msf{T}} \subset \ov{\msf{G}}$. Then, we have
\[
    \spc(x) = \ov{h}_1 \ov{t} \ov{h}_2 \in \Img(\msf{W}(0) \to \msf{\ov{G}}). 
\]
By \Cref{cor:explicit_underlying}, there is $y \in \msf{W}(0)$ such that $\ov{h}_1 \ov{t} \ov{h}_2  = x_\lambda y^{-1}$.

\begin{lem}
    Let $\mfr{m}^\flat$ be the maximal ideal of $C^{\flat+}$ and let 
    \[
        \msf{H}(\mfr{m}^\flat) = \Ker(\msf{H}(C^{\flat +}) \to \msf{H}(\kappa))
    \]
    for an algebraic group $H$ over $\ov{k}$. By modifying a Cartan decomposition, we may assume that 
    \[  
        h_1 \in \msf{N}(\mfr{m}^\flat) ,\quad
        \ov{t}=x_{\lambda} ,\quad
        \ov{h}_2y \in \msf{M}. 
    \]
\end{lem}
\begin{proof}
    First, $\ov{t}$ lies in the same $\msf{T} \rtimes W$-orbit as $x_\lambda$ by \cite[Proposition 6.2.3]{BK05}, so we may assume that $\ov{t}=x_{\lambda}$. Then, $(\ov{h}_1, \ov{h}_2y)$ stabilizes $x_{\lambda}$, so we have 
    \[
        \ov{h}_1\in \msf{\ov{P}} ,\quad \ov{h}_2y \in \msf{P}
    \]
    by \Cref{lem:stab_xlambda}. Since $\ov{t}=x_{\lambda}$, we have 
    \[
        \alpha(t) \in \mfr{m}^\flat \quad (\alpha \in \Phi_{\msf{N}}) ,\quad
         \alpha(t) \in (C^{\flat+})^\times \quad (\alpha \in \Phi_{\msf{M}}). 
    \]
    Here, $\Phi_{(-)}$ denotes the set of roots inside. Let 
    $
        \ov{h}_2y = nm
    $
    be the Levi decomposition and take a lift $\wtd{n} \in \msf{N}(C^{\flat +})$ of $n$. Then,
    \[
        t\wtd{n}t^{-1}\in \msf{N}(\mfr{m}^\flat), 
    \]
    so we may replace $(h_1,h_2)$ with $(h_1t\wtd{n}t^{-1},\wtd{n}^{-1}h_2)$ to assume that $\ov{h}_2y \in \msf{M}$ without changing $\ov{h}_1$. Moreover, since $\ov{h}_1 \in \msf{\ov{P}}$, $h_1$ lies in the big cell
    \[
        \msf{N}(\mfr{m}^\flat) \times \msf{\ov{P}}(C^{\flat+}) \subset \msf{G}(C^{\flat +}). 
    \]
    Let $h_1 = n\ov{p}$. We may replace $(h_1,h_2)$ with $(n, t^{-1}\ov{p}th_2)$ to get $h_1 \in \msf{N}(\mfr{m}^\flat)$ and $\ov{h}_2y \in \msf{M}$. 
\end{proof}

Fix this Cartan decomposition. Since $g\sigma(g)^{-1} = \mu(\pi^\flat) \prod_{\alpha \in \Phi_{\mu<0}} i_\alpha(u_\alpha^\flat)w$, we have
\[
    th_2\sigma(h_2)^{-1}\sigma(t)^{-1} = h_1^{-1} \mu(\pi^\flat)\prod_{\alpha \in \Phi_{\mu<0}} i_\alpha(u_\alpha^\flat)w\sigma(h_1).
\]
Since $\ov{h}_2y \in \msf{M}$, $\ov{h}_2 \sigma(\ov{h}_2)^{-1} \in \msf{U}_{\mu<0} \msf{M} w$. In particular, $h_2\sigma(h_2)^{-1}$ lies in the big cell $\msf{N}(\mfr{m}^\flat) \times \msf{\ov{P}}(C^{\flat+})w$. Let us write $h_2\sigma(h_2)^{-1} = n\ov{p}w$ with $n \in \msf{N}(\mfr{m}^\flat)$ and $\ov{p} \in  \msf{\ov{P}}(C^{\flat+})$. Then, we have 
\begin{equation}
    t\ov{p}\Ad(w\sigma)(t)^{-1} w = \Ad(t)(n)^{-1}h_1^{-1}\mu(\pi^\flat)w \prod_{\alpha\in \Phi_{\mu<0}} \Ad(w^{-1})(i_\alpha(u_\alpha^\flat))\sigma(h_1).\label{eq:fund}
\end{equation}
The left-hand side lies in $\msf{\ov{P}} w$ and the right-hand side lies in $\msf{P}  w \sigma(\msf{P})$ by \Cref{lem:MNprop}.

\begin{lem}\textup{(\cite[Theorem 1.4]{Ric92})}
    We have $\msf{P}w\sigma(\msf{P}) \cap \ov{\msf{P}}w\sigma(\msf{P}) = w \sigma(\msf{P})$.
\end{lem}

Thus, each side of \eqref{eq:fund} lies in $w\sigma(\msf{P})$. Comparing the image of \eqref{eq:fund} through the (Levi) quotient $w\sigma(\msf{P}) \to \msf{M}$, we get
\[
    tm\Ad(w\sigma)(t)^{-1} = \mu(\pi^\flat)
\]
where $m$ is the Levi part of $\ov{p}$. First, we have $m\in \msf{T}(C^{\flat+})$. Let $\nu^\flat$ be the additive valuation on $C^\flat$ normalized so that $\nu^\flat(\pi^\flat)=1$ and let $\nu_\alpha = \nu^\flat(\alpha(t))$ for $\alpha \in X^*(\msf{T})$. Then, we have 
\[
    \nu_{\alpha} - q\nu_{(w\sigma)^{-1}\alpha} = \langle \alpha, \mu \rangle. 
\]
As in the proof of \Cref{lem:lambdadef}, we see $\nu_\alpha = \langle \alpha, \lambda \rangle$. Fix $\pi^{\flat, 1/e} \in C^{\flat +}$. Then, it follows that $ (\pi^{\flat,1/e})^{-e\lambda}t \in \msf{T}(C^{\flat+})$. By replacing $(t,h_2)$ with $((\pi^{\flat,1/e})^{e\lambda}, (\pi^{\flat,1/e})^{-e\lambda}th_2)$, we may assume 
\[
    t=(\pi^{\flat,1/e})^{e\lambda}. 
\]
Then, we have $m=1$. Since both sides of \eqref{eq:fund} lie in $w\sigma(\msf{P})$, we have
\begin{equation} \label{eq:condition_on_unipotent}
    \ov{p}\in \msf{\ov{N}}\cap \Ad(w\sigma)(\msf{N}) = \msf{U}_{\mu < 0} ,\quad
    \Ad(t)(n)^{-1}h_1^{-1} \in \msf{N} \cap \Ad(w\sigma)(\msf{N}). 
\end{equation}
In particular, the latter commutes with $\mu(\pi^\flat)$ by \Cref{lem:MNprop}. Then \eqref{eq:fund} is rewritten as
\begin{equation}
    \Ad((\pi^{\flat,1/e})^{qw\sigma(e\lambda)})(\ov{p})=\Ad((\pi^{\flat,1/e})^{e\lambda})(n)^{-1}h_1^{-1}\prod_{\alpha \in \Phi_{\mu<0}} i_\alpha(u_\alpha^\flat) \Ad(w\sigma)(h_1). \label{eq:fundunip}
\end{equation}

\begin{defi}
    For each rational cocharacter $\nu \in X_*(\msf{T})_{\bb{Q}}$ and a product of root groups $\msf{V}=\prod_{\alpha \in J} \msf{U}_\alpha$ for an ordered subset $J\subset \Phi$, let
    \[
        \msf{V}(\nu) = \biggl\{ \prod_{\alpha \in J} i_\alpha(c_\alpha) \;\bigg\vert\; \nu^\flat(c_\alpha) \geq \langle \alpha, \nu \rangle \biggr\} \subset \msf{V}(C^\flat). 
    \]
    When $J$ is closed under addition, $\msf{V}(\nu)$ is a subgroup and independent of the order of $J$. 
\end{defi}

Now, the left-hand side of \eqref{eq:fundunip} lies in $\msf{U}_{\mu < 0}(qw\sigma\lambda)$. To show $x\in \mfr{W}(0)_\eta$, we need to show that $\prod_{\alpha \in \Phi_{\mu<0}} i_\alpha(u_\alpha^\flat)$ lies in the same subgroup. 

\begin{lem}
    We have $h_1 \in \msf{N}(\lambda)$. 
\end{lem}
\begin{proof}
    We may write $h_1 = h^{-}_1 h^+_1$ with 
    \[
        h^-_1\in \msf{N} \cap \Ad(w\sigma)(\msf{N}) ,\quad
        h^+_1\in \msf{N} \cap \Ad(w\sigma)(\msf{\ov{N}}) = \msf{U}_{\mu > 0}. 
    \]
    By \eqref{eq:condition_on_unipotent}, we have $h^+_1 \in \msf{U}_{\mu > 0}(\lambda)$. By the projection of \eqref{eq:fundunip} along $\pi_w$ (see \eqref{eq:definition_piw}), we have
    \[
        (h_1^-)^{-1}\pi_w(\Ad(w\sigma)(h_1^-h_1^+)) \in (\msf{N} \cap \Ad(w\sigma)(\msf{N}))(\lambda). 
    \]
    Since $h^+_1 \in \msf{U}_{\mu > 0}(\lambda)$, we have 
    \[
        \pi_w(\Ad(w\sigma)(h_1^+)) \in (\msf{N}\cap \Ad(w\sigma)(\msf{N}))(qw\sigma\lambda) = (\msf{N}\cap \Ad(w\sigma)(\msf{N}))(\lambda)
    \]
    since $\msf{N}\cap \Ad(w\sigma)(\msf{N})$ commutes with the image of $\mu$. Thus, we have 
    \[
        (h_1^-)^{-1}\phi_w(h_1^-) \in (\msf{N}\cap \Ad(w\sigma)(\msf{N}))(\lambda). 
    \]
    Since $\lambda - qw\sigma\lambda=\mu$, $(\msf{N}\cap \Ad(w\sigma)(\msf{N}))(\lambda)$ is stable under $\phi_w$. By \Cref{lem:phiwnilp}, it follows that
    \[
        (h_1^-)^{-1} = (h_1^-)^{-1}\phi_w(h_1^-) \cdot \phi_w((h_1^-)^{-1}\phi_w(h_1^-)) \cdots \in (\msf{N}\cap \Ad(w\sigma)(\msf{N}))(\lambda). 
    \]
    Thus, we have $h_1 = h_1^-h_1^+\in \msf{N}(\lambda)$. 
\end{proof}

Now, we have $\Ad(w\sigma)(h_1) \in \Ad(w\sigma)(\msf{N})(qw\sigma\lambda)$. By \eqref{eq:fundunip}, we have 
\[
    \Ad(t)(n)^{-1}h_1^{-1} \cdot \prod_{\alpha \in \Phi_{\mu<0}} i_\alpha(u_\alpha^\flat) \in \Ad(w\sigma)(\msf{N})(qw\sigma\lambda).
\]
By \eqref{eq:condition_on_unipotent}, we have $\prod_{\alpha \in \Phi_{\mu<0}} i_\alpha(u_\alpha^\flat) \in \msf{U}_{\mu<0}(qw\sigma\lambda)$ and we get the claim. 

\section{Group actions on special affinoids} \label{sec:group_action_on_special_affinoids}

\subsection{The inner action} \label{ssec:actinn}

Recall that the stabilizer of $[1]\in X_\mu(b)$ is $G_b(F) \cap \cl{G}(O_{\breve{F}})$. In this section, we show that $\cl{W}(0)$ is stable under $G_b(F) \cap \cl{G}(O_{\breve{F}})$ and compute the action on the reduction $\msf{W}(0)$. 

\begin{defi}
    Let $\tau$ be the unique face of the codominant Weyl chamber of $\cl{A}(\breve{G}, \breve{T})$ that contains $-\lambda$ in its interior and let $\mfr{f} = \mfr{a}\cap \tau$. 
\end{defi}

\begin{lem} \label{lem:charf}
    The facet $\mfr{f}$ is the minimal facet containing $(b\sigma)^ko$ for all $k\geq 0$.  
\end{lem}
\begin{proof}
    Since $\msf{M}$ is the centralizer of $\lambda$, $\tau$ is parallel to $X_*(\msf{Z}_{\msf{M}}/\msf{Z}_{\msf{G}})\otimes \bb{R}$ and stable under $w\sigma$. Since $\mu \in X_*(\msf{Z}_\msf{M})$ by \Cref{lem:lambdapst} and $b\sigma\mfr{a} = \mfr{a}$, $\mfr{f}$ is stable under $b\sigma$. Thus, $\mfr{f}$ contains $(b\sigma)^ko$ for every $k\geq 0$. Let $\mfr{f}_0 \subset \mfr{f}$ be the minimal face containing $(b\sigma)^ko$ for all $k\geq 0$. Let $V \subset X_*(T/Z_G) \otimes \bb{R}$ be the maximal linear subspace parallel to $\mfr{f}_0$. By the minimality of $\mfr{f}_0$, we have $\mfr{f}_0 = b\sigma\mfr{f}_0$ and $V=w\sigma(V)$. Moreover, $\mu \in V$ since $o, b\sigma o \in \mfr{f}_0$. It follows that $\lambda \in V$, so we have $\mfr{f}_0 = \mfr{f}$. 
\end{proof}

\begin{cor}
    The facet $\mfr{f}$ is a minimal $b\sigma$-facet when $\breve{G}$ is simple, or more generally, when every $F$-simple factor of $G$ is $\breve{F}$-simple. 
\end{cor}

This minimality may fail if $G$ is $F$-simple but not $\breve{F}$-simple. 

\begin{exa}
    Suppose that $G$ is adjoint and $F$-simple. We may write 
    \[
        \breve{G} = \breve{G}_0 \times \sigma(\breve{G}_0) \times \cdots \times \sigma^{d-1}(\breve{G}_0)
    \]
    for some $d\geq 1$ and a simple adjoint group $\breve{G}_0$. Let
    \[
        \msf{G} = \msf{G}_0 \times \sigma(\msf{G}_0) \times \cdots \times \sigma^{d-1}(\msf{G}_0)
    \]
    be the corresponding decomposition. Let $\mu_i$ be the projection of $(w\sigma)^i\mu$ to $\msf{G}_0$ for $i \geq 0$. Via the identification $\cl{A}(\breve{G}, \breve{T})^{b\sigma} = \cl{A}(\breve{G}_0, \breve{T}_0)^{(b\sigma)^d}$, $\mfr{f}$ corresponds to a simplex spanned by 
    \[
        -\mu_0 - \mu_1 - \cdots - \mu_i
    \]
    for every $i \geq 0$. It contains a minimal $(b\sigma)^d$-stable simplex spanned by
    \[
        - \mu_0 - \mu_1 - \cdots - \mu_{di-1}
    \]
    for every $i \geq 0$, which can be strictly smaller than $\mfr{f}$. Note that $\mfr{f}$ is minimal if at most one of $\mu_0,\ldots,\mu_{d-1}$ is nontrivial. 
\end{exa}

\begin{lem} \label{lem:fparah}
    Let $G_{b,\mfr{f}} \subset G_b(F)$ be the parahoric subgroup associated to $\mfr{f}$. Then, we have
    \[
        G_{b, \mfr{f}} = G_b(F) \cap \cl{G}(O_{\breve{F}}). 
    \]
\end{lem}
\begin{proof}
    Since $o \in \mfr{f}$, $G_{b,\mfr{f}} \subset G_b(F) \cap \cl{G}(O_{\breve{F}})$. For the converse, let $g\in G_b(F) \cap \cl{G}(O_{\breve{F}})$. %Since $g (b\sigma)^k = (b\sigma)^k g$ and $g\cdot o = o$, we have 
    Then, 
    \[
        g\cdot (b\sigma)^ko = (b\sigma)^k \cdot g\cdot o = (b\sigma)^ko. 
    \]
    By \Cref{lem:charf}, $g$ fixes $\mfr{f}$. Since $g \in \cl{G}(O_{\breve{F}})$, the Kottwitz value of $g$ is trivial and $g\in G_{b,\mfr{f}}$. 
\end{proof}

In particular, $G_{b,\mfr{f}}$ is the stabilizer of $[1] \in X_\mu(b)$ and it acts on $\Spf(R_{\cl{G},\mu})$. Since level structures are insensitive to quasi-isogenies by construction, it even acts on $\mfr{X}_{b}^{\msf{\ov{G}}}$. 

\begin{lem}
    The composition
    \[
        G_{b,\mfr{f}} \to \cl{G}(O_{\breve{F}}) \to \msf{G}(\ov{k})
    \]
    factors through $\msf{\ov{P}}(\ov{k})$. Moreover, the image of
    \[
        G_{b,\mfr{f}} \to \msf{\ov{P}}(\ov{k}) \to \msf{M}(\ov{k})
    \]
    equals $\msf{M}^{w\sigma}$. In particular, the maximal reductive quotient of $G_{b, \mfr{f}}$ is $\msf{M}^{w\sigma}$.
\end{lem}
\begin{proof}
    Since $\mfr{f}=\mfr{a} \cap \tau$, the first claim follows from \cite[Theorem 8.4.19 (2)]{KP23}. Let $g\in G_{b, \mfr{f}}$ and let $p \in \ov{\msf{P}}(\ov{k})$ be its reduction with a Levi part $m \in \msf{M}(\ov{k})$. Since $g\in G_b(F)$, we have
    \[
        \Ad(w\sigma)(g) = \Ad(\mu(-\pi)^{-1})(g). 
    \]
    The reduction of the left-hand side is $\Ad(w\sigma)(p)$, so $g$ lies in the big cell $\Ad(w\sigma)(\msf{N}) \times \msf{M} \times \Ad(w\sigma)(\msf{\ov{N}}) \subset \msf{G}$. By comparing the reductions of the $\msf{M}$-part of both sides using \Cref{lem:MNprop}, we get $\Ad(w\sigma)(m) = m$. Thus, the image of $G_{b,\mfr{f}}\to \msf{M}(\ov{k})$ lies in $\msf{M}^{w\sigma}$. 
    
    Let $\cl{\breve{M}}$ be the Levi subgroup of $\cl{\breve{G}}$ corresponding to $\msf{M}$. Since $(w\sigma)^N = \sigma^N$ for some $N\geq 1$, $\breve{\cl{M}}^{w\sigma}$ is a reductive group scheme over $O_F$ and $\breve{\cl{M}}^{w\sigma}(k) = \msf{M}^{w\sigma}$. Since $\mu$ is central in $\cl{\breve{M}}$, $\breve{\cl{M}}^{w\sigma}(O_F) \subset G_{b,\mfr{f}}$, so $G_{b,\mfr{f}} \to \msf{M}^{w\sigma}$ is surjective. 
\end{proof}

\begin{prop} \label{prop:Gbact}
    The $G_{b,\mfr{f}}$-action on $\mfr{X}_{b,\red}^\perf \subset \msf{\ov{G}}^\perf$ is given by
    \[
        j \cdot \ov{g} = \ov{j} \ov{g} ,\quad j \in G_{b, \mfr{f}} ,\quad \ov{g} \in \ov{\msf{G}}. 
    \]
    Here, $\ov{j}$ is the image of $j$ in $\msf{G}(\ov{k})$. In particular, $U_b$ and $\cl{W}(0)$ are stable under $G_{b, \mfr{f}}$. The induced action on $\msf{W}(0)$ is given by the inflation of the $\msf{M}^{w\sigma}$-action, that is 
    \[
        j \cdot h = h \cdot m^{-1} \quad (h \in \msf{W}(0))
    \]
    where $m \in \msf{M}^{w\sigma}$ is the image of $j$. 
\end{prop}
\begin{proof}
    Let $x\in \mfr{X}_{b,\eta}(C,C^+)$ be a geometric point and choose $\pi^\flat$ and $u_\alpha^\flat$ as in \Cref{lem:expBKF}. Let $g\in \msf{G}(C^\flat)$ be the element associated to the level structure of $x$. Let $\ov{g}$ denote the specialization of $g$ in $\msf{\ov{G}}$. The specialization of $x$ is $\ov{g} \in \mfr{X}_{b,\red}^{\perf} \subset \msf{\ov{G}}^{\perf}$. 
    
    For each $j\in G_{b,\mfr{f}}$, $j\cdot x$ corresponds to a deformation
    \[
        (\cl{G}\otimes W_{O_F}(C^{\flat+}), j\mu([\pi^\flat]-\pi)\prod_{\alpha\in \Phi_{\mu<0}} i_\alpha([u_\alpha^\flat])w\sigma(j)^{-1}\cdot \sigma)
    \]
    by \Cref{prop:innactuniv}. Since $j$ preserves the underlying local $\cl{G}$-shtuka, the level structure of $j\cdot x$ is given by $\ov{j}g$ with $\ov{j}$ the image of $j$ in $\msf{G}(\ov{k})$. Take $v_\alpha^\flat \in C^{\flat+}$ so that $(v_\alpha^\flat)^\sharp$ equals the image of $u_\alpha$ under $R_{\cl{G},\mu}\xrightarrow{j} R_{\cl{G},\mu} \to C^+$. As in \Cref{rmk:changechoice}, there is $\iota \in \cl{G}(W_{O_F}(C^{\flat+}))$ that is trivial modulo $[\varpi]$ for some pseudo-uniformizer $\varpi \in C^{\flat+}$ and satisfies
    \[
        \iota j\mu([\pi^\flat]-\pi)\prod_{\alpha\in \Phi_{\mu<0}} i_\alpha([u_\alpha^\flat])w\sigma(\iota j)^{-1} =  \mu([\pi^\flat]-\pi)\prod_{\alpha\in \Phi_{\mu<0}} i_\alpha([v^\flat_\alpha])w.
    \]
    The prescribed depth-zero level structure for $j \cdot x$ is $\ov{\iota}\ov{j}g$. Since $\iota$ is trivial modulo $[\varpi]$, the specialization of $j\cdot x$ is $\ov{j} \ov{g} \in \mfr{X}_{b,\red}^{\perf} \subset \msf{\ov{G}}^{\perf}$. Thus, we get the first claim. 

    Let $m \in \msf{M}^{w\sigma}$ be the Levi part of $\ov{j}$. Since $\ov{j} \in \msf{\ov{P}}(\ov{k})$, we have $\ov{j}\cdot x_{\lambda} = x_{\lambda}\cdot m$. Since $\msf{M}^{w\sigma}$ acts on $\msf{W}(0)$ by the right translation, the image of
    \[
        \msf{W}(0) \to \mfr{X}_{b, \red} \subset \msf{\ov{G}} ,\quad g \mapsto x_{\lambda} g^{-1}
    \]
    is stable under $G_{b, \mfr{f}}$. By \Cref{thm:role_of_integralmodel}, $\cl{W}(0)$ is also stable under $G_{b, \mfr{f}}$ by setting $\ov{\msf{G}}$ to a toroidal compactification. We will compute the induced action on $\msf{W}(0)$. 

    Suppose that $x$ comes from a geometric point of $\mfr{W}(0)_\eta$. Then, we can write $g=(\pi^{\flat,1/e})^{e\lambda}h^{-1}$ and the specialization of $x$ in $\mfr{W}(0)$ is given by $[h] \in \msf{W}(0)$. Since $j\cdot x \in \mfr{W}(0)_\eta$, the level structure of $j \cdot x$ is also written as $(\pi^{\flat,1/e})^{e\lambda}h'^{-1}$ for $h' \in \msf{G}(C^{\flat+})$ and we have
    \[
        (\pi^{\flat,1/e})^{e\lambda}h'^{-1} = \ov{\iota} \ov{j}(\pi^{\flat,1/e})^{e\lambda}h^{-1} =
        \ov{\iota} (\ov{j}m^{-1})(\pi^{\flat,1/e})^{e\lambda}(hm^{-1})^{-1}.
    \] 
    Since $\iota$ is trivial modulo $[\varpi]$ and $\ov{j}m^{-1} \in \ov{\msf{N}}$, we have 
    \[
        hm^{-1} \equiv h' \in (\msf{G}/\msf{N})(C^{\flat+}/\pi^{\flat, 1/p^N})
    \]
    for sufficiently large $N$. In particular, the specialization of $j\cdot x$ in $\mfr{W}(0)$ is $[h]\cdot m^{-1} \in \msf{W}(0)$. 
\end{proof}

\subsection{The Galois action} \label{ssec:Weil}
In this section, we study the Galois action on $\mfr{W}(0)$. 

First, we study the inertial action on $\mfr{W}(0)$. Let $I_e \subset I_F$ be the Galois group of $\breve{F}_e$. Fix a primitive $e$-th root of unity $\zeta_e \in \breve{F}$ as $e$ is prime to $p$. This choice induces an isomorphism 
\[
    I_F/I_e \cong \bb{Z}/e\bb{Z}, \quad \tau(\pi^{1/e}) = \zeta_{e}^\tau \pi^{1/e}. 
\]

\begin{prop} \label{prop:inertact}
    The right action of $I_F/I_e$ on $\mfr{X}_{b,\eta} \otimes \breve{F}_e$ stabilizes $\mfr{W}(0)_\eta$. The induced action on $\msf{W}(0)$ is 
    \[
        g^\tau = g \cdot (\zeta_e^\tau)^{e\lambda} ,\quad \tau \in I_F/I_e \cong \bb{Z}/e\bb{Z} ,\quad g \in \msf{W}(0). 
    \]
    %For each $\tau \in I_F/I_e \cong \bb{Z}/e\bb{Z}$ and $g \in \msf{W}(0)$, we have $g^\tau = g \cdot (\zeta_e^\tau)^{e\lambda}$. 
\end{prop}
\begin{proof}
    The first claim follows since $\mfr{W}(0)_\eta$ is the base change of $\cl{W}(0)$. We compute the action on $\msf{W}(0)$ by looking at the specialization map. 

    Let $x\in \mfr{W}(0)_{\eta}(C,C^+)$ be a geometric point and choose $\pi^\flat, u_\alpha^\flat\in C^\flat$ as in \Cref{lem:expBKF}. Let $g=(\pi^{\flat, 1/e})^{e\lambda}h^{-1}$ be the element representing the depth-zero level structure at $x$. Since the inertial action only changes the choice of $\pi^{1/e}$ to $\tau(\pi^{1/e})$ preserving $g$, 
    \[
        g = \tau(\pi^{1/e})^{e\lambda}(h(\zeta_e^\tau)^{e\lambda})^{-1}
    \]
    implies that the specialization of $x^\tau$ is $[h]\cdot (\zeta_e^\tau)^{e\lambda} \in \msf{W}(0)$. Thus, the claim follows. 
\end{proof}
\begin{rmk}
    Since $\lambda - qw\sigma \lambda = \mu$, we have 
    \[
        (\zeta_e^\tau)^{e\lambda}\cdot \Ad(w\sigma)((\zeta_e^\tau)^{-e\lambda}) = (\zeta_e^\tau)^{e\mu} = 1 \in \msf{T}.
    \]
    Thus, $(\zeta_e^\tau)^{e\lambda} \in \msf{Z}_{\msf{M}^{w\sigma}}$ and the \textit{left} action of $I_F$ on $\msf{W}(0)$ is induced from that of $\msf{Z}_{\msf{M}^{w\sigma}}$ via
    \[
        I_F \twoheadrightarrow I_F / I_e \to \msf{Z}_{\msf{M}^{w\sigma}} ,\quad \tau \mapsto (\zeta_e^\tau)^{e\lambda}. 
    \]
\end{rmk}

Next, we study a Weil descent datum. Let us review the construction. 

\begin{const}\textup{(cf.\ \cite[(3.1.6)]{PR24})}
    There is a map
    \begin{equation} \label{eq:first_WeilDescent}
        \cl{M}^\ints_{\cl{G},b,\mu} \to \sigma^*\cl{M}^\ints_{\cl{G},b, \sigma\mu}
    \end{equation}
    sending a local $\cl{G}$-shtuka $\cl{P}$ with a quasi-isogeny $\iota\colon \cl{P}\vert_{\cl{Y}_{[r,\infty)}} \cong \cl{E}^b\vert_{\cl{Y}_{[r,\infty)}}$ to a local $\cl{G}$-shtuka $\sigma^*\cl{P}$ with a quasi-isogeny
    \[
        \sigma^*\cl{P}\vert_{\cl{Y}_{[r/q,\infty)}} \xrightarrow{\sigma^*\iota} \cl{E}^{\sigma(b)}\vert_{\cl{Y}_{[r/q,\infty)}}\xrightarrow{b} \cl{E}^{b}\vert_{\cl{Y}_{[r/q,\infty)}}. 
    \]
    The induced map on the underlying space is 
    \[
        X_\mu(b) \to \sigma^* X_{\sigma \mu}(b) ,\quad g \mapsto b\sigma(g). 
    \]
    Let $E$ be the reflex field of $(G, \mu)$ and let $f = [E \colon F]$. The Weil descent datum of $\cl{M}^\ints_{\cl{G},b,\mu}$ is the $f$-fold iteration of \eqref{eq:first_WeilDescent}
    \[
        \cl{M}^\ints_{\cl{G},b,\mu} \to (\sigma^f)^*\cl{M}^\ints_{\cl{G},b, \mu}. 
    \]
\end{const}

By passing to tubular neighborhoods, we get a map of deformation spaces
\[
    (\cl{M}^\ints_{\cl{G},b,\mu})^\wedge_{/[g]} \to \sigma^*(\cl{M}^\ints_{\cl{G},b, \sigma\mu})^\wedge_{/[b\sigma(g)]}. 
\]
Let $R$ be a perfectoid $O_{\breve{F}}$-algebra and let $R_1=W_{O_F}(R^\flat) / \sigma(\xi_R)W_{O_F}(R^\flat)$. The map of deformations of prismatic $(\cl{G},\mu)$-displays is given by
\[
    (\cl{G} \otimes W_{O_F}(R^\flat),X\sigma) \mapsto (\cl{G} \otimes W_{O_F}(R_1^\flat),\sigma(X)\cdot \sigma). 
\]
For each $f\geq 0$, let 
\[  
    b_f = b \sigma(b) \cdots \sigma^{f-1}(b) ,\quad
    R_f=W_{O_F}(R^\flat) / \sigma^{f}(\xi_R)W_{O_F}(R^\flat). 
\]
The Weil descent datum of $\cl{M}^\ints_{\cl{G},b,\mu}$ restricts to tubular neighborhoods
\[
    (\cl{M}^\ints_{\cl{G},b,\mu})^\wedge_{/[1]} \to (\sigma^f)^*(\cl{M}^\ints_{\cl{G},b,\mu})^\wedge_{/[b_f]}. 
\]
\begin{defi}
    Let $f_0$ be the minimum multiple of $f$ such that $[b_{f_0}]\in G_b(F) \cdot [1]$. 
\end{defi}

Then, we have an $f_0$-fold Weil descent datum on $\coprod_{j \in G_b(F)/G_{b,\mfr{f}}} (\cl{M}^\ints_{\cl{G},b,\mu})^\wedge_{/[j]}$. Let
\[
    b_{f_0} = j_0 b_0 ,\quad j_0 \in G_b(F) ,\quad b_0 \in \cl{G}(O_{\breve{F}}). 
\]
Then, we get an $f_0$-fold Weil descent datum 
\begin{equation}
    (\cl{M}^\ints_{\cl{G},b,\mu})^\wedge_{/[1]} \to (\sigma^{f_0})^*(\cl{M}^\ints_{\cl{G},b,\mu})^\wedge_{/[b_{f_0}]} \xrightarrow{j_0^{-1}}   (\sigma^{f_0})^*(\cl{M}^\ints_{\cl{G},b,\mu})^\wedge_{/[1]}. \label{eq:Weildesc}
\end{equation}
The map of deformations of prismatic $(\cl{G},\mu)$-displays is given by
\begin{equation} \label{eq:BKF_WeilDescent}
    (\cl{G} \otimes W_{O_F}(R^\flat),X\sigma) \mapsto (\cl{G} \otimes W_{O_F}(R_{f_0}^\flat),b_0\sigma^{f_0}(X)\sigma(b_0)^{-1} \cdot \sigma). 
\end{equation}
Note that we have $b_0\sigma^{f_0}(b)\sigma(b_0)^{-1} = b$ since $b_{f_0}\sigma^{f_0}\cdot b\sigma = b\sigma \cdot b_{f_0}\sigma^{f_0}$.

\begin{lem} \label{lem:underlying_Weil_descent}
    Let $\msf{\ov{G}}$ be a toroidal compactification associated to a $\sigma$-stable and $W$-stable complete fan in $X_*(\msf{T})$. The Weil descent datum \eqref{eq:Weildesc} extends to a Weil descent datum 
    \[
        \sigma^{f_0} \colon \mfr{X}_b \to (\sigma^{f_0})^*\mfr{X}_b, 
    \]
    and the induced Weil descent datum on the underlying space is
    \[
        \mfr{X}_{b,\red}^\perf \to (\sigma^{f_0})^*\mfr{X}_{b,\red}^\perf ,\quad 
        \ov{g} \mapsto b_0\sigma^{f_0}(\ov{g}) ,\quad \ov{g} \in \msf{\ov{G}}. 
    \]
\end{lem}
\begin{proof}
    The $\sigma$-action on $\msf{G}$ extends to $\ov{\msf{G}}$, so the first claim follows from the functoriality of the construction of $\mfr{X}_b$. By \eqref{eq:BKF_WeilDescent}, the Weil descent datum
    \[
        \mfr{X}_{b ,\eta} \to (\sigma^{f_0})^*\mfr{X}_{b , \eta} 
    \]
    sends a depth-zero level structure $g_x\in \msf{G}(C^\flat)$ to $b_0\sigma^{f_0}(g_x)$ for every geometric point $x\in \mfr{X}_{b,\eta}(C,C^+)$. By passing to the specialization, we get the second claim. 
\end{proof}

From now on, we choose $\msf{\ov{G}}$ as in \Cref{lem:underlying_Weil_descent}. 

\begin{defi}
    For each multiple $f$ of $f_0$, let
    \[  
        \sigma^f = (\sigma^{f_0})^{f / f_0} \colon \mfr{X}_b \to (\sigma^f)^* \mfr{X}_b
    \]
    be the $f$-fold Weil descent datum of $\mfr{X}_b$. This is obtained as in \eqref{eq:Weildesc} from the decomposition $ b_f = j_{0, f} b_{0, f}$ where
    \[
        b_{0, f} = b_0 \sigma^{f_0}(b_0) \cdots \sigma^{f - f_0}(b_0) \in \cl{G}(O_{\breve{F}}) ,\quad 
        j_{0, f} = j_0^{f/f_0} \in G_b(F)
    \]
\end{defi}

\begin{prop} \label{prop:Weilf1}
    Let $f_1$ be the minimum multiple of $f_0$ such that $U_b \cap \sigma^{f_1}(U_b) \neq \emptyset$. Then, we have $U_b = \sigma^{f_1}(U_b)$. In particular, 
    \[
        \coprod_{j \in G_b(F)/G_{b,\mfr{f}}} j \cdot \cl{W}(0) \hookrightarrow \cl{M}_{G,b,\mu, \cl{G}(1)}
    \]
    is an open immersion stable under the $f_1$-fold Weil descent datum and 
    \[
        \coprod_{\tau \in \langle \sigma^{f} \rangle/\langle \sigma^{f_1} \rangle} \coprod_{j \in G_b(F)/G_{b,\mfr{f}}} (\tau \times j )\cdot \cl{W}(0) \hookrightarrow \cl{M}_{G,b,\mu, \cl{G}(1)}
    \]
    is an open immersion. 
\end{prop}

\begin{proof}
    By construction, 
    \[
        U_b, \sigma^{f_0}(U_b),\ldots, \sigma^{f_1-f_0}(U_b) \subset \mfr{X}_{b, \red}
    \]
    are disjoint. Since each of them is open, $f_1$ is finite and well-defined. It is enough to show $U_b = \sigma^{f_1}(U_b)$, as the remaining claims follow from \Cref{lem:fparah} and \Cref{thm:role_of_integralmodel}. 

    Since $f_1$ is a multiple of $f$, $\sigma^{f_1}(\mu)$ is conjugate to $\mu$. Since both of them are dominant, we have $\sigma^{f_1}(\mu)=\mu$. Moreover, as $U_b \cap \sigma^{f_1}(U_b) \neq \emptyset$, $x_\lambda$ and $x_{\sigma^{f_1}(\lambda)}$ are in the same $\msf{G}\times \msf{G}$-orbit. Since $\lambda$ and $\sigma^{f_1}(\lambda)$ are dominant by \Cref{lem:lambdadom}, $x_\lambda = x_{\sigma^{f_1}(\lambda)}$ by \cite[Proposition 6.2.3 (ii)]{BK05}. Let $b_1 = b_{0, f_1}$. For $g\in \msf{W}(0)$, we have 
    \[
        \sigma^{f_1}(x_\lambda g^{-1}) = b_1 x_\lambda \sigma^{f_1}(g)^{-1}
    \]
    by \Cref{lem:underlying_Weil_descent}. Suppose that $ b_1 x_\lambda \sigma^{f_1}(g)^{-1} = x_\lambda h^{-1}$ for $h \in \msf{W}(0)$. Then, $(b_1, \sigma^{f_1}(g)^{-1} h)$ lies in the stabilizer of $x_\lambda$, so the reduction of $b_1$ lies in $\ov{\msf{P}}$. Moreover, $b_1 \sigma^{f_1}(b)\sigma(b_1)^{-1} = b$ implies
    \[
        \mu(-\pi)^{-1} b_1\mu(-\pi) = w \sigma(b_1)\sigma^{f_1}(w)^{-1}. 
    \]
    Let $m\in \msf{M}$ be the Levi part of the reduction of $b_1$. By \Cref{lem:MNprop}, the reduction of the above equation implies
    \[
        w \sigma(m) \sigma^{f_1}(w)^{-1} \in m\cdot \Ad(w\sigma)(\ov{\msf{N}}). 
    \]
    Since $\msf{M}$ is stable under $w\sigma$, the Bruhat cell $\Ad(w\sigma)(\ov{\msf{P}}) \cdot w \sigma^{f_1}(w)^{-1}\cdot \Ad(w\sigma)(\ov{\msf{P}})$ contains the identity. By considering the Bruhat decomposition for $\Ad(w\sigma)(\ov{\msf{P}}) \backslash \msf{G} / \Ad(w\sigma)(\ov{\msf{P}})$, we get $w \sigma^{f_1}(w)^{-1} \in \msf{M}$. Thus, we have
    \[
        w \sigma(m) \sigma^{f_1}(w)^{-1} = m. 
    \]
    Now, for every $g \in \msf{W}(0)$, let $h=\sigma^{f_1}(g)m^{-1}$. Then, $b_1x_\lambda\sigma^{f_1}(g^{-1}) = x_\lambda h^{-1}$ and we have 
    \begin{align*}
        h^{-1}\sigma(h) &\in m \cdot \sigma^{f_1}(\msf{U}_{\mu<0}\cdot w)\cdot \sigma(m)^{-1} \\
        &= \msf{U}_{\mu<0} \cdot m \sigma^{f_1}(w) \sigma(m)^{-1}\\
        &= \msf{U}_{\mu<0} \cdot w. 
    \end{align*}
    Thus, $h \in \msf{W}(0)$, so we have $U_b = \sigma^{f_1}(U_b)$. 
\end{proof}

\begin{cor} \label{cor:W0_WeilDescent}
    The $f_1$-fold Weil descent datum $\sigma^{f_1}$ induces the one for $\mfr{W}(0)$. Let $m\in \msf{M}$ be the Levi part of the reduction of $b_{0, f_1}$. Then, the induced Weil descent datum on $\msf{W}(0)$ is given by
    \[
        \msf{W}(0) \to (\sigma^{f_1})^* \msf{W}(0) ,\quad g\mapsto \sigma^{f_1}(g)m^{-1}. 
    \] 
\end{cor}
\begin{proof}
    Let $x \in \mfr{W}(0)_\eta(C,C^+)$ be a geometric point and choose $\pi^\flat$ and $u_\alpha^\flat$ as in \Cref{lem:expBKF}. Let $b_1 = b_{0, f_1}$. Then, $\sigma^{f_1}(x)$ corresponds to a deformation
    \[
        (\cl{G}\otimes W_{O_F}(C_{f_1}^{\flat+}), b_1 \mu([\pi^{\flat, p^{ f_1}}]-\pi)\prod_{\alpha\in \Phi_{\mu<0}} i_\alpha([u_\alpha^{\flat, p^{f_1}}])w\sigma(b_1)^{-1}\cdot \sigma)
    \]
    by \eqref{eq:BKF_WeilDescent}. Let $g=(\pi^{\flat,1/e})^{e\lambda}h^{-1}\in \msf{G}(C^{\flat})$ be the depth-zero level structure of $x$. % Then, the one for $\sigma^{f_1}(x)$ is given by $b_1^{-1} (\pi^{\flat,1/e})^{e\lambda}h^{-1}$. 
    
    As in the proof of \Cref{prop:Gbact}, take $v_\alpha^\flat \in C^{\flat+}$ and $\iota \in \cl{G}(W_{O_F}(C^{\flat+}))$ that is trivial modulo $[\varpi]$ for some pseudo-uniformizer $\varpi \in C^{\flat+}$ so that 
    \[
        (\iota b_1 )\mu([\pi^{\flat, p^{f_1}}]-\pi)\prod_{\alpha\in \Phi_{\mu<0}} i_\alpha([u_\alpha^{\flat, p^{f_1}}])w\sigma(\iota b_1)^{-1} =  \mu([\pi^{\flat, p^{f_1}}]-\pi)\prod_{\alpha\in \Phi_{\mu<0}} i_\alpha([v^\flat_\alpha])w.
    \]
    The depth-zero level structure for $\sigma^{f_1}(x)$ is $\ov{\iota} b_1 \sigma^{f_1}(g)$. As $\sigma^{f_1}(x) \in \mfr{W}(0)_\eta$, we can write
    \[
        \ov{\iota}b_1 \sigma^{f_1}(g) = (\pi^{\flat,p^{f_1}/e})^{e\lambda}h'^{-1}. 
    \]
    Since the reduction of $b_1$ lies in $\ov{\msf{P}}$ (see the proof of \Cref{prop:Weilf1}), we have
    \[
        \sigma^{f_1}(h) m^{-1} \equiv h' \in (\msf{G}/\msf{N})(C^{\flat+}/\pi^{\flat, 1/p^N})
    \]
    for sufficiently large $N$. Thus, the specialization of $\sigma^{f_1}(x)$ in $\mfr{W}(0)$ is $\sigma^{f_1}(h) m^{-1} \in \msf{W}(0)$. 
\end{proof}

When $G$ is split, the Weil descent datum becomes quite simple. 

\begin{prop} \label{prop:split_WeilDescent}
    When $G$ is split, $E = F$ and $\mu$ can be defined over $O_F$. Then, $\cl{U}_{\mu < 0}$ is defined over $O_F$ and we can take $b \in G(F)$. Then, we can take
    \[
        b_0 = 1 ,\quad f_0 = f_1 = 1. 
    \]
    Let $O_{F_e} = O_F[\pi^{1/e}]$. Then, $\mfr{W}(0)$ can be defined over $O_{F_e}$ and the induced Weil descent datum on $\mfr{W}(0)$ equals the one introduced in \Cref{cor:W0_WeilDescent}. 
\end{prop}
\begin{proof}
    The claim is immediate from the construction. In this case, $R_{\cl{G}, \mu}$ (resp.\ $R_{\cl{G}, \mu, \lambda}$) is defined over $O_F$ (resp.\ $O_{F_e}$). Moreover, the Lang torsor $\msf{W}(0) \to \msf{U}_{\mu < 0}$ is defined over $k$. Thus, $\mfr{W}(0)$ can be defined over $O_{F_e}$. 
\end{proof}

\section{Nearby cycles of special affinoids} \label{sec:nearby}

In this section, we explain a general strategy to extract subrepresentations of the middle \'{e}tale cohomology of local Shimura varieties from special affinoids $\cl{W}(0)$. 

Let $\ell \neq p$ be another prime and let $\Lambda$ be a $\bb{Z}_{\ell}$-algebra such that
\[
    \Lambda = \bb{Z}/\ell^m\bb{Z} \quad (m\geq 1) \quad \text{or} \quad \Lambda = \Qla. 
\]
Let $\bb{C}_p$ be the completion of an algebraic closure of $\breve{F}_e$ and let $O_{\bb{C}_p}$ be its ring of integers. Moreover, $H^i$ denotes the $i$-th \'{e}tale cohomology group of schemes or adic spaces, and $R\Gamma$ denotes the \'{e}tale cohomology complex. We use $H_c^i$ and $R\Gamma_c$ for compactly supported variants. 

\subsection{Specialization of the \'{e}tale cohomology}

In this section, we review the specialization map of the \'{e}tale cohomology for smooth affine formal schemes over $O_{\bb{C}_p}$, following \cite[Section 4]{Tsu16} and \cite[Section 2]{Mie16}. 

First, suppose $\Lambda = \bb{Z}/\ell^m$ for some $m \geq 1$. Let $\mfr{X}$ be a $p$-adic smooth affine formal scheme over $O_{\bb{C}_p}$. Then, we have a morphism of \'{e}tale sites
\[
    \lambda_{\mfr{X}} \colon \mfr{X}_{\eta, \et} \to \mfr{X}_{\red, \et}. 
\]
The derived pushforward along $\lambda_{\mfr{X}}$ is denoted by $R\Psi_{\mfr{X}}$. By adjunction, we have the unit
\[
    \spc^* \colon \Lambda \to R\Psi_{\mfr{X}} \Lambda. 
\]

\begin{prop}\textup{(\cite[Corollary 4.29]{Mie14})}
    There is a functorial isomorphism
    \[
        \xi_{\mfr{X}} \colon R\Gamma_{c}(\mfr{X}_\red, R\Psi_{\mfr{X}}\Lambda) \cong R\Gamma_{c}(\mfr{X}_\eta, \Lambda). 
    \]
\end{prop}

This isomorphism enables us to define the specialization map for \textit{compactly supported} \'{e}tale cohomology and the \textit{cospecialization} map for \'{e}tale cohomology. 

\begin{defi}
    The composition with $\spc^*$ induces a specialization map
    \[
        R\Gamma_{c}(\mfr{X}_\red, \Lambda) \to R\Gamma_{c}(\mfr{X}_\red, R\Psi_{\mfr{X}}\Lambda) \cong R\Gamma_{c}(\mfr{X}_\eta, \Lambda), 
    \]
    which is also denoted by $\spc^*$ by abuse of notation. 
\end{defi}

\begin{defi}
    Let $d$ be the dimension of $\mfr{X}$. The trace map
    \[
        R\Gamma_c(\mfr{X}_\red, R\Psi_{\mfr{X}}\Lambda) \cong R\Gamma_c(\mfr{X}_\eta, \Lambda) \to \Lambda(-d)[-2d]
    \]
    induces a map
    \[
        \cosp^* \colon R\Psi_{\mfr{X}} \Lambda \to \pi_\red^! \Lambda(-d)[-2d] \cong \Lambda, 
    \]
    where $\pi_\red \colon \mfr{X}_\red \to \Spec(\ov{k})$ is the structure map. The induced map
    \[
        R\Gamma(\mfr{X}_\eta, \Lambda) \cong R\Gamma(\mfr{X}_\red, R\Psi_{\mfr{X}} \Lambda) \to R\Gamma(\mfr{X}_\red, \Lambda)
    \]
    is the cospecialization map, which is also denoted by $\cosp^*$ by abuse of notation. 
\end{defi}

\begin{prop}\textup{(\cite[Corollary 2.7]{Mie16})} \label{prop:sp_cosp_commute}
    The diagram
    \begin{center}
        \begin{tikzcd}
            H_c^i(\mfr{X}_\red, \Lambda) \ar[r, "\spc^*"] \ar[d] &  H_c^i(\mfr{X}_\eta, \Lambda) \ar[d] \\
            H^i(\mfr{X}_\red, \Lambda) &  H^i(\mfr{X}_\eta, \Lambda) \ar[l, "\cosp^*"] 
        \end{tikzcd}
    \end{center}
    commutes and can be defined even when $\Lambda = \Qla$ by taking projective limits. 
\end{prop}

\begin{rmk}
    Due to the Artin vanishing, this diagram is meaningful only for $i = d$. 
\end{rmk}

\subsection{Computation at depth zero}

In this section, we explain how to extract a subrepresentation of the middle \'{e}tale cohomology
\[
    H_c^r(\cl{M}_{G,b,\mu, \cl{G}(1), \bb{C}_p}, \Lambda). 
\]
Recall $\dim \msf{W}(0) = \dim \cl{M}_{G, b, \mu, \cl{G}(1)} = r$. By \Cref{prop:sp_cosp_commute}, we have a commutative diagram
\begin{center}
    \begin{tikzcd}
        H_c^r(\msf{W}(0),\Lambda) \ar[r, "\spc^*"] \ar[d] & H_c^r(\cl{W}(0)_{\bb{C}_p}, \Lambda) \ar[r] \ar[d] & H_c^r(\cl{M}_{G,b,\mu, \cl{G}(1), \bb{C}_p}, \Lambda) \ar[d] \\
        H^r(\msf{W}(0),\Lambda) & H^r(\cl{W}(0)_{\bb{C}_p}, \Lambda) \ar[l, "\cosp^*"] &
        H^r(\cl{M}_{G,b,\mu, \cl{G}(1),\bb{C}_p}, \Lambda). \ar[l]
    \end{tikzcd}
\end{center}
The second square is obtained formally from the inclusion. 
As a corollary, we get the following. 

\begin{cor}\textup{(\cite[Theorem 2.8]{Mie16})}
    If a submodule $V \subset H_c^r(\msf{W}(0), \Lambda)$ maps injectively to $H^r(\msf{W}(0),\Lambda)$, the composition
    \[
        V \hookrightarrow H_c^r(\msf{W}(0),\Lambda) \to H_c^r(\cl{M}_{G,b,\mu, \cl{G}(1), \bb{C}_p}, \Lambda)
    \]
    is injective. 
\end{cor}

The above diagram is equivariant under $\msf{G}^\sigma \times G_{b,\mfr{f}} \times I_F$ by the functoriality of the construction. Its action on $H_c^r(\msf{W}(0),\Lambda)$ is induced from the $\msf{G}^\sigma \times \msf{M}^{w\sigma}$-action on $\msf{W}(0)$ by \Cref{prop:Gbact} and \Cref{prop:inertact}. Moreover, by \Cref{prop:Weilf1}, 
\[
    \coprod_{j \in G_b(F)/G_{b,\mfr{f}}} j \cdot \cl{W}(0)_{\bb{C}_p} \hookrightarrow \cl{M}_{G, b,\mu, \cl{G}(1), \bb{C}_p}
\]
is an open subspace stable under the $f_1$-fold Weil descent datum. Thus, it induces a map
\[
    \cInd_{G_{b,\mfr{f}}}^{G_b(F)} H_c^r(\msf{W}(0), \Lambda) \to H_c^r(\cl{M}_{G, b,\mu, \cl{G}(1), \bb{C}_p}, \Lambda)
\]
of smooth $\msf{G}^\sigma \times G_b(F) \times W_{E_1}$-representations where $E_1$ is the unramified extension of $F$ of degree $f_1$. By abuse of notation, let $\sigma$ denote an arithmetic Frobenius element of $W_{F[\pi^{1/e}]}$. Applying \eqref{eq:Weildesc} to $f_1$ instead of $f_0$, the action of $\sigma^{f_1}$ on $H_c^r(\msf{W}(0),\Lambda) \subset \cInd_{G_{b,\mfr{f}}}^{G_b(F)} H_c^r(\msf{W}(0), \Lambda)$ is given by 
\[
    H_c^r(\msf{W}(0),\Lambda)\xrightarrow{\sigma^{f_1}_*} H_c^r(\msf{W}(0),\Lambda) \xrightarrow{j_{0, f_1}^{-1}} \cInd_{G_{b,\mfr{f}}}^{G_b(F)} H_c^r(\msf{W}(0), \Lambda).
\] 
Here, $\sigma^{f_1}_*$ denotes the pushforward along the $f_1$-fold Weil descent datum in \Cref{cor:W0_WeilDescent}. By the Frobenius reciprocity, we get a map
\[
    \Ind^{W_{E}}_{W_{E_{1}}} \cInd_{G_{b,\mfr{f}}}^{G_b(F)} H_c^r(\msf{W}(0), \Lambda) \to H_c^r(\cl{M}_{G, b, \mu, \cl{G}(1), \bb{C}_p}, \Lambda).
\]

\begin{prop} \label{lem:subrep}
    Let $V\subset H_c^r(\msf{W}(0),\Lambda)$ be a $\sigma^{f_1}_*$-stable $\msf{G}^\sigma \times \msf{M}^{w\sigma}$-subrepresentation such that $V\to H^r(\msf{W}(0),\Lambda)$ is injective. Then, $\Ind^{W_{E}}_{W_{E_1}} \cInd_{G_{b,\mfr{f}}}^{G_b(F)} V \to H_c^r(\cl{M}_{G,b,\mu,\cl{G}(1), \bb{C}_p}, \Lambda)$ is injective. More precisely, we have the following commutative diagram. 
    \begin{center}
        \begin{tikzcd}
            \Ind^{W_{E}}_{W_{E_1}} \cInd_{G_{b,\mfr{f}}}^{G_b(F)} V \ar[r,hook] \ar[rd, hook] & H_c^r(\cl{M}_{G,b,\mu,\cl{G}(1), \bb{C}_p}, \Lambda) \ar[d] \\
            & \Ind^{W_{E}}_{W_{E_1}} \Ind_{G_{b,\mfr{f}}}^{G_b(F)} H^r(\msf{W}(0),\Lambda). 
        \end{tikzcd}
    \end{center}
\end{prop}
\begin{proof}
    The map
    \[
        H_c^r(\cl{M}_{G,b,\mu, \cl{G}(1), \bb{C}_p}, \Lambda) \to \Ind^{W_{E}}_{W_{E_1}} \Ind_{G_{b,\mfr{f}}}^{G_b(F)} H^r(\msf{W}(0),\Lambda)
    \]
    can be constructed by repeating the previous argument for $H^r(\msf{W}(0),\Lambda)$ instead of $H_c^r(\msf{W}(0),\Lambda)$. For the claim, it is enough to show that the composition
    \[
        \Ind^{W_{E}}_{W_{E_1}} \cInd_{G_{b,\mfr{f}}}^{G_b(F)} H_c^r(\msf{W}(0), \Lambda) \to H_c^r(\cl{M}_{G,b,\mu, \cl{G}(1), \bb{C}_p}, \Lambda) \to \Ind^{W_{E}}_{W_{E_1}} \Ind_{G_{b,\mfr{f}}}^{G_b(F)} H^r(\msf{W}(0),\Lambda)
    \]
    is induced from $ H_c^r(\msf{W}(0), \Lambda) \to H^r(\msf{W}(0),\Lambda)$. For any two distinct elements $s, t \in W_E / W_{E_1} \times G_b(F) / G_{b, \mfr{f}}$, we have 
    \[
        s \cdot \cl{W}(0)_{\bb{C}_p} \cap t \cdot \cl{W}(0)_{\bb{C}_p} = \emptyset
    \]
    by \Cref{prop:Weilf1}. Thus, the composition 
    \[
        s\cdot H_c^r(\msf{W}(0), \Lambda) \to  H_c^r(\cl{M}_{G,b,\mu, \cl{G}(1), \bb{C}_p}, \Lambda) \to t \cdot H^r(\msf{W}(0),\Lambda)
    \]
    is zero, so we get the claim. 
\end{proof}

Even if $V$ is not stable under $\sigma^{f_1}_*$, the same argument works by replacing $\Ind^{W_{E}}_{W_{E_1}} \cInd_{G_{b,\mfr{f}}}^{G_b(F)} V$ with $\cInd_{G_{b,\mfr{f}}}^{G_b(F)} V$. After fixing a central character, we get the following. 

\begin{prop} \label{prop:induction_with_central_character}
    Let $V\subset H_c^r(\msf{W}(0),\Qla)$ be a $\msf{G}^\sigma \times \msf{M}^{w\sigma}$-subrepresentation such that 
    \[
        V\to H^r(\msf{W}(0),\Qla)
    \]
    is injective and $V$ has a central character $\chi^0 \colon Z_G(F) \cap G_{b, \mfr{f}} \to \Qlax$. Let $\chi \colon Z_G(F) \to \Qlax$ be an extension of $\chi^0$. Then, we have a commutative diagram
    \begin{center}
        \begin{tikzcd}
            \cInd_{Z_G(F)G_{b,\mfr{f}}}^{G_b(F)} (\chi \boxtimes V) \ar[r, hook] \ar[rd, hook] & H_c^r(\cl{M}_{G,b,\mu, \cl{G}(1), \bb{C}_p}, \Qla) \otimes_{Z_G(F)} \chi \ar[d] \\
            & \Ind_{Z_G(F)G_{b,\mfr{f}}}^{G_b(F)} (\chi \boxtimes H^r(\msf{W}(0),\Qla)_{\chi^0}).
        \end{tikzcd}
    \end{center}
    Here, $H^r(\msf{W}(0),\Qla)_{\chi^0}$ denotes the coinvariant $H^r(\msf{W}(0),\Qla) \otimes_{Z_G(F) \cap G_{b,\mfr{f}}} \chi^0$. 
\end{prop}
\begin{proof}
    The sequence
    \[
        V \to H_c^r(\cl{M}_{G,b,\mu, \cl{G}(1), \bb{C}_p}, \Qla) \to \prod_{z \in Z_G(F) / Z_G(F) \cap G_{b, \mfr{f}}} z \cdot H^r(\msf{W}(0),\Qla)
    \]
    induces a sequence
    \[
        \chi \boxtimes V \to H_c^r(\cl{M}_{G,b,\mu, \cl{G}(1), \bb{C}_p}, \Qla) \otimes_{Z_G(F)} \chi \to \chi \boxtimes H^r(\msf{W}(0),\Qla)_{\chi^0}
    \]
    of $Z_G(F) G_{b, \mfr{f}}$-representations. Here, the composition is injective since $H^r(\msf{W}(0),\Qla)_{\chi^0}$ is the $\chi^0$-isotypic component of $H^r(\msf{W}(0),\Qla)$. It is enough to show that the composition
    \[
        \cInd_{Z_G(F)G_{b,\mfr{f}}}^{G_b(F)} (\chi \boxtimes V)  \to H_c^r(\cl{M}_{G,b,\mu, \cl{G}(1), \bb{C}_p}, \Qla) \otimes_{Z_G(F)} \chi \to \Ind_{Z_G(F)G_{b,\mfr{f}}}^{G_b(F)} (\chi \boxtimes H^r(\msf{W}(0),\Qla)_{\chi^0})
    \]
    is induced from the injection $\chi \boxtimes V \hookrightarrow \chi \boxtimes H^r(\msf{W}(0),\Qla)_{\chi^0}$. For any two distinct elements $s, t \in G_b(F) / Z_G(F) G_{b, \mfr{f}}$, we have 
    \[
        s Z_G(F) \cdot \cl{W}(0)_{\bb{C}_p} \cap t Z_G(F) \cdot \cl{W}(0)_{\bb{C}_p} = \emptyset
    \]
    by \Cref{prop:Weilf1}. Thus, the composition 
    \[
        s\cdot (\chi \boxtimes V) \to  H_c^r(\cl{M}_{G,b,\mu, \cl{G}(1), \bb{C}_p}, \Qla) \otimes_{Z_G(F)} \chi \to t \cdot (\chi \boxtimes H^r(\msf{W}(0),\Qla)_{\chi^0}). 
    \]
    is zero, so we get the claim. 
\end{proof}

\section{Contribution of depth-zero regular supercuspidal representations} \label{sec:dep0reg}

In this section, we construct explicit Jacquet-Langlands pairs of depth-zero regular supercuspidal representations in the cohomology of local Shimura varieties. 

\subsection{Depth-zero regular supercuspidal representations}

In this section, we review the construction of depth-zero regular supercuspidal representations in \cite{DR09}.  
First, recall the definition of unramified elliptic regular pairs. 

\begin{defi}
    Let $S \subset G$ be an unramified elliptic maximal torus with an integral model $\cl{S} \subset \cl{G}$ and a special fiber $\msf{S} \subset \msf{G}^\sigma$. Let $\theta \colon S(F) \to \Qlax$ be a depth-zero regular character, in the sense that 
    \begin{enumerate}
        \item $\theta$ is trivial on $\cl{S}(O_F) \cap \cl{G}(1)$, and 
        \item the induced character $\theta^0$ on $\msf{S}$ is in general position in the sense of \cite[Definition 5.15]{DL76}.
    \end{enumerate} 
    Then, we say that the pair $(S, \theta)$ is an unramified elliptic regular pair. 
\end{defi}

\begin{rmk}
    Let $\chi = \theta\vert_{Z_G(F)}$. Since $S(F) = Z_G(F) \cl{S}(O_F)$ when $S$ is elliptic and unramified, the pair $(S, \theta)$ is uniquely determined by the triple $(\msf{S}, \theta^0, \chi)$. 
\end{rmk}

Local class field theory attaches a tame $L$-parameter $\varphi_{\theta}\colon W_F \to {}^LS(\ov{\bb{Q}}_\ell)$ to each $(S, \theta)$. Since $S$ is unramified, there is a natural $L$-embedding 
\[
    {}^Lj \colon {}^LS \to {}^L G
\]
and the composition $\varphi_{(S, \theta)} = {}^Lj \circ \varphi_\theta$ is the tame regular semisimple elliptic $L$-parameter (TRSELP) attached to $(S, \theta)$. Here, ${}^Lj \vert_{\widehat{S} \rtimes I_F}$ is written as
\[
    \widehat{S}\rtimes I_F \cong \widehat{T} \rtimes I_F \subset \LG. 
\]
By abuse of notation, let $\sigma$ denote an arithmetic Frobenius element in $W_F$. Then, ${}^Lj(\sigma) = \widehat{v} \rtimes \sigma$ for some $\widehat{v} \in N_{\widehat{G}}(\widehat{T})$. Let us denote the corresponding element of the Weyl group by 
\[
    v \in W_o = N_{\breve{\cl{G}}}(\breve{\cl{T}}) / \breve{\cl{T}}. 
\]
Take $p_v \in \cl{G}(O_{\breve{F}})$ so that $\cl{S} = \Ad(p_v)(\cl{T})$. Then, $\dot{v} = p_v^{-1}\sigma(p_v) \in N_{\breve{\cl{G}}}(\breve{\cl{T}})$ is a lift of $v$. 
\begin{thm} \textup{(\cite{DR09}, \cite{Kal14})} \label{thm:Depth0_LLC}
    Let $S_{\varphi}$ be the centralizer of $\varphi_{(S, \theta)}$ in $\widehat{G}$. Then, $S_\varphi \cong \widehat{T}^{\widehat{v} \sigma = 1}$. In particular, $S_\varphi$ is abelian and the set of its algebraic irreducible representations is 
    \[
        \Irr(S_\varphi) = X^*(\widehat{T}^{\widehat{v}\sigma = \id}) \cong X_*(T)/(1 - v \sigma)X_*(T). 
    \]
    This set is bijective to the union of $L$-packets $\Pi_{\varphi_{(S, \theta)}}(G_b)$ over all $b\in B(G)_\bas$. 
\end{thm}

\begin{rmk}
    Note that \cite{DR09} parametrizes the union of $L$-packets over pure inner forms and \cite{Kal14} extends the parametrization to extended pure inner forms. 
\end{rmk}

Now, we review the construction of the depth-zero regular supercuspidal representation attached to $\nu \in X_*(T)$. Since $S$ is elliptic, $X_*(T/Z_G)^{v \sigma= \id} = 0$, so $\nu(\pi)v \sigma$ has a unique fixed point 
\[
    x_\nu \in \cl{A}(\breve{G}, \breve{T}).
\]
Let $J_\nu$ be the facet of $\cl{A}(\breve{G}, \breve{T})$ such that $x_\nu \in J_\nu^\circ$ and let $C_\nu$ be an alcove with $J_\nu \subset \ov{C_\nu}$. Since $J_\nu$ is fixed by $\nu(\pi)v \sigma$, we get an expression
\begin{equation} \label{eq:basic_decomposition}
    \nu(\pi)v \sigma = w_\nu y_\nu \sigma
\end{equation}
where $w_\nu \in W_a$ fixes $J_\nu$ and $y_\nu \in \wtd{W}$ is an element such that $y_\nu \sigma$ stabilizes $C_\nu$. Fix a decent lift $u_\nu \in N_G(T)(\breve{F})$ of $y_\nu \in \wtd{W}$ and let $\sigma_\nu = u_\nu \sigma$.

Let $G_\nu$ be the parahoric subgroup of $G(\breve{F})$ attached to $J_\nu$. We can take $p_\nu \in G_\nu$ so that 
\[
    p_\nu^{-1} \sigma_\nu(p_\nu) = \nu(\pi)\dot{v}u_\nu^{-1}
\]
(see \cite[p.814]{DR09}). Let $T_\nu = \Ad(p_\nu)(T)$. Then, $T_\nu$ is stable under $\sigma_\nu$ and there is a natural isomorphism $T_\nu^{\sigma_\nu} \cong S$. Then, $\theta$ induces a depth-zero character of $T_\nu^{\sigma_\nu}(F)$. 

By the Deligne-Lusztig induction from $\theta$, we get an irreducible representation
\[
    \kappa_\nu \in \Irr(Z_G(F)G_\nu^{\sigma_\nu}). 
\]
By taking the compact induction, we get the depth-zero regular supercuspidal representation 
\[
    \pi_\nu = \cInd_{Z_G(F)G_\nu^{\sigma_\nu}}^{G^{\sigma_\nu}(F)} \kappa_\nu
\]
attached to $\nu$. Note that $\pi_\nu$ depends only on the class $[\nu] \in X_*(T)/(1 - v \sigma)X_*(T)$. 

\begin{lem}\textup{(\cite[Lemma 4.5.1]{DR09})} \label{lem:cInd=Ind}
    For each $\nu$, we have $\cInd_{Z_G(F)G_\nu^{\sigma_\nu}}^{G^{\sigma_\nu}(F)} \kappa_\nu \cong \Ind_{Z_G(F)G_\nu^{\sigma_\nu}}^{G^{\sigma_\nu}(F)} \kappa_\nu$. 
\end{lem}

\subsubsection{Representations of interest}

We consider when $\pi_\nu$ appears in the contribution of our special affinoid $\cl{W}(0)$. We will choose $\nu$ so that 
\[
    J_\nu = \mfr{f} \subset \mfr{a} = C_\nu. 
\]
Then, $w_\nu \in W_o$ and \eqref{eq:basic_decomposition} implies that $(w_\nu^{-1}\nu)(\pi) \cdot w_\nu^{-1}v\sigma$ fixes $\mfr{a}$. In particular, $w_\nu^{-1} \nu$ is dominant. On the other hand, $\nu(\pi)v\sigma$ fixes $\mfr{f}$, so $\nu$ is also dominant, so $w_\nu^{-1}\nu = \nu$. As we see in \Cref{ssec:isotcoho}, $\nu$ is also required to be a weight of the highest representation $V_{\mu}$. Thus, we need 
\[
    \nu = \mu ,\quad w_\nu \mu = \mu. 
\]
Then, $\mu(\pi) \cdot w_\nu^{-1}v\sigma$ fixes $\mfr{a}$, so we have
\[
    v = w_\nu w ,\quad y_\nu = b \in \wtd{W}. 
\]
Since $w_\nu$ fixes $\mfr{f}$, $w_\nu \lambda = \lambda$. Let $W_{\msf{M}} \subset W_o$ be the Weyl group of $\msf{M}$. Then, $w_\nu \in W_{\msf{M}}$. 
In conclusion, the relevant case is 
\[
    \nu = \mu ,\quad v = w_{M} w \quad (w_{M} \in W_{\msf{M}}). 
\]
In this case, we have 
$
    \msf{S} \subset \msf{M}^{w\sigma} \subset \msf{G}^\sigma
$. 
Thus, we take $(S, \theta)$ under the following conditions. 

\begin{ass} \label{ass:triple}
    Let $(S, \theta)$ be an unramified elliptic regular pair such that $S \subset \breve{M}^{w\sigma}$. 
\end{ass}

\begin{lem}
    Under \Cref{ass:triple}, $\mfr{f}$ is a minimal $b\sigma$-facet, and $(\msf{S},\theta)$ is $(\msf{G}^\sigma, \msf{M}^{w\sigma})$-superregular in the sense of \textup{\cite[{}10.3]{BR03}}.  
\end{lem}
\begin{proof}
    Since $S$ is elliptic, $Z_{\breve{M}^{w\sigma}} / Z_G$ is anisotropic. As in the proof of \Cref{lem:charf}, $\mfr{f}$ is parallel to $X_*(\msf{Z}_{\msf{M}} / \msf{Z}_\msf{G}) \otimes \bb{R}$, so $\mfr{f}^{b\sigma}$ is a point. Thus, we get the first claim. The second claim is immediate from the definition since $(\msf{S}, \theta)$ is in general position in $\msf{G}^\sigma$. 
\end{proof}

\begin{defi}
    Let
    \[
        \pi^0_{(S,\theta)} = \varepsilon_{(\msf{G}^\sigma, \msf{S})} R_{\msf{S}}^{\msf{G}^\sigma}(\theta^0) ,\quad \rho^0_{(S,\theta)} = \varepsilon_{(\msf{M}^{w\sigma}, \msf{S})} R_{\msf{S}}^{\msf{M}^{w\sigma}}(\theta^0)
    \]
    be irreducible cuspidal representations obtained by the twisted Deligne-Lusztig induction from $\theta^0$. Here, $\varepsilon_{(\msf{G}^\sigma, \msf{S})}$ and $\varepsilon_{(\msf{M}^{w\sigma}, \msf{S})}$ are the signs in \cite[Theorem 8.3]{DL76}. 
\end{defi}

\begin{lem} \label{lem:middegrho}
    The natural map
    \[
        H_c^r(\msf{W}(0), \bb{\ov{Q}}_{\ell}) \otimes_{\bb{\ov{Q}}_{\ell}[\msf{M}^{w\sigma}]} \rho^0_{(S,\theta)} \to H^r(\msf{W}(0), \bb{\ov{Q}}_{\ell}) \otimes_{\bb{\ov{Q}}_{\ell}[\msf{M}^{w\sigma}]} \rho^0_{(S,\theta)}
    \]
    is an isomorphism. Moreover, it is isomorphic to $\pi^0_{(S,\theta)}$. 
\end{lem}
\begin{proof}
    Since $\msf{W}(0)$ is the parabolic Deligne-Lusztig variety attached to $(\msf{P}, w)$, the first claim follows from \cite[Théorème 10.7]{BR03} as $(\msf{S}, \theta)$ is $(\msf{G}^\sigma, \msf{M}^{w\sigma})$-superregular. The second claim follows from the composability of twisted Deligne-Lusztig inductions proved in \cite[(10.2)]{BR03}. 
\end{proof}

\begin{defi}
    Let $\chi = \theta\vert_{Z_G(F)}$. Let
    \[
        \pi_{(S, \theta)} = \cInd_{Z_G(F)\cl{G}(O_F)}^{G(F)} (\chi \otimes \pi^0_{(S,\theta)}), \quad 
        \rho_{(S, \theta)} = \cInd_{Z_G(F)G_{b,\mfr{f}}}^{G_b(F)} (\chi \otimes \rho^0_{(S,\theta)})
    \]
    be depth-zero regular supercuspidal irreducible representations. 
\end{defi}

\begin{prop}
    In the parametrization of \Cref{thm:Depth0_LLC}, we have 
    \[
        \pi_{(S, \theta)} \cong \pi_0 ,\quad \rho_{(S, \theta)} \cong \pi_\mu. 
    \]
    In particular, $\pi_{(S, \theta)}$ is generic. 
\end{prop}
\begin{proof}
    It is easy to check that when $\nu = 0$, we have
    \[
        x_0 = o ,\quad
        y_0 = 1 ,\quad 
        J_0 = \{o\} \subset \mfr{a} = C_0
    \]
    and when $\nu = \mu$, we have 
    \[
        x_\mu = \mfr{f}^{b\sigma = 1} ,\quad
        y_\mu = b ,\quad 
        J_\mu = \mfr{f} \subset \mfr{a} = C_\mu. 
    \]
    The first claim can be checked directly from the construction of $\pi_\nu$. The second claim is a general fact, which holds when $x_\nu$ is hyperspecial (see \cite[Lemma 6.2.1]{DR09}). 
\end{proof}

\subsection{The contribution from our special affinoid} \label{ssec:DeBacker-Reeder}

In this section, we study the contribution of $\pi_{(S, \theta)}$ and $\rho_{(S, \theta)}$ in the middle cohomology of local Shimura varieties. 

\begin{defi}
    Let 
    \[
        H_c^r(G,b,\mu) = \colim_{K \subset G(F)} H_c^r(\cl{M}_{G,b,\mu,K,\bb{C}_p}, \ov{\bb{Q}}_\ell)
    \]
    be the middle cohomology of the tower of local Shimura varieties attached to $(G,b,\mu)$, equipped with a natural action of $G(F) \times G_b(F) \times W_E$. 
\end{defi}

\begin{rmk}
    The diagonal action $Z_G(F) \subset G(F) \times G_b(F)$ is trivial, since it is already trivial at the level of the geometry of local Shimura varieties. 
\end{rmk}

\begin{defi}
    Let $\theta \circ \lambda$ denote the tame character 
    \[
        I_F \to I_F/I_e \cong \bb{Z}/e\bb{Z} \xrightarrow{\tau \mapsto (\zeta_e^\tau)^{e\lambda}} \msf{Z}_{\msf{M}^{w\sigma}} \xrightarrow{\theta^0} \ov{\bb{Q}}_\ell^\times.
    \]
\end{defi}

\begin{lem} \label{lem:depthzero_LCF}
    As a character of $I_F$, we have $(\theta \circ \lambda)^{-1} = \mu \circ \varphi_{(S, \theta)}\vert_{I_F}$. Here, $\varphi_{(S, \theta)}\vert_{I_F}$ is regarded as a homomorphism $I_F \to \widehat{T}(\Qla)$ corresponding to $\theta^0$ via local class field theory. 
\end{lem}
\begin{proof}
    Recall the description of depth-zero local class field theory in \cite[{}4.3]{DR09}. Let $N \geq 1$ be an integer such that $(w\sigma)^N = \sigma^N$ and let $k_N$ be the finite extension of $k$ of degree $N$. We have a map
    \[
        I_F \to k_N^\times, \quad \tau \mapsto \zeta_{q^N-1}^\tau = \tau(\pi^{1/(q^N-1)}) / \pi^{1/(q^N-1)}. 
    \]
    Here, $\zeta_{q^N-1}^\tau$ is independent of the choice of $\pi^{1/(q^N-1)}$ in the algebraic closure of $F$. By fixing an identification $\cl{S} \cong \cl{\breve{T}}^{v\sigma}$, $\mu \circ \varphi_{(S, \theta)}\vert_{I_F}$ is given by the composition
    \[
        I_F \to k_N^\times \xrightarrow{\mu} \msf{T}(k_N) \xrightarrow{\prod_{0 \leq i < N} (v\sigma)^i} \msf{S}(k) \xrightarrow{\theta^0} \Qlax. 
    \]
    Since $(v\sigma)^i \mu = (w\sigma)^i \mu \in X_*(\msf{Z}_{\msf{M}})$, the composition $k_N^\times \to \msf{T}(k_N) \to \msf{S}(k)$ factors as 
    \[
        \sum_{0 \leq i < N} (qw\sigma)^i \mu \colon k_N^\times \to \msf{Z}_{\msf{M}^{w\sigma}}(k) \subset \msf{S}(k). 
    \]
    By construction, $\sum_{0 \leq i < N} (qw\sigma)^i \mu = - (q^N - 1) \lambda$ and $(\zeta_{q^N-1}^\tau)^{(q^N - 1) \lambda} = (\zeta_{e}^\tau)^{e\lambda}$ for every $\tau \in I_F$. Thus, we get the claim. 
\end{proof}

\begin{prop} \label{prop:computcoh}
    There is a finite-dimensional smooth representation $\sigma_{(S, \theta)}$ of $W_E$ such that we have
    \[
        \pi_{(S, \theta)} \boxtimes \rho^{\vee}_{(S,\theta)} \boxtimes \sigma_{(S, \theta)} \subset H_c^r(G,b,\mu) \otimes_{Z_G(F) \times Z_G(F)} \chi \boxtimes \chi^{-1}
    \]
    and $\sigma_{(S, \theta)}\vert_{I_F}$ decomposes into Frobenius twists of $\mu \circ \varphi_{(S, \theta)}\vert_{I_F}$. In particular, $\sigma_{(S, \theta)}$ is tame. 
\end{prop}
\begin{proof}
    By \Cref{lem:middegrho}, 
    \begin{equation} \label{eq:subrep_finite_level}
        \pi^0_{(S,\theta)} \boxtimes \rho^{0, \vee}_{(S,\theta)} \subset H_c^r(\msf{W}(0), \bb{\ov{Q}}_{\ell})
    \end{equation}
    is an irreducible subrepresentation that maps injectively to $H^r(\msf{W}(0), \Qla)$. By \Cref{prop:inertact}, the action of $\tau \in I_F$ is equal to the action of $\zeta_e^{\tau} \in \msf{Z}_{\msf{M}^{w\sigma}}$. Thus, 
    \[
        \pi^0_{(S,\theta)} \boxtimes \rho^{0, \vee}_{(S,\theta)} \boxtimes (\theta\circ \lambda)^{-1} \subset H_c^r(\msf{W}(0),\ov{\bb{Q}}_\ell)
    \]
    as a representation of $\msf{G}^\sigma \times \msf{M}^{w\sigma} \times I_F$. Applying \Cref{prop:induction_with_central_character} to \eqref{eq:subrep_finite_level} with a central character $\chi^{-1}$, we get
    \[
        (\chi \otimes \pi^0_{(S, \theta)}) \boxtimes \rho^{\vee}_{(S,\theta)}  \boxtimes \mu \circ \varphi_{(S, \theta)}\vert_{I_F} \subset H_c^r(\cl{M}_{G, b,\mu, \cl{G}(1), \bb{C}_p}) \otimes_{Z_G(F)} \chi^{-1}
    \]
    by \Cref{lem:depthzero_LCF}. Moreover, by applying \Cref{prop:induction_with_central_character} to the span of $f_1$-fold Frobenius twists of \eqref{eq:subrep_finite_level} under \Cref{cor:W0_WeilDescent}, it follows that the span under the smooth $W_E$-action
    \[
        W_E \cdot \pi^0_{(S, \theta)} \boxtimes \rho^{\vee}_{(S,\theta)} \subset H_c^r(\cl{M}_{G, b,\mu, \cl{G}(1), \bb{C}_p}) \otimes_{Z_G(F)} \chi^{-1}
    \]
    is a finite direct sum of $\pi^0_{(S, \theta)} \boxtimes \rho^{\vee}_{(S,\theta)}$ as a $\msf{G}^\sigma \times G_b(F)$-representation. By setting 
    \[
        \sigma_{(S, \theta)} = \Hom_{\msf{G}^\sigma \times G_b(F)}(\pi^0_{(S, \theta)} \boxtimes \rho^{\vee}_{(S,\theta)}, W_E \cdot \pi^0_{(S, \theta)} \boxtimes \rho^{\vee}_{(S,\theta)}), 
    \]
    we have 
    \[ 
        (\chi \otimes \pi^0_{(S, \theta)}) \boxtimes \rho^{\vee}_{(S,\theta)}  \boxtimes \sigma_{(S, \theta)} \subset H_c^r(\cl{M}_{G, b,\mu, \cl{G}(1), \bb{C}_p}) \otimes_{Z_G(F)} \chi^{-1}.
    \]
    By taking the compact induction from $Z_G(F) \cl{G}(O_F)$ to $G(F)$, we get the claim. 
\end{proof}

\subsubsection{Computation of the Galois representation}

We would like to compute the Galois representation $\sigma_{(S, \theta)}$ more explicitly. For this, we assume that $G$ is split, to simplify the Weil descent datum as in \Cref{prop:split_WeilDescent}. Precisely, we make the following assumption. 

\begin{ass} \label{ass:split_assumption}
    We assume that $G$ is split and $b$ takes the form
    \[
        b = \mu(-\pi) w \in G(F),
    \]
    where $\mu \colon \bb{G}_m \to G$ is defined over $F$ and $w \in N_{\cl{G}}(\cl{T})(O_F)$. We define $\msf{M}^{w\sigma}$ and $\msf{W}(0)$ over $k$. 
\end{ass}

\begin{lem} \label{lem:sigma_and_GM}
    Let $\sigma_{\msf{W}(0)}$ be the Weil descent datum of $\msf{W}(0)$. Then, we have 
    \[
        \sigma_{\msf{W}(0)}(g \cdot x \cdot m) = g \cdot \sigma_{\msf{W}(0)}(x) \cdot w^{-1} m w
    \]
    for every $g \in \msf{G}^\sigma$, $m \in \msf{M}^{w\sigma}$ and $x \in \msf{W}(0)(\ov{k})$. In particular, 
    \[  
        \Fr \circ (g \times m) = (g \times w m w^{-1}) \circ \Fr
    \]
    as actions on $H_c^r(\msf{W}(0), \Qla)$ for every $g \in \msf{G}^\sigma$ and $m \in \msf{M}^{w\sigma}$. Here, $\Fr$ denotes the geometric Frobenius relative to $k$. 
\end{lem}
\begin{proof}
    By \Cref{cor:W0desc}, $\sigma_{\msf{W}(0)}$ equals the restriction of the Weil descent datum of $\msf{G}^\sigma$. Then, the first claim follows from $\sigma(g) = g$ and $\sigma(m) = w^{-1} m w$. Since the actions of the \textit{arithmetic} Frobenius and $g$ (resp.\ $m$) on $H_c^r(\msf{W}(0), \Qla)$ are the pushforward along $\sigma_{\msf{W}(0)}$ and the left (resp.\ right) multiplication by $g$ (resp.\ $m^{-1}$), the second claim follows from the first one. 
\end{proof}

\begin{lem} \label{lem:defd}
    Under \Cref{ass:split_assumption}, for each $i \geq 0$, let
    \[
        \rho_{(S,\theta)}^{0,w^{i}} = \rho_{(S,\theta)}^{0} \circ \Ad(w^{-i}) ,\quad b^i = \mu_i(-\pi)w^i
    \]
    with $\mu_i \in X_*(T)$. For every $i \geq 0$, we have 
    \[
        \Fr^i (\pi_{(S,\theta)}^0 \boxtimes \rho_{(S,\theta)}^{0,\vee}) = \pi_{(S,\theta)}^0 \boxtimes (\rho_{(S,\theta)}^{0,w^i})^\vee
    \]
    inside $H_c^r(\msf{W}(0), \Qla)$, and 
    $
        \rho_{(S,\theta)}^{0,w^i} \cong \rho_{(S,\theta)}^{0}
    $
    if and only if $w^i \in T(F)$. In this case, $\mu_i \in X_*(Z_G)$. 
\end{lem}
\begin{proof}
    The first claim follows from \Cref{lem:sigma_and_GM}. For the second claim, if $w^i \in T(F)$, $w^i \in \cl{T}(O_F) \subset G_{b,\mfr{f}}$, so $\rho_{(S,\theta)}^{0,w^i} \cong \rho_{(S,\theta)}^{0}$. We will prove the converse. 
    
    Suppose $\rho_{(S,\theta)}^{0,w^i} \cong \rho_{(S,\theta)}^{0}$. Since $(\msf{S}, \theta)$ is in general position and 
    \[
        \rho_{(S,\theta)}^{0,w^i} \cong \varepsilon_{(\msf{M}^{w\sigma}, \msf{S})} R_{\Ad(w^{i})(\msf{S})}^{\msf{M}^{w\sigma}}(\theta \circ \Ad(w^{-i})), 
    \]
    it follows from \cite[Theorem 6.8]{DL76} that $\Ad(w^{i})(\msf{S},\theta) = \Ad(v)(\msf{S},\theta)$ for some $v\in \msf{M}^{w\sigma}$. 
    
    Let $\wtd{w}= \Ad(\ov{p}_w)(w)$. Since $\sigma(w) = w$, $\wtd{w}\in \msf{G}^\sigma$ and $\Ad(w)\vert_{\msf{M}^{w\sigma}}$ extends to $\Ad(\wtd{w})\vert_{\msf{G}^\sigma}$ via the embedding $\Ad(\ov{p}_w) \colon \msf{M}^{w\sigma} \subset \msf{G}^\sigma$. Then, $\wtd{w}^{-i} v \in \msf{G}^\sigma$ stabilizes $(\msf{S},\theta)$. Since $(\msf{S},\theta)$ is in general position in $\msf{G}^\sigma$, $\wtd{w}^{-i} v \in \msf{S}$. In particular, $\wtd{w}^i \in \msf{M}^{w\sigma}$, so $w^i \in \msf{M}$. Since $w^i \in \msf{M}$, $w^i\mu_i = \mu_i$ and $b^{ij} = (j\mu_i)(-\pi)w^{ij}$ for $j\geq 0$. Since $b^{ij}$ fixes $\mfr{a}$, $j\mu_i$ is minuscule for every $j$, so $\mu_i \in X_*(\msf{Z}_G)$. Since $b^i$ fixes $\mfr{a}$, $w^i = 1$ in $W_o$, so $w^i \in T(F)$. 
\end{proof}

For each $d \geq 1$, let $F_d$ be the unramified extension of degree $d$ over $F$. By abuse of notation, let $\Fr$ also denote a geometric Frobenius element of $W_{F[\pi^{1/e}]}$. 

\begin{thm} \label{prop:computsplit}
    Let $d\geq 1$ be the minimum positive integer such that $w^d \in T(F)$. Under \Cref{ass:split_assumption}, $w^{-d} \circ \Fr^d$ acts by a scalar on $\pi_{(S,\theta)}^0 \boxtimes \rho_{(S,\theta)}^{0,\vee} \subset H_c^r(\msf{W}(0), \Qla)$. Let $\alpha_d$ be the eigenvalue of $w^{-d} \circ \Fr^d$ and consider the tame character
    \[
        \xi_{(S, \theta)} \colon W_{F_d} \to \Qlax ,\quad \xi_{(S, \theta)}\vert_{I_F} = \mu \circ \varphi_{(S, \theta)}\vert_{I_F} ,\quad \xi_{(S, \theta)}(\Fr^d) = \chi(\mu_d(-\pi))\alpha_d. 
    \]
    Then, we have an injection
    \[
        \pi_{(S, \theta)} \boxtimes \rho^{\vee}_{(S,\theta)} \boxtimes \Ind_{W_{F_d}}^{W_F} \xi_{(S, \theta)} \subset H_c^r(G,b,\mu) \otimes_{Z_G(F) \times Z_G(F)} \chi \boxtimes \chi^{-1}. 
    \]
\end{thm}

Here, $w^{-d} \circ \Fr^d$ denotes the composition of the actions of $w^{-d} \in \msf{M}^{w\sigma}$ and $\Fr^d$. Note that these two actions commute by \Cref{lem:sigma_and_GM}.

\begin{proof}
    By \Cref{lem:sigma_and_GM}, $w^{- d} \circ \Fr^d$ is equivariant under $\msf{G}^\sigma \times \msf{M}^{w\sigma}$, so it acts on $\pi_{(S,\theta)}^0 \boxtimes \rho_{(S,\theta)}^{0,\vee}$ by a scalar by Schur's lemma. We will prove the second claim. By \Cref{lem:defd}, 
    \[
        \bigoplus_{0 \leq i < d} \pi_{(S,\theta)}^0 \boxtimes (\rho_{(S,\theta)}^{0,w^i})^\vee \subset H_c^r(\msf{W}(0),\ov{\bb{Q}}_\ell)
    \]
    is stable under $\sigma_*$ and maps injectively to $H^r(\msf{W}(0),\ov{\bb{Q}}_\ell)$. Applying \Cref{prop:induction_with_central_character}, we get 
    \begin{equation} \label{eq:big_subrep_first}
        \bigoplus_{0\leq i < d} \pi_{(S,\theta)}^0 \boxtimes \cInd_{Z_G(F)G_{b,\mfr{f}}}^{G_b(F)} (\chi \otimes \rho_{(S,\theta)}^{0,w^i})^\vee \subset H_c^r(\cl{M}_{G,b,\mu,\cl{G}(1), \bb{C}_p}, \ov{\bb{Q}}_\ell) \otimes_{Z_G(F)} \chi^{-1}. 
    \end{equation}
    As recalled in \eqref{eq:Weildesc}, the Weil descent datum on $\cl{M}_{G,b,\mu, \cl{G}(1), \breve{F}_e }$ restricts to $\cl{W}(0)_{\breve{F}_e}$ as 
    \[
        \cl{W}(0)_{\breve{F}_e} \xrightarrow{\sigma} \cl{W}(0)_{\breve{F}_e} \xrightarrow{b} \cl{M}_{G,b,\mu,\cl{G}(1), \breve{F}_e }
    \]
    using the Weil descent datum on $\cl{W}(0)$. %, which can be defined over $O_{F_e}$ as in \Cref{prop:split_WeilDescent}. 
    Thus, we have the commutative diagram
    \begin{center}
        \begin{tikzcd}
            H_c^r(\msf{W}(0),\ov{\bb{Q}}_\ell) \ar[r, "\iota"] \ar[d, "\Fr^{-1}"] & H_c^r(\cl{M}_{G,b,\mu,\cl{G}(1), \bb{C}_p }, \ov{\bb{Q}}_\ell) \ar[d, "\Fr^{-1}"] \\
            H_c^r(\msf{W}(0),\ov{\bb{Q}}_\ell) \ar[r, "b \circ \iota"] & H_c^r(\cl{M}_{G,b,\mu,\cl{G}(1), \bb{C}_p }, \ov{\bb{Q}}_\ell) 
        \end{tikzcd}
    \end{center}
    Thus, for each $i \geq 0$, $\Fr^i$ sends $\pi_{(S,\theta)}^0 \boxtimes \rho_{(S, \theta)}^\vee$ to $\pi_{(S,\theta)}^0 \boxtimes \cInd_{Z_G(F)G_{b,\mfr{f}}}^{G_b(F)} (\chi \otimes \rho_{(S,\theta)}^{0,w^i})^\vee $, and as 
    \[
        \Fr^{-d} \circ \iota = \mu_d(-\pi)\circ \iota \circ (w^{d}\circ \Fr^{-d}), 
    \]
    $\Fr^{-d}$ acts on $\pi_{(S,\theta)}^0 \boxtimes \rho_{(S, \theta)}^\vee$ by the scalar $\chi(\mu_d(-\pi))^{-1} \alpha^{-1}_d$. Note that the actions of $w^d$ and $\Fr$ commute by \Cref{lem:sigma_and_GM}. By \Cref{lem:depthzero_LCF}, \eqref{eq:big_subrep_first} is rewritten as
    \begin{equation} \label{eq:refined_inclusion}
        (\chi \otimes \pi_{(S,\theta)}^0) \boxtimes \rho_{(S,\theta)}^\vee \boxtimes \Ind_{W_{F_d}}^{W_F} \xi_{(S, \theta)} \subset H_c^r(\cl{M}_{G,b,\mu,\cl{G}(1), \bb{C}_p}, \ov{\bb{Q}}_\ell) \otimes_{Z_G(F) \times Z_G(F)} \chi \boxtimes \chi^{-1}. 
    \end{equation}
    Now, since
    \[
        H_c^r(\cl{M}_{G,b,\mu,\cl{G}(1), \bb{C}_p}, \ov{\bb{Q}}_\ell) = H_c^r(G, b, \mu)^{\cl{G}(1)}, 
    \]
    $H_c^r(\cl{M}_{G,b,\mu,\cl{G}(1), \bb{C}_p}, \ov{\bb{Q}}_\ell)$ is a direct summand of $H_c^r(G, b, \mu)$ as a $Z_G(F) \times Z_G(F)$-module. In particular, \eqref{eq:refined_inclusion} provides an inclusion
    \[
        (\chi \otimes \pi_{(S,\theta)}^0) \boxtimes \rho_{(S,\theta)}^\vee \boxtimes \Ind_{W_{F_d}}^{W_F} \xi_{(S, \theta)} \subset H_c^r(G, b, \mu) \otimes_{Z_G(F) \times Z_G(F)} \chi \boxtimes \chi^{-1}. 
    \]
    Since $\pi_{(S, \theta)} = \cInd_{Z_G(F) \cl{G}(O_F)}^{G(F)} (\chi \otimes \pi_{(S,\theta)}^0)$ is irreducible, we get the claim by taking the compact induction from $Z_G(F)\cl{G}(O_F)$ to $G(F)$. 
\end{proof}

\begin{rmk}
    As the image of $w^d$ in $\msf{M}^{w\sigma}$ is torsion, $\alpha_{d}$ is pure of weight $rd$, i.e. $\lvert \iota(\alpha_{d}) \rvert = q^{rd/2}$ for every isomorphism $\iota \colon \ov{\bb{Q}}_\ell \cong \bb{C}$ by the general theory of weights and \Cref{lem:middegrho}. However, we leave the determination of the exact value of $\alpha_d$, due to the lack of references on Frobenius eigenvalues on the cohomology of parabolic Deligne-Lusztig varieties. 
\end{rmk}

\begin{exa}
    The particular case studied in \cite{Yos10} (cf.\ \cite[Section 5]{Mie16}) is 
    \[
        G = \GL_n, \quad 
        \mu(t) = \begin{pmatrix}
            t^{-1} & & & \\
            & 1 & & \\
            & & \ddots & \\
            & & & 1 
        \end{pmatrix}
        ,\quad 
        w = \begin{pmatrix}
            & & & 1 \\
            1 & & & \\
            & \ddots & & \\
            & & 1 & 
        \end{pmatrix}. 
    \]
    In this case, $d = n$ and $w^d = 1$. The computation of $\alpha_d$ is done in \cite{DM85}, and we have
    \[
        \alpha_d = (-1)^{n-1} q^{n(n-1)/2} ,\quad 
        \xi_{(S, \theta)}(\Fr^n) = \chi(-\pi)^{-1}  (-1)^{n-1} q^{n(n-1)/2}. 
    \]
    Here, $\chi$ is identified with a character of $F^\times$. 
\end{exa}

\subsection{Expectation for the isotypic part of the cohomology} \label{ssec:isotcoho}

Let
\[
    R\Gamma_c(G,b,\mu) = \colim_{K \subset G(F)} R\Gamma_c(\cl{M}_{G,b,\mu,K,\bb{C}_p}, \ov{\bb{Q}}_\ell[r](\tfrac{r}{2}))
\]
be the cohomology of the tower of local Shimura varieties attached to $(G,b,\mu)$. Let 
\[
    R\Gamma(G,b,\mu)[\pi] = R\Hom_{G(F)}(R\Gamma_c(G,b,\mu), \pi)^{\sm}
\]
be the isotypic part of $\pi$ as a complex of $G_b(F) \times W_E$-representations. In this section, we collect general expectations for $R\Gamma(G,b,\mu)[\pi_{(S, \theta)}]$ that follows from the Fargues-Scholze machinery. 

For an irreducible smooth representation $\pi$, let $\varphi_\pi^\FS$ be the semisimple $L$-parameter constructed by \cite{FS24}. Let $r_{-\mu}$ be the highest weight representation of ${}^LG$ with extreme weight $-\mu$ (see \cite[Lemma 2.1.2]{Kot84}). As a special case of \cite[Theorem 3.5, Corollary 3.8]{HJ25}, we get the following. 

\begin{prop} \label{prop:Hansen_Johansson}
    If $\varphi_{\pi_{(S, \theta)}}^\FS = \varphi_{(S,\theta)}$, we have
    \begin{equation}
    R\Gamma(G,b,\mu)[\pi_{(S, \theta)}] = \bigoplus_{\delta \in \Irr(S_{\varphi}, \chi_b)} \rho_\delta[d_\delta] \boxtimes \Hom_{S_{\varphi}}(\delta^{\vee}, r_{-\mu} \circ \varphi_{(S, \theta)}). \label{eq:Kotdec} 
    \end{equation}
    Here, $\rho_\delta$ is an irreducible supercuspidal representation of $G_b(F)$ and $d_\delta \in \bb{Z}$. 
\end{prop}
\begin{proof}
    The claim follows from \cite[Theorem 3.5, Corollary 3.8]{HJ25} by applying to the local Shimura datum $(G_b,b^{-1}, \mu^{-1})$ and $\pi_{(S, \theta)} \in \Irr(G(F))$, since $\varphi_{(S,\theta)}$ is supercuspidal and $S_{\varphi}$ is abelian. 
\end{proof}

The assumption $\varphi_{\pi_{(S, \theta)}}^\FS = \varphi_{(S, \theta)}$ is an instance of a comparison conjecture between Fargues-Scholze parameters and Kaletha's local Langlands correspondence for regular supercuspidal representations. An immediate corollary of \Cref{prop:computcoh} is a Jacquet-Langlands type result for Fargues-Scholze parameters. 

\begin{cor} \label{cor:JLforFS}
    For every $(S, \theta)$ under \Cref{ass:triple}, we have $\varphi_{\pi_{(S, \theta)}}^\FS = \varphi_{\rho_{(S, \theta)}}^\FS$. 
\end{cor}
\begin{proof}
    By \cite[\! IX.3]{FS24}, \Cref{prop:computcoh} implies 
    \[
        \rho_{(S, \theta)} \subset H^0(i^{b*}T_{r_{-\mu}}i^1_{*}\pi_{(S, \theta)}).
    \]
    Then, the claim follows from \cite[\! IX.4.1, IX.5.2]{FS24}. 
\end{proof}

\begin{rmk}
    When $p$ and $q$ are sufficiently large so that endoscopic character identities are available (see \cite[Section 3.4]{Kal11}), \cite[Theorem 1.0.2]{HKW22} implies that $\rho_{(S, \theta)}$ contributes to the Euler characteristic of $R\Gamma(G,b,\mu)[\pi_{(S, \theta)}]$. As a result, $\varphi_{\pi_{(S, \theta)}}^\FS = \varphi_{\rho_{(S, \theta)}}^\FS$ follows from the result loc. cit. in this case.  
\end{rmk}

The general expectation (see e.g. \cite{HJ25}) can be formatted as follows in our specific case.  

\begin{conj}
    Let $(S, \theta)$ be an unramified elliptic regular pair and let $\pi_{(S, \theta)} = \pi_0$ in the parametrization of \Cref{thm:Depth0_LLC}. Then, the following hold. 
    \begin{enumerate}
        \item \textup{(Comparison conjecture)} We have $\varphi^\FS_{\pi_{(S, \theta)}} = \varphi_{(S, \theta)}$. 
        \item \textup{(Kottwitz conjecture)} The parametrization $\delta \mapsto \rho_\delta$ in \Cref{prop:Hansen_Johansson} equals the one $\delta \mapsto \pi_\delta$ in \Cref{thm:Depth0_LLC}. 
        \item \textup{(Vanishing conjecture)} For every $\delta \in \Irr(S_\varphi, \chi_b)$, we have $d_\delta = 0$ in \eqref{eq:Kotdec}.  
    \end{enumerate}
\end{conj}

What we can say about this conjecture from \Cref{prop:computcoh} is as follows. 

\begin{cor}
    If $\varphi_{\pi_{(S, \theta)}}^\FS = \varphi_{(S, \theta)}$, $\rho_\delta = \rho_{(S, \theta)}$ and $d_\delta=0$ for some $\delta$ in \eqref{eq:Kotdec}. 
\end{cor}
\begin{proof}
    It follows from \Cref{prop:computcoh} due to the shift in the definition of $R\Gamma_c(G,b,\mu)$. 
\end{proof}

\begin{rmk}
    It would be interesting to find the contribution of other $\pi_\delta$ in $R\Gamma_c(G,b,\mu)$ by constructing special affinoids in $\cl{M}_{G,b,\mu,\cl{G}(1), \bb{C}_p}$ other than $\cl{W}(0)_{\bb{C}_p}$. Then, the vanishing conjecture $d_\delta = 0$ would follow for all $\delta$ by a counting argument, under the working hypothesis $\varphi_{\pi_{(S, \theta)}}^\FS = \varphi_{(S, \theta)}$. 
\end{rmk}

\section{Moduli of depth-zero level structures of $p$-divisible groups} \label{sec:level_structure_pdivgrp}

In this section, we apply the construction in \Cref{prop:replevel} to classical settings such as $p$-divisible groups and Shimura varieties. 

\subsection{Equivariant compactifications of $\GL_n$} \label{ssec:cptGLn}

In this section, we introduce several equivariant compactifications of $\GL_n$ following \cite{Kau00}. 

Let $\msf{T}\subset \msf{B}\subset \GL_n$ be the standard Borel pair of $\GL_n$. The Weyl group of $\msf{T}$ in $\msf{G}$ is isomorphic to the symmetric group $S_n$. Let $e_1,\ldots,e_n$ be the standard basis of $X_*(\msf{T})$. 

\begin{defi}
    For $0\leq \ell \leq n$, let $\sigma_\ell$ be the cone of $X_*(\msf{T})\otimes_{\bb{Z}} \bb{R}$ generated by $-\sum_{1\leq j \leq i} e_j$ ($1\leq i \leq \ell$) and $\sum_{i \leq j \leq n} e_j$ ($\ell +1 \leq i \leq n$). That is,
    \[
        \sigma_\ell=\biggl\{\sum_{1\leq i \leq n} a_ie_i \;\bigg\vert\; a_1\leq a_2 \leq \cdots \leq a_\ell \leq 0 \leq a_{\ell+1} \leq a_{\ell+2} \leq \cdots \leq a_n\biggr\}. 
    \]
\end{defi}

\begin{defi} \label{defi:toroidal_KGLn}
    Let $\Sigma^+$ be the fan generated by $\sigma_\ell$ ($0 \leq \ell \leq n$) spanning $X_*(\msf{T})_\bb{R}^+$ and let $\Sigma$ be the $S_n$-stable complete fan of $X_*(\msf{T}) \otimes_{\bb{Z}} \bb{R}$ spanned by translations of $\Sigma^+$. The toroidal compactification of $\GL_n$ associated to the $S_n$-stable complete fan $\Sigma$ is denoted by $\KGL_n$. 
\end{defi}

In fact, Kausz first introduced the compactification $\KGL_n$ via successive blowups of the following ``naive'' compactification of $\GL_n$. 

\begin{defi}
    Let 
    \[
        \GL_n \subset \Mat_n = \Spec(k[x_{ij}]) \subset \Proj(k[x_{00},x_{ij}]) = \bb{P}(\Mat_n \oplus \triv)
    \]
    be equivariant embeddings of $\GL_n$. Here, $(i, j)$ runs through $1 \leq i, j \leq n$ and $\GL_n \times \GL_n$ acts on $\Mat_n$ via the matrix multiplication. 
\end{defi}

\begin{prop}\textup{(\cite[Section 4]{Kau00})}
    For $1\leq \ell \leq n$, let $I_\ell\subset k[x_{ij}]$ be the homogeneous ideal generated by all $\ell\times \ell$-minors of the matrix $(x_{ij})_{1\leq i, j \leq n}$. Then, there is a successive blowup sequence
    \[
        \KGL_n=X^{(n)}\to X^{(n-1)} \to \cdots \to X^{(0)}=\bb{P}(\Mat_n \oplus \triv), 
    \]
    where the center of $X^{(\ell)} \to X^{(\ell-1)}$ is the strict transform of $V^+(I_\ell) \cup V^+(x_{00}, I_{\ell+1}) \subset \bb{P}(\Mat_n \oplus \triv)$ in $X^{(\ell-1)}$. 
\end{prop}
\begin{proof}
    In the first part of \cite[Section 4]{Kau00}, $\KGL_n$ is first constructed as the successive blowup as in the statement. We will explain the identification of this blow up as a toroidal compactification. 
    
    Let $\msf{U}$ (resp.\ $\msf{\ov{U}}$) be the unipotent radical of (resp.\ the opposite to) the standard Borel subgroup and let $\msf{\ov{T}}^+$ be the toric variety associated to the fan $\Sigma^+$. For each $\alpha \in S_n$, let $n_\alpha$ be the associated permutation matrix. Then, \cite[Proposition 4.1, Lemma 4.4]{Kau00} constructs an open covering $\{X(\alpha, \beta)\}_{\alpha, \beta\in S_n}$ of $\KGL_n$ such that we have an isomorphism
    \[
        n_\alpha \cdot \ov{\msf{U}} \times \msf{\ov{T}}^+ \times \msf{U} \cdot n_\beta^{-1} \cong X(\alpha, \beta)
    \]
    extending the inclusion of the associated big cell
    \[
        n_\alpha \cdot \ov{\msf{U}} \times \msf{T} \times \msf{U} \cdot n_\beta^{-1} = X(\alpha, \beta) \cap \GL_n
    \]
    (see the expression in \cite[p.560]{Kau00}). Let $Y$ denote the toroidal compactification of $\GL_n$ as in \Cref{defi:toroidal_KGLn}. Then, we similarly have an open embedding $X(\alpha, \beta) \hookrightarrow Y$ extending the inclusion of the big cell. Since $Y$ is separated, the open embeddings are compatible for any two choices of $(\alpha, \beta)$, so we can glue them to an open embedding $\KGL_n \hookrightarrow Y$. Since $\KGL_n$ is proper and $Y$ is irreducible, we get $\KGL_n \cong Y$. 
\end{proof}

\begin{rmk}
    The main interest of \cite{Kau00} is a moduli interpretation of $\KGL_n$. We will not use it in our application, but it would be interesting if one can study flat moduli spaces of type $\KGL_n$ more explicitly from this moduli interpretation. 
\end{rmk}

\subsection{Local shtukas associated to $p$-divisible groups} \label{ssec:pdiv_vs_shtuka}

In this section, we explain how to attach local shtukas to $p$-divisible groups. First, we recall the Dieudonn\'{e} theory over perfectoid rings. 

\begin{thm}\textup{(\cite[Theorem 17.5.2]{SW20})}
    Let $R$ be a perfectoid ring. There is a covariant fully faithful functor
    \[
        \bb{D} \colon \{ \text{$p$-divisible groups over $R$} \} \to \{ \text{BKF-modules over $R$} \}
    \]
    with the essential image consisting of BKF-modules $(M, \varphi_M)$ with $M \subset \varphi_M(\varphi^*M) \subset \xi_R^{-1} M$. 
\end{thm}

From now on, let $F = \bb{Q}_p$ and $\cl{G} = \GL_n$ over $\bb{Z}_p$. 

\begin{defi}
    Let $R$ be a $p$-adically complete ring and let $X$ be a $p$-divisible group of height $n$ over $R$. Then, we define the local $\cl{G}$-shtuka $\cl{P}_X$ over $\Spd(R)_{/\bb{Z}_p}$ so that $\cl{P}_X(S, S^+)$ is the restriction of the local $\cl{G}$-shtuka associated to $\bb{D}(X_{S^+})$ via \Cref{prop:GBKF} for every perfectoid Huber pair $(S, S^+)$ over $R$. Equivalently, for every perfectoid $R$-algebra $S$, 
    \[
        \cl{P}_X \times_{\Spd(R)} \Spd(S) \cong \Isom(W(S^\flat)^{\oplus n}, \bb{D}(X_S))
    \]
    via \Cref{prop:GBKF}. 
\end{defi}

\begin{defi}
    A full flag of $p$-divisible groups of height $n$ is a chain of isogenies
    \[
        X_1 \xrightarrow{\alpha_1} X_2 \xrightarrow{\alpha_2} \cdots \xrightarrow{\alpha_n} X_1
    \]
    of $p$-divisible groups such that each $\alpha_i$ is of degree $p$ and $\alpha_n \circ \alpha_{n-1} \circ \cdots \circ \alpha_1 = p$. 
\end{defi}

Let $\cl{I}$ be the standard Iwahori subgroup scheme of $\cl{G}$. Then, an $\cl{I}$-torsor over a $p$-torsion free $\bb{Z}_p$-algebra $A$ corresponds to a flag
\[
    M_1 \xrightarrow{\alpha_1} M_2 \xrightarrow{\alpha_2} \cdots \xrightarrow{\alpha_n} M_1
\]
of finite projective $A$-modules of rank $n$ such that $\Cok(\alpha_i)$ is an invertible $A / pA$-module and $\alpha_n \circ \alpha_{n-1} \circ \cdots \circ \alpha_1 = p$. 

\begin{defi}
    Let $R$ be a $p$-adically complete ring and let $(X_\bullet, \alpha_\bullet)$ be a full flag of $p$-divisible groups of height $n$ over $R$. Then, we define the local $\cl{I}$-shtuka $\cl{P}_{(X_\bullet, \alpha_\bullet)}$ over $\Spd(R)_{/\bb{Z}_p}$ so that $\cl{P}_{(X_\bullet, \alpha_\bullet)}(S, S^+)$ is the restriction of the $\cl{I}$-BKF module associated to 
    \[
        \bb{D}(X_1) \xrightarrow{\bb{D}(\alpha_1)} \bb{D}(X_2) \xrightarrow{\bb{D}(\alpha_2)} \cdots \xrightarrow{\bb{D}(\alpha_n)} \bb{D}(X_1). 
    \]
\end{defi}

\subsection{Interpretation of Drinfeld level structures} \label{ssec:Drinfeld}
In this section, we study the relation between Drinfeld level structures and the construction of \Cref{prop:replevel}. 

\begin{defi}
    Let $X$ be a $p$-divisible group of height $n$ over a ring $R$. A depth-zero Drinfeld level structure is a homomorphism $\alpha\colon \bb{F}_p^{\oplus n} \to X[p](R)$ such that $\Img(\alpha)$ forms a full set of sections of $X[p]$ in the sense of \cite[Section 1.8]{KM85}. 
\end{defi}

\begin{prop}\textup{(\cite[Proposition 1.9.1]{KM85})}
    Let $X$ be a $p$-divisible group of height $n$ over a ring $R$. The moduli space of depth-zero Drinfeld level structures is represented by a closed subscheme $X_0^\Dr \subset X[p]^{\oplus n}$. In particular, $X^\Dr_0$ is finite over $R$. 
\end{prop}

\begin{defi}
    Let $X$ be a $p$-divisible group of height $n$ over a ring $R$. Let $X_{0,\flat}^\Dr$ denote the Zariski closure of $X_{0}^\Dr[\tfrac{1}{p}]$ in $X_0^\Dr$. 
\end{defi}

Now, suppose that $R$ is $p$-adically complete and admits a good cover. 
Then, we can take the flat moduli space of depth-zero level structures of type $\bb{P}(\Mat_n \oplus \triv)$ on $\cl{P}_X$. We will compare it with $X_{0,\flat}^\Dr$. In the proof, we use the prismatic Dieudonn\'{e} theory developed in \cite{ALB23}. Precisely, we take the naive dual of the contravariant prismatic Dieudonn\'{e} theory developed in \cite{ALB23}, so that it is compatible with the perfectoid Dieudonn\'{e} theory. 

\begin{prop} \label{prop:Drlevel}
    Let $R$ be a $p$-adically complete ring admitting a good cover and let $X$ be a $p$-divisible group of height $n$ over $R$. Let $\mfr{Y}$ be the flat moduli space of depth-zero level structures of type $\bb{P}(\Mat_n \oplus \triv)$ on $\cl{P}_X$ and let $X_{0, \flat}^{\Dr,\wedge}$ denote the $p$-adic completion of $X_{0, \flat}^{\Dr}$. Then, there is a morphism $\mfr{Y} \to X_{0, \flat}^{\Dr,\wedge}$ inducing $\mfr{Y}^\diamond \cong (X_{0, \flat}^{\Dr,\wedge})^\diamond$. 
\end{prop}
\begin{proof}
    First, we construct a map 
    \[
        (X_{0, \flat}^{\Dr,\wedge})^\diamond \to [\bb{P}(\Mat_n \oplus \triv)^\diamond/(L_W^+\cl{G})^\diamondsuit] \times_{[\ast/(L_W^+\cl{G})^\diamondsuit]} \Sht_{\cl{G}}. 
    \]
    Let $(S,S^+)$ be a perfectoid Huber pair over $X_{0, \flat}^{\Dr,\wedge}$. Since $X_{0}^\Dr$ is affine, we get $\Spec(S^+) \to X_{0}^\Dr$. It corresponds to a depth-zero Drinfeld level structure $\alpha\colon \bb{F}_p^{\oplus n} \to X[p](S^+)$. 
    
    By \cite[Theorem 5.4]{ALB23}, the torsion Dieudonn\'{e} theory attaches to $\alpha$ a homomorphism 
    \[
        \bb{D}(\alpha) \colon (S^{\flat+})^{\oplus n} \to \bb{D}(X_{S^+})/p
    \]
    and $\alpha$ is uniquely determined by $\bb{D}(\alpha)$. By identifying $\bb{D}(\alpha)$ with a section of
    \[
        \cl{P}_X\vert_{S^{\flat+}} \times^{\GL_n} \Mat_n \cong \Hom((S^{\flat+})^{\oplus n}, \bb{D}(X_{S^+})/p),
    \]
    we get a map
    \[
        \Spd(S,S^+) \to [\Mat_n^\diamond/(L_W^+\cl{G})^\diamondsuit] \hookrightarrow [\bb{P}(\Mat_n \oplus \triv)^\diamond/(L_W^+\cl{G})^\diamondsuit]. 
    \]
    
    Let $Y= [\bb{P}(\Mat_n \oplus \triv)^\diamond/(L_W^+\cl{G})^\diamondsuit] \times_{\lbrack \ast/(L_W^+\cl{G})^\diamondsuit \rbrack, \cl{P}_X} \Spd(R)$ and $Z=\Sht_{\cl{G}^+} \times_{\Sht_{\cl{G}}, \cl{P}_X} \Spd(R)_\eta$. The above construction provides $(X_{0, \flat}^{\Dr,\wedge})^\diamond \to Y$ and it is injective since the association $\alpha \mapsto \bb{D}(\alpha)$ is injective. To see that the generic fiber maps to $Z$, suppose that $(S,S^+)$ is in characteristic $0$. Then, $\xi_{S^+,0}$ is a non-zero-divisor in $S^{\flat+}$. Using the data of torsion Dieudonn\'{e} modules as in \cite[Definition 5.1]{ALB23}, we see that $\bb{D}(\alpha)[1/\xi_{S^+,0}]$ is a Frobenius-equivariant isomorphism, so $\bb{D}(\alpha)$ induces $\Spd(S,S^+) \to \Sht_{\cl{G}^+,\eta}$ by \Cref{lem:finetdesc}. Thus, $(X_{0}^\Dr)^\diamond \to Y$ sends the generic fiber to $Z$. Since $R$ is Noetherian, the generic fiber is dense in $(X_{0, \flat}^{\Dr,\wedge})^\diamond$ by \cite[Lemma 4.4]{Lou20}, so $(X_{0, \flat}^{\Dr,\wedge})^\diamond \to Y$ factors through $\mfr{Y}^\diamond$. Since $X_0^\Dr$ is finite over $R$, $(X_{0, \flat}^{\Dr,\wedge})^\diamond \to \mfr{Y}^\diamond$ is a proper injection, so it is a closed immersion. By comparing degrees over $\Spd(R)_\eta$, it is easy to see that $(X_{0, \flat}^{\Dr,\wedge})^\diamond \to \mfr{Y}^\diamond$ is an isomorphism on the generic fiber. Since $\mfr{Y}^\diamond$ is the $v$-closure of $Z$ in $Y$, we get $(X_{0, \flat}^{\Dr,\wedge})^\diamond \cong \mfr{Y}^\diamond$. Since $\mfr{Y}$ admits a maximal good cover, it is uniquely represented by a morphism $\mfr{Y} \to X_{0,\flat}^{\Dr, \wedge}$ by \Cref{lem:uniqmapmaxlgood}. 
\end{proof}

Next, we study the relation between $\Gamma_1(p)$-level structures on full flags of $p$-divisible groups and the construction of \Cref{prop:replevel}. First, we recall $\Gamma_1(p)$-level structures. Our main reference is \cite[Section 3.3]{HR12}. 

\begin{defi}
    A $\Gamma_1(p)$-level structure on a full flag $(X_\bullet, \alpha_\bullet)$ of $p$-divisible groups of height $n$ is the set of homomorphisms $\gamma_i \colon \bb{F}_p \to \Ker(\alpha_i)$ for every $1 \leq i \leq n$ such that $\Img(\gamma_i)$ forms a full set of sections of $\Ker(\alpha_i)$. 
\end{defi}

\begin{prop}\textup{(Oort-Tate)}
    Let $(X_\bullet, \alpha_\bullet)$ be a full flag of $p$-divisible groups of height $n$ over a $p$-adically complete ring $R$. The moduli space of $\Gamma_1(p)$-level structures is represented by a closed subscheme $(X_\bullet, \alpha_\bullet)_{\Gamma_1(p)} \subset \Ker(\alpha_1)\times \cdots \times \Ker(\alpha_n)$ that is finite and flat over $R$. 
\end{prop}

Here, each $\gamma_i$ is bijective to an Oort-Tate generator of $\Ker(\alpha_i)$ via the map $\gamma_i \mapsto \gamma_i(1)$. 

\begin{defi}
    Let $(X_\bullet, \alpha_\bullet)$ be a full flag of $p$-divisible groups of height $n$ over a ring $R$. Let $(X_\bullet, \alpha_\bullet)_{\Gamma_1(p), \flat}$ denote the Zariski closure of $(X_\bullet, \alpha_\bullet)_{\Gamma_1(p)}[\tfrac{1}{p}]$ in $(X_\bullet, \alpha_\bullet)_{\Gamma_1(p)}$. 
\end{defi}

Now, the maximal reductive quotient $\msf{I}$ of $\cl{I}$ is the diagonal maximal torus of $\GL_n$. We take an equivariant compactification of $\msf{I}$ as follows: 
\[
    \msf{I} = \bb{G}_m^{n} \subset \std_1 \times \cdots \times \std_n \subset (\bb{P}^1)^n = \msf{\ov{I}}. 
\]
Here, $\std_i \in \Rep(\msf{I})$ denotes the standard representation of the $i$-th component of $\msf{I} \cong \bb{G}_m^n$ and the $i$-th component of $\msf{\ov{I}}$ is $\bb{P}(\std_i \oplus \triv)$ as an $\msf{I}$-variety.

\begin{prop} \label{prop:Gamma_1(p)-level}
    Let $R$ be a $p$-adically complete ring admitting a good cover and let $(X_\bullet, \alpha_\bullet)$ be a full flag of $p$-divisible groups of height $n$ over $R$. Let $\mfr{Y}$ be the flat moduli space of depth-zero level structures of type $\ov{\msf{I}}$ on $\cl{P}_{(X_\bullet, \alpha_\bullet)}$ and let $(X_\bullet, \alpha_\bullet)_{\Gamma_1(p), \flat}^{\wedge}$ denote the $p$-adic completion of $(X_\bullet, \alpha_\bullet)_{\Gamma_1(p), \flat}$. Then, there is a morphism $\mfr{Y} \to (X_\bullet, \alpha_\bullet)_{\Gamma_1(p), \flat}^{\wedge}$ inducing $\mfr{Y}^\diamond \cong ((X_\bullet, \alpha_\bullet)_{\Gamma_1(p), \flat}^{\wedge})^\diamond$. 
\end{prop}
\begin{proof}
    First, we construct a map 
    \[
        ((X_\bullet, \alpha_\bullet)_{\Gamma_1(p), \flat}^{\wedge})^\diamond \to [\msf{\ov{I}}^\diamond/(L_W^+\cl{I})^\diamondsuit] \times_{[\ast/(L_W^+\cl{I})^\diamondsuit]} \Sht_{\cl{I}}. 
    \]
    Let $(S,S^+)$ be a perfectoid Huber pair over $(X_\bullet, \alpha_\bullet)_{\Gamma_1(p), \flat}^{\wedge}$. Since $(X_\bullet, \alpha_\bullet)_{\Gamma_1(p), \flat}$ is affine, we get $\Spec(S^+) \to (X_\bullet, \alpha_\bullet)_{\Gamma_1(p), \flat}$. It corresponds to a $\Gamma_1(p)$-level structure $(\gamma_i)_{1 \leq i \leq n}$ over $S^+$. 
    
    By \cite[Theorem 5.4]{ALB23}, the torsion Dieudonn\'{e} theory attaches to $\gamma_i$ a homomorphism 
    \[
        \bb{D}(\gamma_i) \colon S^{\flat+} \to \bb{D}(\Ker(\alpha_i)) \cong \Cok(\bb{D}(\alpha_i))
    \]
    and $\gamma_i$ is uniquely determined by $\bb{D}(\gamma_i)$. By identifying $\bb{D}(\gamma_i)$ with a section of
    \[
        \cl{P}_{(X_\bullet, \alpha_\bullet)}\vert_{S^{\flat+}} \times^{\cl{I}_{k}} \std_i \cong \Hom(S^{\flat+}, \Cok(\bb{D}(\alpha_i))),
    \]
    we get a map
    \[
        \Spd(S,S^+) \to [(\std_1 \times \cdots \times \std_n)^\diamond/(L_W^+\cl{I})^\diamondsuit] \hookrightarrow [\msf{\ov{I}}^\diamond/(L_W^+\cl{I})^\diamondsuit]. 
    \]
    
    Let $Y= [\msf{\ov{I}}^\diamond/(L_W^+\cl{I})^\diamondsuit] \times_{\lbrack \ast/(L_W^+\cl{I})^\diamondsuit \rbrack, \cl{P}_{(X_\bullet, \alpha_\bullet)}} \Spd(R)$ and $Z=\Sht_{\cl{I}^+} \times_{\Sht_{\cl{I}}, \cl{P}_{(X_\bullet, \alpha_\bullet)}} \Spd(R)_\eta$. The above construction provides $((X_\bullet, \alpha_\bullet)_{\Gamma_1(p), \flat}^{\wedge})^\diamond \to Y$ and it is injective since the association $\gamma_i \mapsto \bb{D}(\gamma_i)$ is injective. To see that the generic fiber maps to $Z$, suppose that $(S,S^+)$ is in characteristic $0$. Then, $\xi_{S^+,0}$ is a non-zero-divisor in $S^{\flat+}$. Using the data of torsion Dieudonn\'{e} modules as in \cite[Definition 5.1]{ALB23}, we see that $\bb{D}(\gamma_i)[1/\xi_{S^+,0}]$ is a Frobenius-equivariant isomorphism, so $(\bb{D}(\gamma_i))_{1\leq i \leq n}$ induces $\Spd(S,S^+) \to \Sht_{\cl{I}^+,\eta}$ by \Cref{lem:finetdesc}. Thus, $((X_\bullet, \alpha_\bullet)_{\Gamma_1(p), \flat}^{\wedge})^\diamond \to Y$ sends the generic fiber to $Z$. Then, the claim follows by the same argument as in \Cref{prop:Drlevel}. 
    %Since $R$ is Noetherian, the generic fiber is dense in $((X_\bullet, \alpha_\bullet)_{\Gamma_1(p), \flat}^{\wedge})^\diamond$ by \cite[Lemma 4.4]{Lou20}, so $((X_\bullet, \alpha_\bullet)_{\Gamma_1(p), \flat}^{\wedge})^\diamond \to Y$ factors through $\mfr{Y}^\diamond$. Since $X_0^\Dr$ is finite over $R$, $(X_{0, \flat}^{\Dr,\wedge})^\diamond \to \mfr{Y}^\diamond$ is a proper injection, so it is a closed immersion. By comparing degrees over $\Spd(R)_\eta$, it is easy to see that $(X_{0, \flat}^{\Dr,\wedge})^\diamond \to \mfr{Y}^\diamond$ is an isomorphism on the generic fiber. Since $\mfr{Y}^\diamond$ is the $v$-closure of $Z$ in $Y$, we get $(X_{0, \flat}^{\Dr,\wedge})^\diamond \cong \mfr{Y}^\diamond$. Since $\mfr{Y}$ admits a maximal good cover, it is uniquely represented by a morphism $\mfr{Y} \to X_{0,\flat}^{\Dr, \wedge}$ by \Cref{lem:uniqmapmaxlgood}. 
\end{proof}

\subsection{Depth-zero Cartier duality} \label{ssec:dual}
In this section, we explain a duality on the underlying spaces of flat moduli spaces of depth-zero level structures of type $\KGL_n$, derived from the Chevalley involution. %The result in this section will not be used in other parts of this paper. 
First, let us recall the description of the Dieudonn\'{e} module of the Cartier dual. 

\begin{prop}\textup{(\cite[Proposition 4.73]{ALB23}, \cite[Corollary 1.3]{Mon22})} \label{prop:Dieudonne_Cartier_dual}
    Let $X$ be a $p$-divisible group over a perfectoid ring $R$ and let $X^\vee=\Hom(X,\mu_{p^\infty})$ be the Cartier dual of $X$. Then, 
    \[
        \bb{D}(X^\vee) \cong \bb{D}(X)^\vee\{1\}. 
    \]
    Here, $\{1\}$ denotes the Breuil-Kisin twist. % and the Frobenius of $\bb{D}(X)^\vee=\Hom(\bb{D}(X),W(R^\flat))$ is given by $\varphi_{\bb{D}(X)^\vee}(f)(a) = \sigma(f(\varphi_{\bb{D}(X)}^{-1}(a)))$ for $a\in \bb{D}(X)$ and $f\in \bb{D}(X)^\vee$. 
\end{prop}

Now, let $R$ be a complete adic ring with an ideal of definition $I\subset R$ containing $p$ and let $X$ be a $p$-divisible group over $R$ with Cartier dual $X^\vee$. Suppose that $R$ admits a good cover as a $p$-adic ring. Then, let $\mfr{Y}_X$ be the flat moduli space of depth-zero level structures of type $\KGL_n$ on $\cl{P}_X$. Let $X_{0, \flat}^{\Dr,\wedge}$ and $\mfr{Y}_X^\wedge$ denote the $I$-adic completion. By the functoriality of \Cref{prop:replevel} along $\KGL_n \to \bb{P}(\Mat_n \oplus \triv)$ and \Cref{prop:Drlevel}, we get a morphism
\[
    \mfr{Y}_X^\wedge \to X_{0, \flat}^{\Dr,\wedge}. 
\]
There is no direct relation between $X_0^\Dr$ and $(X^\vee)_0^\Dr$, but after passing to $\mfr{Y}_X$, we get the following relation. Let $X_{0,\flat, \red}^{\Dr}$ denote the underlying reduced scheme of $X_{0,\flat}^{\Dr, \wedge}$. 

\begin{prop} \label{prop:Cartdual}
    Suppose that $X_{0,\flat, \red}^{\Dr}$ lies in the zero sections of $X[p]$ and the same holds true for $X^\vee$ in the above situation. Then, there is an isomorphism 
    \[
        \mfr{Y}_{X, \red}^{\wedge, \perf} \cong \mfr{Y}_{X^\vee, \red}^{\wedge, \perf}. 
    \]
\end{prop}
\begin{proof}
    Let $\ov{\PGL}_n$ be the wonderful compactification of $\PGL_n$. Let $\ov{\mfr{Y}}_X$ (resp.\ $\ov{\mfr{Y}}_{X^\vee}$) be the flat moduli space of depth-zero level structures of type $\ov{\PGL}_n$ on $\cl{P}_X\times^{\GL_n} \PGL_n$ (resp.\ $\cl{P}_{X^\vee}\times^{\GL_n} \PGL_n$). Since $\KGL_n$ is toroidal, there is a natural morphism $\KGL_n \to \ov{\PGL}_n$ and the functoriality of \Cref{prop:replevel} induces a proper map
    \[
        \mfr{Y}_X \to \mfr{\ov{Y}}_X. 
    \]
    As above, let $\mfr{\ov{Y}}_X^\wedge$ (resp.\ $\ov{\mfr{Y}}_{X^\vee}^\wedge$) denote the $I$-adic completion. 

    Recall the construction of $(X_{0,\flat}^{\Dr, \wedge})^\diamond \to [\Mat_{n}^\diamond/(L_W^+\cl{G})^\diamondsuit]$ in the proof of \Cref{ssec:Drinfeld}. Then, the assumption on $X_{0,\flat, \red}^{\Dr}$ implies that $(X_{0,\flat,\red}^\Dr)^\diamond \to [\Mat_{n}^\diamond/(L_W^+\cl{G})^\diamondsuit]$ factors through $[\{0\}/(L_W^+\cl{G})^\diamondsuit]$. By \cite[Proposition 10.1]{Kau00}, the fiber of $\KGL_n \to \bb{P}(\Mat_n \oplus \triv)$ over $\{0\}$ is isomorphic to $\ov{\PGL}_n$, so $(\mfr{Y}_{X, \red}^{\wedge})^\diamond \to [\KGL_n^\diamond/(L_W^+\cl{G})^\diamondsuit]$ factors through $[\ov{\PGL}_n/(L_W^+\cl{G})^\diamondsuit]$. Since the diagram
    \begin{center}
        \begin{tikzcd}
            \lbrack \ov{\PGL}_n^\diamond/(L_W^+ \cl{G})^\diamondsuit \rbrack \ar[r] \ar[d] & \lbrack \ov{\PGL}_n^\diamond/(L_W^+ \PGL_n)^\diamondsuit \rbrack \ar[d] \\
            \lbrack \ast / (L_W^+ \cl{G})^\diamondsuit \rbrack \ar[r] & \lbrack \ast / (L_W^+ \PGL_n)^\diamondsuit \rbrack 
        \end{tikzcd}
    \end{center}
    is Cartesian, it follows that $(\mfr{Y}_{X, \red}^{\wedge, \perf})^\diamond \to (\mfr{\ov{Y}}_{X, \red}^{\wedge, \perf})^\diamond$ is injective. On the other hand, $\mfr{Y}_X^\diamond \to \ov{\mfr{Y}}_X^\diamond$ is proper and surjective on the generic fibers, so it is surjective as the generic fiber of $\ov{\mfr{Y}}_X$ is dense. In particular, its base change $(\mfr{Y}_{X, \red}^{\wedge, \perf})^\diamond \to (\mfr{\ov{Y}}_{X, \red}^{\wedge, \perf})^\diamond$ is also surjective, so we get $(\mfr{Y}_{X, \red}^{\wedge, \perf})^\diamond \cong (\mfr{\ov{Y}}_{X, \red}^{\wedge, \perf})^\diamond$. The same holds true for $X^\vee$, so by \cite[Proposition 18.3.1]{SW20}, it is enough to show 
    \[
        \mfr{\ov{Y}}_{X, \red}^{\wedge, \perf} \cong \mfr{\ov{Y}}_{X^\vee, \red}^{\wedge, \perf}. 
    \] 

    Let $\theta\colon \GL_n \to \GL_n$ be the Chevalley involution sending $g$ to $(g^{t})^{-1}$. In passing from rank $n$ vector bundles to $\GL_n$-torsors, taking the dual is interpreted as the pullback along $\theta$. Moreover, twists by line bundles are interpreted as the pushforward along $\GL_n \times \bb{G}_m \to \GL_n$. Thus, we have 
    \[
        \cl{P}_{X^\vee}\times^{\GL_n} \PGL_n \cong \theta^*(\cl{P}_{X}\times^{\GL_n} \PGL_n)
    \]
    by \Cref{prop:Dieudonne_Cartier_dual}, since the contribution of the Breuil-Kisin twist is killed after passing to $\PGL_n$. Since $\theta$ acts on $X_*(\msf{T})$ by $-1$, $\theta$ induces an involution of $\ov{\PGL}_n$. Then, by the functoriality of depth-zero level structures, we get $\mfr{\ov{Y}}_X \cong \mfr{\ov{Y}}_{X^\vee}$. We get the claim by restricting to the special fibers. 
\end{proof}

\subsection{Depth-zero integral models over universal deformation rings} \label{ssec:LTn_KGLn}

In this section, we study the case where a $p$-divisible group $X$ over $R$ is the universal deformation of a formal $p$-divisible group.

\begin{defi}
    Let $\bb{X}$ be a one-dimensional formal $p$-divisible group of height $n$ over $\bb{\ov{F}}_p$ and let $\LT_n$ be the universal deformation space of $\bb{X}$. Let $\LT_n = \Spf(R_{\bb{X}})$ and let $\bb{X}^\univ$ be the universal deformation of $\bb{X}$ over $R_{\bb{X}}$. Let $\mfr{m} \subset R_{\bb{X}}$ be the maximal ideal. 
\end{defi}

\begin{defi}
    Let $R_{\bb{X}}^1$ be the $p$-torsion free regular local ring representing $(\bb{X}^\univ)_0^\Dr$ (see \cite{Dr74}) and let $\LT_n^1= \Spf(R_{\bb{X}}^1)$ equipped with the $\mfr{m}$-adic topology. 
\end{defi}

\begin{prop} \label{prop:LTn1}
    As a $p$-adic ring, $R_{\bb{X}}$ admits a good cover. Let $\mfr{Y}$ be the flat moduli space of depth-zero level structures of type $\bb{P}(\Mat_n \oplus \triv)$ on $\cl{P}_{\bb{X}^\univ}$ and let $\mfr{Y}^\wedge$ be the $\mfr{m}$-adic completion. Then, we have $\mfr{Y}^\wedge \cong \LT_n^1$. 
\end{prop}
\begin{proof}
    First, $R_{\bb{X}}$ is isomorphic to $\breve{\bb{Z}}_p\llbracket T_1, \ldots, T_{n-1} \rrbracket$, so it admits a good cover as a $p$-adic ring by \cite[Proposition 5.5]{Tak26_rel} and \cite[Lemma 5.7]{Tak26_rel}. Here, we endow $R_{\bb{X}}^1$ with the $p$-adic topology. Then, it is enough to show $\mfr{Y} \cong \Spf(R_{\bb{X}}^1)$. 

    By \Cref{prop:Drlevel}, we get $\mfr{Y} \to \Spf(R_{\bb{X}}^1)$ inducing $\mfr{Y}^\diamond \to \Spd(R_{\bb{X}}^1)$. By \cite[Proposition 18.3.1]{SW20}, we have $\mfr{Y}_{\red}^{\perf} \cong \Spf(R^1_{\bb{X}})_\red^\perf$, so $\mfr{Y}$ has a unique closed point. It follows that $\mfr{Y}$ is affine and $\Gamma(\mfr{Y}, \cl{O}) \to \Gamma(\mfr{Y}_\eta, \cl{O}^+)$ is injective by \cite[Lemma 5.30]{Tak26_rel}. On the other hand, $R_{\bb{X}}^1$ is normal, so $R_{\bb{X}}^1 \cong \Gamma(\Spf(R_{\bb{X}}^1)_\eta, \cl{O}^+)$. Since $\mfr{Y}_\eta \cong \Spf(R_{\bb{X}}^1)_\eta$, we get $R_{\bb{X}}^1 \cong \Gamma(\mfr{Y}, \cl{O})$.  
\end{proof}

Next, we study the duality on the underlying space. Take $0 < d < n$ and let $\bb{X}_d$ be a formal $p$-divisible group of height $n$ and dimension $d$ over $\bb{\ov{F}}_p$ such that $\bb{X}_{n - d} = \bb{X}_d^\vee$ is also formal. Let $R_{d}$ be the universal deformation ring of $\bb{X}_d$ and let $\bb{X}_d^\univ$ be the universal deformation over $R_{d}$. Let $\mfr{m}_d \subset R_{d}$ be the maximal ideal. By the Cartier duality, $R_d$ is the universal deformation ring of $\bb{X}_{n-d}$ and $\bb{X}_{n-d}^\univ = (\bb{X}_d^\univ)^\vee$ is the universal deformation of $\bb{X}_{n - d}$.  

\begin{prop}
    As a $p$-adic ring, $R_{d}$ admits a good cover. Let $\mfr{Y}_d$ (resp.\ $\mfr{Y}_{n - d}$) be the flat moduli space of depth-zero level structures of type $\KGL_n$ on $\cl{P}_{\bb{X}_d^\univ}$ (resp.\ $\cl{P}_{\bb{X}_{n-d}^\univ}$) and let $\mfr{Y}_d^\wedge$ (resp.\ $\mfr{Y}_{n-d}^\wedge$) be the $\mfr{m}_d$-adic completion. Then, there is an isomorphism
    \[
        \mfr{Y}_{d, \red}^{\wedge, \perf} \cong \mfr{Y}_{n - d, \red}^{\wedge, \perf}. 
    \]
\end{prop}
\begin{proof}
    Since $R_{d}$ is isomorphic to $\breve{\bb{Z}}_p\llbracket T_1, \ldots, T_{d(n-d)} \rrbracket$, it admits a good cover as a $p$-adic ring by \cite[Proposition 5.5]{Tak26_rel} and \cite[Lemma 5.7]{Tak26_rel}. Since $\bb{X}_d$ and $\bb{X}_{n-d}$ are formal, depth-zero Drinfeld level structures on them are trivial. Thus, the assumption of \Cref{prop:Cartdual} holds, so the claim follows. 
\end{proof}

%Our motivation for this duality is the study of the underlying perfect scheme of $\cl{M}_{G, b, \mu, \cl{G}(1)}^\ints$ constructed in \Cref{defi:full_depth0_integral_model}. In the next section, we will study the underlying space thoroughly for Lubin-Tate spaces. 

\begin{defi} \label{defi:LTn_KGLn}
    Let $\LT_n^{\KGL_n} = \mfr{Y}_1^\wedge$ denote the $\mfr{m}$-adic completion of the flat moduli space of depth-zero level structures of type $\KGL_n$ on $\cl{P}_{\bb{X}^\univ}$. 
\end{defi}

\subsection{Depth-zero integral models of local Shimura varieties} \label{ssec:for_local_Shimura_varieties}

In this section, we explain the extension of \Cref{const:Xb} to the whole integral model of $\cl{M}_{G, b, \mu, \cl{G}(1)}$. 
Here, we consider an arbitrary (not necessarily hyperspecial) parahoric subgroup $\cl{G}$ and an arbitrary (not necessarily basic) $\sigma$-conjugacy class $[b] \in B(G, \mu)$. Instead, we assume the following. 

\begin{ass} \label{ass:integralmodel}
    The $v$-sheaf theoretic integral model $\cl{M}^\ints_{\cl{G}, b, \mu}$ is representable by a formal scheme that locally admits a \textit{good algebraization} together with a local shtuka, in the sense that for every quasicompact open subset $U \subset \cl{M}_{\cl{G}, b, \mu}^\ints$, we have
    \[
        \cl{M}_{\cl{G}, b, \mu}^\ints\vert_U \cong S^\wedge_{/Z}, \quad Z \subset S \otimes_{O_{\breve{F}}} \ov{k}
    \]
    for a normal $O_{\breve{F}}$-scheme $S$ and its closed $\ov{k}$-subscheme $Z$ such that the local $\cl{G}$-shtuka over $\cl{M}_{\cl{G}, b, \mu}^\ints\vert_U$ extends to the $p$-adic completion $\mfr{S} = S^\wedge$ and $\mfr{S}$ admits a good cover. 
\end{ass}

\begin{rmk}
    This condition is satisfied when $\cl{M}^\ints_{\cl{G}, b, \mu}$ admits a $p$-adic uniformization to a canonical integral model of a Shimura variety at hyperspecial level (cf.\ \cite[Corollary 4.9]{Mie20}). We think that the above algebraicity condition is only technical due to the limitation of \Cref{prop:replevel} and should essentially be irrelevant. 
\end{rmk}

%Then, as in \cite[Corollary 4.9]{Mie20}, for every quasicompact open subset $U \subset \cl{M}_{\cl{G}, b, \mu}^\ints$, we have
%\[
%    \cl{M}_{\cl{G}, b, \mu}^\ints\vert_U \cong S^\wedge_{/Z}, \quad Z \subset S \otimes_{O_{\breve{F}}} \ov{k}
%\]
%for a smooth $O_{\breve{F}}$-scheme $S$ and its closed $\ov{k}$-subscheme $Z$ such that the local $\cl{G}$-shtuka over $\cl{M}_{\cl{G}, b, \mu}^\ints\vert_U$ extends to the $p$-adic completion $\mfr{S} = S^\wedge$. 
Now, \Cref{ass:integralmodel} allows us to apply \Cref{prop:replevel} to $\mfr{S}$ and take the completion of the output $\mfr{S}(1)$ along $Z$ to get a $p$-torsion free formal scheme
\[
    \cl{M}(1)_U = \mfr{S}(1)^\wedge_{/Z}
\]
over $\cl{M}_{\cl{G}, b, \mu}^\ints\vert_U$. This construction depends on the algebraization, but as $\cl{M}(1)_U$ is $p$-torsion free, 
\[
    \cl{M}(1)_U^\diamond = \mfr{S}(1)^\diamond \times_{\mfr{S}^\diamond} (\cl{M}_{\cl{G}, b, \mu}^\ints\vert_U)^\diamond \subset \lbrack \msf{\ov{G}}^{\diamond} /(L_W^+\cl{G})^\diamondsuit \rbrack \times_{\lbrack \ast/(L_W^+\cl{G})^\diamondsuit \rbrack, \cl{P}} (\cl{M}_{\cl{G}, b, \mu}^\ints\vert_U)^\diamond
\]
is the $v$-closure of 
\[
    \cl{M}(1)_{U, \eta}^\diamond = \Sht_{\cl{G}^+}\times_{\Sht_{\cl{G}}, \cl{P}} (\cl{M}_{\cl{G}, b, \mu}^\ints\vert_U)_\eta^\diamond \subset \lbrack \msf{\ov{G}}^{\diamond} /(L_W^+\cl{G})^\diamondsuit \rbrack \times_{\lbrack \ast/(L_W^+\cl{G})^\diamondsuit \rbrack, \cl{P}} (\cl{M}_{\cl{G}, b, \mu}^\ints\vert_U)^\diamond
\]
by \cite[Lemma 4.4]{Lou20}. In particular, $\cl{M}(1)_U^\diamond$ is uniquely determined. 

Let $\mfr{S}(1)^\norm$ be the normalization of $\mfr{S}(1)$. Since $\mfr{S}(1)_\eta$ is normal, it is equal to the normalization inside the generic fiber. Let
\[
    \cl{M}(1)_U^\norm = (\mfr{S}(1)^\norm)^{\wedge}_{/Z}. 
\]
Thanks to the excellence of $\mfr{S}$, $\cl{M}(1)_U^\norm$ is the normalization of $\cl{M}(1)_U$. Precisely speaking, for each affine open subspace $\Spf(R) \subset \cl{M}(1)_U$, its inverse image is given by $\Spf(R^+) \subset \cl{M}(1)_U^\norm$ where $R^+ \subset R[\tfrac{1}{p}]$ is the integral closure of $R$. In particular, the natural map
\[
    \cl{M}(1)_{U, \eta}^\norm \to \cl{M}(1)_{U, \eta}
\]
is an isomorphism. 

\begin{lem} \label{lem:mapfromnormal}
    Let $\mfr{X}$ and $\mfr{Y}$ be formal schemes formally of finite type over $O_{\breve{F}}$ and let $f \colon \mfr{X}_\eta \to \mfr{Y}_\eta$ be a morphism. If $\mfr{X}$ is $p$-torsion free and normal, every map 
    \[
        \mfr{X}^\diamond \to \mfr{Y}^\diamond 
    \]
    extending $f^\diamond \colon \mfr{X}_\eta^\diamond \to \mfr{Y}_\eta^\diamond$ is uniquely representable by a morphism $\mfr{X} \to \mfr{Y}$. 
\end{lem}
\begin{proof}
    By following the argument given in the proof of \cite[Theorem 2.16]{AGLR22}, we may assume that $\mfr{X}$ and $\mfr{Y}$ are affine. Let $\mfr{X} = \Spf(R)$ and $\mfr{Y} = \Spf(S)$. By \cite[Theorem 7.4.1]{dJ95}, we have a continuous map
    \[
        S \to \Gamma(\mfr{Y}_\eta, \cl{O}^+) \to \Gamma(\mfr{X}_\eta, \cl{O}^+) = R. 
    \] 
    The induced map $\Spf(R) \to \Spf(S)$ is an extension of
    \[
        \mfr{X}_\eta \xrightarrow{f} \mfr{Y}_\eta \to \Spf(S). 
    \]
    It represents the given map $\mfr{X}^\diamond \to \mfr{Y}^\diamond$ since $\Spd(S)$ is separated and $\mfr{X}_\eta^\diamond \subset \mfr{X}^\diamond$ is dense. 
\end{proof}

\begin{cor} \label{cor:indep_normalization}
    The $p$-torsion free normal formal scheme $\cl{M}(1)_U^\norm$ is independent of the choice of $S$ and $Z$. 
\end{cor}
\begin{proof}
    Suppose that we have another choice $Z' \subset S'$ and let $\cl{M}(1)_{U, S'}$ and $\cl{M}(1)_{U, S'}^\norm$ be the outputs for this choice. Then, we have the following diagram. 
    \begin{center}
        \begin{tikzcd}
            (\cl{M}(1)_U^\norm)^\diamond \ar[d, "p"] & (\cl{M}(1)_{U, S'}^\norm)^\diamond \ar[d, "q"] \\
            \cl{M}(1)_U^\diamond \ar[r, "\cong"] & \cl{M}(1)_{U, S'}^\diamond. 
        \end{tikzcd}
    \end{center}
    By \Cref{lem:mapfromnormal}, we get a morphism
    \[
        \cl{M}(1)_U^\norm \to \cl{M}(1)_{U, S'}
    \]
    and it can be uniquely lifted to $f \colon \cl{M}(1)_U^\norm \to \cl{M}(1)_{U, S'}^\norm$ due to the normality. In the same way, we get $g \colon \cl{M}(1)_{U, S'}^\norm \to \cl{M}(1)_{U}^\norm$. By $p \circ g \circ f = p$, $q \circ f \circ g = q$ and the normality, we see that $g$ is the inverse to $f$. 
\end{proof}

Now, we get $\cl{M}(1)_U^\norm$ independently of the choice of a good algebraization. 

\begin{defi} \label{defi:full_depth0_integral_model}
    Under \Cref{ass:integralmodel}, we define an integral model
    \[
        \cl{M}^{\msf{\ov{G}}, \ints}_{G, b, \mu, \cl{G}(1)}
    \]
    of $\cl{M}_{G, b, \mu, \cl{G}(1)}$ as a $p$-torsion free normal formal scheme adic over $\cl{M}^\ints_{\cl{G}, b, \mu}$ equipped with a $\msf{G}^\sigma$-action such that for every quasicompact open subset $U \subset \cl{M}_{\cl{G}, b, \mu}^\ints$, we have
    \[
        \cl{M}^{\msf{\ov{G}}, \ints}_{G, b, \mu, \cl{G}(1)}\vert_U \cong \cl{M}(1)_U^\norm. 
    \]
    When $\msf{\ov{G}}$ is clear from the context, we omit it from the notation. 
\end{defi}

From now on, we assume that $\cl{G}$ is reductive and $[b]$ is basic, and fix an equivariant compactification $\msf{\ov{G}}$. We will explain the relation between $\cl{M}^{\ints}_{G, b, \mu, \cl{G}(1)}$ and $\mfr{X}_b$ (recall \Cref{const:Xb}). 

\begin{lem} \label{lem:relation_between_Xb_Mint}
    Let $\mfr{X}_{b}^{\pre, \norm}$ be the normalization of $\mfr{X}_{b}^\pre$ and let $\mfr{X}_{b}^\norm$ be the $\mfr{m}$-adic completion of $\mfr{X}_b^{\pre, \norm}$. Under \Cref{ass:integralmodel}, the completion of $\cl{M}^{\ints}_{G, b, \mu, \cl{G}(1)}$ at the fiber of $[1] \in (\cl{M}_{\cl{G}, b, \mu}^\ints)_\red$ is isomorphic to $\mfr{X}_{b}^\norm$. In particular, the dimension of the underlying space of $\cl{M}^{\ints}_{G, b, \mu, \cl{G}(1)}$ is at least $r$ when $[b]$ is basic and $\msf{\ov{G}}$ is a toroidal compactification. 
\end{lem}
\begin{proof}
    Since $\mfr{X}_{b, \eta}^{\pre}$ is normal, $\mfr{X}_{b, \eta}^\norm \to \mfr{X}_{b, \eta}$ is an isomorphism. Let $U \subset \cl{M}_{\cl{G}, b, \mu}^\ints$ be a quasicompact open neighborhood of $[1]$. Then, by the same argument for the independence of $\cl{M}(1)^\norm_U$ (see \Cref{cor:indep_normalization}), it follows that $\mfr{X}_b^\norm \cong (\mfr{S}(1)^\norm)^\wedge_{/[1]}$. Then, the last claim follows from \Cref{thm:role_of_integralmodel} since $\dim \msf{W}(0) = r$. 
\end{proof}

\begin{rmk}
    In the algebraization of $\cl{M}_{\cl{G}, b, \mu}^\ints\vert_U$, $S$ can be often taken to be $r$-dimensional. For example, it holds when the algebraization comes from a $p$-adic uniformization. In that case, the underlying space of $\cl{M}^{\ints}_{G, b, \mu, \cl{G}(1)}$ is at most $r$-dimensional. It is interesting that the fiber of $\cl{M}^{\ints}_{G, b, \mu, \cl{G}(1)} \to \cl{M}^{\ints}_{\cl{G}, b, \mu}$ at $[1] \in X_\mu(b)$ is $r$-dimensional when $\msf{\ov{G}}$ is a toroidal compactification. 
\end{rmk}

\subsection{Depth-zero integral models of Shimura varieties} \label{ssec:for_Shimura_varieties}

In this section, we apply \Cref{prop:replevel} to canonical integral models of Shimura varieties, following \cite[Section 6.6]{Tak26_rel}. 

Let $(\mbf{G}, \mbf{X})$ be a Shimura datum satisfying the axiom (SV5). Let $\msf{K}=\msf{K}_p\msf{K}^p\subset \mbf{G}(\bb{A}_f)$ be a compact open subgroup such that $\msf{K}^p\subset \mbf{G}(\bb{A}_f^p)$ is a sufficiently small (or neat) compact open subgroup and $\msf{K}_p\subset \mbf{G}(\bb{Q}_p)$ is a parahoric subgroup. Let $\mbf{E}$ be the reflex field of $(\mbf{G},\mbf{X})$. Let $v$ be a place of $\mbf{E}$ over $p$ and set $E=\mbf{E}_v$. The Shimura variety $\Sh_\msf{K}(\mbf{G},\mbf{X})$ is expected to admit an integral model $\scr{S}_{\msf{K}}$ over $O_E$. 
We refer to \cite[Section 6.6]{Tak26_rel} for the theory of canonical integral models $\{\scr{S}_{\msf{K}}\}_{\msf{K}^p}$ due to Pappas-Rapoport \cite{PR24}. 

Let $\cl{G}$ be a parahoric subgroup scheme of $G =\mbf{G}_{\bb{Q}_p}$ with $\cl{G}(\bb{Z}_p)=\msf{K}_p$. We take $\cl{G}^+$ as in \Cref{ssec:level_structure_in_char0}. Let $\msf{K}_p^+ = \cl{G}^+(\bb{Z}_p)$ and let $\msf{K}^+ = \msf{K}^+_p \msf{K}^p$. To consider integral models of $\Sh_{\msf{K}^+}(\mbf{G},\mbf{X})_E$, we make the following assumption. 

\begin{ass} \label{ass:global_assumption}
    Suppose that $\msf{K}_p$ is contained in a hyperspecial subgroup $\msf{K}'_p$ and there exists a system of canonical integral models $\{\scr{S}_{\msf{K}'}\}_{\msf{K}^p}$ of Shimura varieties in the sense of Pappas-Rapoport. Here, $\msf{K}' = \msf{K}'_p \msf{K}^p$ and $G$ is assumed to be unramified. 
\end{ass}
\begin{rmk} \label{rmk:PR_canonical_models}
    The existence of canonical integral models at hyperspecial levels is known in the abelian type (when $p > 2$) by \cite{IKY25}, where integral models constructed by \cite{Kis10} are shown to be canonical. For exceptional Shimura varieties, \cite{MY26} verifies that integral models constructed by \cite{BST25} for sufficiently large $p$ are canonical. 
\end{rmk}

Under this assumption, by \cite[Theorem 6.32]{Tak26_rel}, we can construct an integral model $\scr{S}_{\msf{K}}$ of $\Sh_{\msf{K}}(\mbf{G},\mbf{X})_E$ as an algebraic space. Precisely, $\scr{S}_{\msf{K}}$ is a normal flat separated algebraic space over $O_E$ that is proper and surjective over $\scr{S}_{\msf{K}'}$. Moreover, \cite[Theorem 6.32]{Tak26_rel} (2) implies that there is a universal local $\cl{G}$-shtuka $\cl{P}_{\msf{K}}$ on $\scr{S}_\msf{K}^\wedge$ and $\scr{S}_{\msf{K}}^\wedge$ admits a maximal good cover. 

Now, we fix an equivariant compactification $\msf{\ov{G}}$ over the residue field of $E$ so that \Cref{prop:replevel} can be applied over $O_E$ as in \Cref{rmk:defined_O_F}. 

\begin{prop} \label{prop:depth_zero_integral_models}
    Under \Cref{ass:global_assumption}, for every $\msf{K}^p$, there is a flat separated algebraic space $\scr{S}_{\msf{K}^+}^{\msf{\ov{G}}}$ of finite type over $O_E$ with $(\scr{S}_{\msf{K}^+}^{\msf{\ov{G}}})_E \cong \mrm{Sh}_{\msf{K}^+}(\mbf{G},\mbf{X})_E$ such that 
    \begin{enumerate}
        \item the $\msf{G}^\sigma$-action on $\mrm{Sh}_{\msf{K}^+}(\mbf{G},\mbf{X})_E$ over $\mrm{Sh}_{\msf{K}}(\mbf{G},\mbf{X})_E$ extends to $\scr{S}_{\msf{K}^+}^{\msf{\ov{G}}}$ over $\scr{S}_{\msf{K}}$, 
        \item $\scr{S}_{\msf{K}^+}^{\msf{\ov{G}}}$ is proper and surjective over $\scr{S}_{\msf{K}}$, and
        \item $(\scr{S}_{\msf{K}^+}^{\msf{\ov{G}}})^\wedge$ is the flat moduli space of $\cl{G}^+$-level structures of type $\msf{\ov{G}}$ on $\cl{P}_\msf{K}$. 
    \end{enumerate}
    Moreover, the system $\{\scr{S}_{\msf{K}^+}^{\msf{\ov{G}}}\}_{\msf{K}^p}$ admits finite \'{e}tale prime-to-$p$ Hecke correspondences and satisfies the extension property
    \[
        (\varprojlim_{\msf{K}^p} \Sh_{\msf{K}^+}(\mbf{G},\mbf{X}))(R[\tfrac{1}{p}]) = (\varprojlim_{\msf{K}^p} \scr{S}_{\msf{K}^+})(R)
    \]
    for every discrete valuation ring $R$ of mixed characteristic over $O_E$. 
\end{prop}
\begin{proof}
    Since $\scr{S}_\msf{K}^\wedge$ admits a maximal good cover, we may apply \Cref{prop:replevel} and \Cref{rmk:defined_O_F} to $\scr{S}_\msf{K}^\wedge$. Let $\cl{S}_{\msf{K^+}}$ be the flat moduli space of $\cl{G}^+$-level structures of type $\msf{\ov{G}}$ on $\cl{P}_\msf{K}$. We follow the proof of \cite[Theorem 6.32]{Tak26_rel} to algebraize $\cl{S}_{\msf{K^+}}$. 

    Let $U\subset \scr{S}_{\msf{K}'}$ be an arbitrary affine open subset and let $U=\Spec(A)$. Let $A[\tfrac{1}{p}] \to B$ be the finite \'{e}tale morphism representing $\Sh_{\msf{K}^+}(\mbf{G},\mbf{X})_E \to \Sh_{\msf{K}'}(\mbf{G},\mbf{X})_E$ over $U_\eta$. Let $B^+ \subset B$ be the integral closure of $A$. Let $\cl{S}_{\msf{K}^+,U} \subset \cl{S}_{\msf{K}^+}$ be the open formal subscheme over $U^\wedge \subset \scr{S}_{\msf{K}'}^\wedge$. Since $\cl{S}_{\msf{K}^+,U}$ admits a maximal good cover by construction, it follows from \cite[Lemma 5.30]{Tak26_rel} that $\Gamma(\cl{S}_{\msf{K}^+,U,\eta}, \cl{O}^{\circ \circ}) \hookrightarrow \Gamma(\cl{S}_{\msf{K}^+,U}, \cl{O}) \hookrightarrow \Gamma(\cl{S}_{\msf{K}^+,U,\eta}, \cl{O}^{+}) \cong (B^+)^{\wedge}$. Thus, there is a finite $A$-subalgebra $B_0\subset B^+$ such that $B_0[\tfrac{1}{p}]=B$ and there is a homomorphism $B_0^\wedge \to \Gamma(\cl{S}_{\msf{K}^+,U}, \cl{O})$. By \cite[Proposition 2.19]{Tak26_rel}, $\cl{S}_{\msf{K}^+,U} \to \Spf(B_0^\wedge)$ is a formal modification. By \cite[Theorem 3.2]{Art70}, there is a unique proper morphism $\scr{S}_{\msf{K}^+,U} \to \Spec(B_0)$ of algebraic spaces that is an isomorphism outside $V(p)$ and induces $\cl{S}_{\msf{K}^+,U} \to \Spf(B_0^\wedge)$ by taking the $p$-adic completion. It is easy to see from the uniqueness of $\scr{S}_{\msf{K}^+,U}$ that $\scr{S}_{\msf{K}^+,U}$ is independent of the choice of $B_0 \subset B$ and its construction is functorial in $U$. In particular, we can glue them to obtain an algebraic space $\scr{S}_{\msf{K}^+}$ such that $\scr{S}_{\msf{K}^+}[\tfrac{1}{p}] \cong \Sh_{\msf{K}^+}(\mbf{G},\mbf{X})_E$ and $\scr{S}_{\msf{K}^+}^\wedge \cong \cl{S}_{\msf{K}^+}$. 

    The finite \'{e}tale $\msf{G}^\sigma$-torsor $\cl{S}_{\msf{K}^+, \eta} \to \cl{S}_{\msf{K}, \eta}$ is the restriction of $\mrm{Sh}_{\msf{K}^+}(\mbf{G},\mbf{X})_E \to \mrm{Sh}_{\msf{K}}(\mbf{G},\mbf{X})_E$ since both come from level structures of $\cl{P}_{\msf{K}}$. By the full faithfulness in \cite[Theorem 2.26]{AY26}, $\cl{S}_{\msf{K^+}} \to \cl{S}_{\msf{K}}$ extends to $\scr{S}_{\msf{K^+}} \to \scr{S}_{\msf{K}}$, and the $\msf{G}^\sigma$-action on $\cl{S}_{\msf{K}^+}$ extends to $\scr{S}_{\msf{K}^+}$. Moreover, since $\scr{S}_{\msf{K}}$ and $\scr{S}_{\msf{K^+}}$ are both proper over $\scr{S}_{\msf{K}'}$, $\scr{S}_{\msf{K^+}} \to \scr{S}_{\msf{K}}$ is proper and surjective. Thus, $\scr{S}_{\msf{K}^+}$ satisfies (1)--(3). We will show that the system $\{\scr{S}_{\msf{K}^+}\}$ satisfies the remaining conditions. 

    First, we verify the extension property. Take $ \{y_{\msf{K}^p}\} \in (\varprojlim_{\msf{K}^p} \Sh_{\msf{K}^+}(\mbf{G},\mbf{X}))(R[\tfrac{1}{p}])$ and let $\{x_{\msf{K}^p}\}\in (\varprojlim_{\msf{K}^p} \Sh_{\msf{K}}(\mbf{G},\mbf{X}))(R[\tfrac{1}{p}])$ be the projection of  $\{y_{\msf{K}^p}\}$. By the extension property of $\{\scr{S}_{\msf{K}}\}$, $\{x_{\msf{K}^p}\}$ uniquely lifts to $(\varprojlim_{\msf{K}^p} \scr{S}_{\msf{K}})(R)$. By the valuative criterion for properness (see \cite[Tag 0A40]{stacks-project} for algebraic spaces), we have 
    \[
        \scr{S}_{\msf{K}^+}(R) \cong \scr{S}_{\msf{K}}(R) \times_{\Sh_{\msf{K}}(\mbf{G},\mbf{X})(R[1/p])} \Sh_{\msf{K}^+}(\mbf{G},\mbf{X})(R[\tfrac{1}{p}]). 
    \]
    Thus, $\{y_{\msf{K}^p}\}$ uniquely lifts to $(\varprojlim_{\msf{K}^p} \scr{S}_{\msf{K}^+})(R)$. 
    
    Next, we show that $\{\scr{S}_{\msf{K}^+}^{\msf{\ov{G}}}\}_{\msf{K}^p}$ admits finite \'{e}tale prime-to-$p$ Hecke correspondences. Let $\msf{K}^{p}_i$ ($i=0,1$) be sufficiently small (or neat) compact open subgroups with $g \msf{K}^p_0 g^{-1} \subset \msf{K}^p_1$ for $g \in \mbf{G}(\bb{A}_f^p)$. Let $\msf{K}_i=\msf{K}_p\msf{K}^p_i$ (resp.\ $\msf{K}^+_i=\msf{K}^+_p\msf{K}^p_i$) for $i=0,1$. By assumption, we have a finite \'{e}tale prime-to-$p$ Hecke correspondence $\scr{S}_{\msf{K}_0} \to \scr{S}_{\msf{K}_1}$ associated with $g$, and $\cl{P}_{\msf{K}_0}$ is isomorphic to the pullback of $\cl{P}_{\msf{K}_1}$ along the Hecke correspondence. Then, we have a closed immersion $i \colon \cl{S}_{\msf{K}^+_0}^\diamond \to \cl{S}_{\msf{K}^+_1}^\diamond \times_{(\scr{S}_{\msf{K}_1}^\wedge)^\diamond} (\scr{S}_{\msf{K}_0}^\wedge)^\diamond$ by the functoriality of \Cref{prop:replevel} and it is represented by a natural isomorphism $\cl{S}_{\msf{K}_0^+, \eta} \cong \cl{S}_{\msf{K}_1^+, \eta} \times_{\scr{S}_{\msf{K}_1, \eta}^\wedge} \scr{S}_{\msf{K}_0, \eta}^\wedge$ on the generic fiber. 
    
    By \cite[Lemma 5.32]{Tak26_rel}, $\cl{S}_{\msf{K}_1^+}\times_{\scr{S}_{\msf{K}_1}^\wedge} \scr{S}_{\msf{K}_0}^\wedge$ admits a maximal good cover since $\scr{S}_{\msf{K}_0} \to \scr{S}_{\msf{K}_1}$ is finite \'{e}tale. In particular, the generic fiber is dense in $\cl{S}_{\msf{K}^+_1}^\diamond \times_{(\scr{S}_{\msf{K}_1}^\wedge)^\diamond} (\scr{S}_{\msf{K}_0}^\wedge)^\diamond$, so $i$ is an isomorphism. Then, by \Cref{lem:uniqmapmaxlgood}, $i$ is uniquely represented by an isomorphism $\cl{S}_{\msf{K}^+_0} \cong \cl{S}_{\msf{K}^+_1} \times_{\scr{S}_{\msf{K}_1}^\wedge} \scr{S}_{\msf{K}_0}^\wedge$. By \cite[Theorem 3.2]{Art70}, it induces $\scr{S}_{\msf{K}^+_0} \cong \scr{S}_{\msf{K}^+_1} \times_{\scr{S}_{\msf{K}_1}} \scr{S}_{\msf{K}_0}$ providing a finite \'{e}tale Hecke correspondence $\scr{S}_{\msf{K}^+_0} \to \scr{S}_{\msf{K}^+_1}$. 
\end{proof}

One particular case of interest is the case where $\cl{G}$ is Iwahori. In this case, $\cl{G}^+$ is taken so that $\msf{G}^\sigma = \cl{G}(O_F) / \cl{G}^+(O_F)$ is a torus, %, and $\msf{\ov{I}}$ is a toric compactification of the torus $\msf{I}$. 
and $\cl{G}^+$-level structures are also known as \textit{$\Gamma_1(p)$-level structures}, which have been studied in several special cases (e.g. \cite{HR12}, \cite{Sha16}). 

\begin{exa} \label{exa:Gamma_1(p)_unitary}
    In \cite[Section 2.1]{Sha16}, some unramified unitary Shimura varieties are considered, where $G = \GL_n \times \bb{G}_m$. Let $\msf{K}_p = \cl{I}_{\GL_n} \times \bb{Z}_p^\times$ be an Iwahori subgroup. In \cite[Section 3.2]{Sha16}, an integral model $\scr{S}_{\msf{K}^+}$ of the Shimura variety is introduced at the level $\msf{K}^+_p = \cl{I}^+_{\GL_n} \times \bb{Z}_p^\times$. There, $\scr{S}_{\msf{K}^+}$ is constructed as the moduli space of $\Gamma_1(p)$-level structures for the universal $p$-divisible group over $\scr{S}_{\msf{K}}$. Then, $\scr{S}_{\msf{K}^+}$ is finite flat over $\scr{S}_{\msf{K}}$, so it is $p$-torsion free. By \Cref{prop:Gamma_1(p)-level}, $(\scr{S}_{\msf{K}^+}^\wedge)^\diamond$ represents the flat moduli space of $\msf{K}^+_p$-level structures of type $\msf{\ov{I}} = (\bb{P}^1)^n$. In particular, $\scr{S}_{\msf{K}^+}$ can be recovered by \Cref{prop:depth_zero_integral_models}, up to absolute weak normalization (cf.\ \cite[Theorem 2.16]{AGLR22}). 
\end{exa}

In a recent preprint, Pappas-Rapoport \cite{PR26} proposed axioms of $\Gamma_1(p)$-level integral models \cite[Conjecture 1.6.1]{PR26}. Our construction satisfies conditions (i) and (iii) loc. cit. 

Condition (ii) loc. cit. expects that $\scr{S}_{\msf{K}^+}$ should be finite over $\scr{S}_{\msf{K}}$ and the coarse quotient $\scr{S}_{\msf{K}^+} /\!/ \msf{I}^\sigma$ should equal $\scr{S}_{\msf{K}}$. In our construction, one way to deduce the finiteness of $\scr{S}_{\msf{K}^+}$ is to find an $\msf{I}$-equivariant affine open subset $\msf{U} \subset \msf{\ov{I}}$ so that \Cref{prop:finite_replevel} can be applied, which is the case in the proof of \Cref{prop:Gamma_1(p)-level}. However, one cannot expect this kind of finiteness in general, at least outside $\Gamma_1(p)$-levels. 

\begin{exa} \label{exa:Harris_Taylor}
    Let $\scr{S}_{\msf{K}}$ be an integral model of a Harris-Taylor Shimura variety considered in \cite[Section III.4]{HT01}, where $G = \GL_n \times \bb{G}_m$ and $\msf{K}_p = \GL_n(\bb{Z}_p) \times \bb{Z}_p^\times$. Then, there is a closed point $x \in \scr{S}_{\msf{K}}(\ov{k})$ such that $(\scr{S}_{\msf{K}})^\wedge_{/x}$ is isomorphic to $\LT_n$. Let $\msf{K}^+_p = (1 + p \, \Mat_n(\bb{Z}_p)) \times \bb{Z}_p^\times$. By \Cref{prop:Drlevel}, the moduli $\scr{S}_{\msf{K}^+}^1$ of depth-zero Drinfeld level structures over $\scr{S}_{\msf{K}}$ is essentially recovered by \Cref{prop:replevel} if we set $\ov{\GL}_n = \bb{P}(\Mat_n \oplus \triv)$. However, if we apply \Cref{prop:depth_zero_integral_models} to $\ov{\GL}_n = \KGL_n$, the completion of $\scr{S}_{\msf{K}^+}^{\KGL_n}$ at the fiber of $x$ recovers $\LT_n^{\KGL_n}$ up to normalization as in \Cref{cor:indep_normalization}. We will show in \Cref{sec:comparison_Yoshida} that $\LT_n^{\KGL_n}$ equals Yoshida's generalized semistable model of $\LT_n^1$ up to normalization. We expect that $\scr{S}_{\msf{K}^+}^{\KGL_n}$ should recover the successive blowup of $\scr{S}_{\msf{K}}^1$ introduced in \cite[Section 4.2]{Yos10} up to normalization, and in particular, the singularity of $\scr{S}_{\msf{K}^+}^{\KGL_n}$ should be milder than that of $\scr{S}_{\msf{K}^+}^1$ as noted in \cite[Remark 4.9]{Yos10}. 
\end{exa}

\section{Comparison with Yoshida's model} \label{sec:comparison_Yoshida}

\subsection{Review of Yoshida's construction} \label{ssec:Yoshida}

In this section, we recall the construction of generalized semistable models of depth-zero Lubin-Tate spaces introduced in \cite{Yos10}. 

Recall the notation in \Cref{ssec:LTn_KGLn}. Let $\alpha^\univ\colon (\bb{F}_p)^{\oplus n} \to \bb{X}^\univ[p]$ be the universal depth-zero Drinfeld level structure over $R^1_{\bb{X}}$. The $i$-th component $\bb{F}_p \to \bb{X}^\univ[p]$ is denoted by $\alpha_i^\univ$. In the following, $(a_i)_{1 \leq i \leq n}$ runs through elements of $\bb{F}_p^{\oplus n} - \{0\}$ and let $(a_i)^{\perp} \subset \bb{P}^{n-1}$ be the $\bb{F}_p$-rational hyperplane defined by $\sum_{1\leq i \leq n} a_i T_i = 0$. 

\begin{defi}
    Let $Y_{(a_i)} \subset \Spec(R^1_{\bb{X}})$ be the vanishing locus of $\sum_{1\leq i \leq n} a_i \alpha_i^\univ$. As in \cite[Section 3.2]{Yos10}, $Y_{(a_i)}$ is represented by an $(n-1)$-dimensional regular local ring in characteristic $p$.
\end{defi}

\begin{defi}
    For each $\bb{F}_p$-rational linear subspace $N\subset \bb{P}^{n-1}$, let $Y_N \subset \Spec(R^1_{\bb{X}})$ be the intersection of $Y_{(a_i)}$ for all $(a_i)$ with $N \subset (a_i)^\perp$. For $1\leq h < n$, let $Y^{[h]} \subset \Spec(R^1_{\bb{X}})$ be the union of $Y_N$ for all $N$ with $\dim N = h-1$.  
\end{defi}

The generalized semistable model of $\LT_{n,\eta}^1$ is constructed as a successive blowup of $\LT_n^1$. First, let $Z_1 \to \Spec(R^1_{\bb{X}})$ be the blowup along the maximal ideal of $R^1_{\bb{X}}$. Inductively on $1\leq h < n$, we define $Z_{h+1} \to Z_h$ to be the blowup along the strict transform of $Y^{[h]}$. Then, we have a sequence
\[Z_n \to Z_{n-1} \to \cdots \to Z_1 \to Z_0 = \Spec(R^1_{\bb{X}}).\]
Note that the strict transform of $Y^{[h]}$ in $Z_{h}$ is a regular closed subscheme of codimension $n-h$ (see \cite[\S 4.1]{Yos10}). In particular, $Z_n \to Z_{n-1}$ is an isomorphism. 

\begin{defi}
    Let $\LT_n^{1,i}$ be the $\mfr{m}$-adic completion of $Z_i$ for each $0\leq i \leq n$. The formal scheme $\LT_n^{1,n-1}$ is the generalized semistable model of $\LT_{n,\eta}^1$ introduced in \cite{Yos10}. 
\end{defi}

Now, $\LT_{n,\red}^{1,1} \cong \bb{P}^{n-1}$ via an explicit coordinate given by the formal parameters of $\alpha^\univ$. Then, $\LT_{n,\red}^{1,n-1}$ is a successive blowup of $\bb{P}^{n-1}$ along rational linear subspaces $N$. In particular, it is isomorphic to a compactified Deligne-Lusztig variety $\msf{\ov{X}}(w)$ (see \Cref{defi:cpt_DL}) by \cite[Section II, Proposition 1.6]{Gen96}. 

\subsection{Review of our construction} \label{ssec:ourconst}

In this section, we specialize the results in \Cref{sec:dep0_integral_model} to $\LT_n^{\KGL_n}$. We set $F=\bb{Q}_p$ and $\cl{G}=\GL_n$. Let $T \subset B \subset G$ be the opposite of the standard Borel pair. Let $\mu(t) = \diag(t^{-1}, 1,\ldots,1)$ and $b = \mu(-p)w$ with $w=\begin{pmatrix} 0 & 1 \\ I_{n-1} & 0 \end{pmatrix}$. Here, $I_{n-1}$ denotes the identity matrix of size $(n-1)\times (n-1)$. By \cite[Theorem 25.1.2]{SW20}, $\cl{M}^\ints_{\cl{G},b,\mu}$ is represented by the moduli space of $p$-divisible groups with quasi-isogenies to $\bb{X}$ (see \cite{RZ96}). %Since formal groups over $\bb{\ov{F}}_p$ are uniquely determined by their heights, 
We may identify the $p$-divisible group at $[1] \in X_\mu(b)$ with $\bb{X}$. Then, we have $(\cl{M}^\ints_{\cl{G},b,\mu})^{\wedge}_{/[1]} \cong \LT_n$. %, that is, $\LT_n^1 \cong \Spf(R_{\cl{G},\mu})$.  
%Recall that $\LT_n^{\KGL_n}$ is the $\mfr{m}$-adic completion of the flat moduli space of depth-zero level structures of type $\KGL_n$ on $\bb{X}^\univ$. 
By construction, $\LT_n^{\KGL_n}$ equals $\mfr{X}_b^{\KGL_n}$ in \Cref{const:Xb}. In the current setting, $e=p^n-1$ and $(e\lambda)(t) = \diag(t, t^{p}, \ldots, t^{p^{n-1}})$. In particular, $\msf{P} = \msf{G}^{e\lambda \geq 0}$ equals $\msf{B}$. The following property of $\KGL_n$ is essential to our argument. 

\begin{lem}
    The stabilizer of the limit point $x_\lambda \in \KGL_n$ is $\msf{\ov{B}}\times \msf{B}$. 
\end{lem}
\begin{proof}
    The map $\KGL_n \to \bb{P}(\Mat_n \oplus \triv)$ sends $x_\lambda$ to $0 \in \Mat_n$, so $x_\lambda$ lies in the fiber of $\KGL_n \to  \bb{P}(\Mat_n \oplus \triv)$ at $0$. The fiber is isomorphic to $\ov{\PGL}_n$ by \cite[Proposition 10.1]{Kau00}. Since $\msf{P} = \msf{B}$, $x_\lambda$ lies in the closed orbit of $\ov{\PGL}_n$, so the claim follows. 
\end{proof}

\begin{defi} \label{defi:cpt_DL}
    Let $\msf{\ov{X}}(w) \subset \msf{G}/\msf{B}$ be the Zariski closure of the Deligne-Lusztig variety $\msf{X}(w)$ as defined in \Cref{rmk:Xwopen}. It is also known as a compactified Deligne-Lusztig variety. 
\end{defi}

\begin{lem} \label{lem:inclusion_cpt_DL}
    The inclusion $(\LT_n^{\KGL_n})^\perf_\red \subset \ov{\msf{G}}^\perf$ induces $\msf{\ov{X}}(w)^\perf \subset (\LT_n^{\KGL_n})^\perf_\red$. 
\end{lem}
\begin{proof}
    By \Cref{rmk:Xwopen}, $\msf{X}(w)^\perf \subset (\LT_n^{\KGL_n})^\perf_\red$ is an open subspace, so the claim follows since $(\LT_n^{\KGL_n})^\perf_\red \subset \ov{\msf{G}}^\perf$ is closed. 
\end{proof}

%The Zariski closure $\msf{\ov{X}}(w)$ of $\msf{X}(w) \subset \msf{G}/\msf{B}$ is known as a compactified Deligne-Lusztig variety. In our case, $\msf{\ov{X}}(w)$ is realized as a successive blowup of $\bb{P}^{n-1}$ and admits a moduli interpretation as a subspace of the flag variety $\msf{G}/\msf{B}$ (see \cite[Section II, Proposition 1.6]{Gen96}). 

The comparison between $\LT_n^{1,n-1}$ and $\LT_n^{\KGL_n}$ will be deduced from \cite[Corollary 6.6]{Lou17} in the end. For this, we calculate and compare their specialization maps. These calculations are done using the prismatic Dieudonn\'{e} theory. As in \Cref{ssec:Drinfeld}, we employ a covariant Dieudonn\'{e} theory and adopt the naive dual of the contravariant one developed in \cite{ALB23}. When $n = 1$, the claim follows from \Cref{prop:LTn1}. From now on, we assume $n\geq 2$. 

\subsection{Specialization in $\LT_n^{\KGL_n}$} \label{sssec:spKGLn}

In this section, we compute the specialization map of $\LT_n^{\KGL_n}$. Let $(C,C^+)$ be an algebraically closed perfectoid field over $\breve{F}$ and take a geometric point $x\in \LT_n^{\KGL_n}(C,C^+)$. Let $p^\flat$ and $u_\alpha^\flat$ be as in \Cref{lem:expBKF}. Let $\nu$ and $\nu^\flat$ be additive valuations on $C$ and $C^\flat$, respectively, such that $\nu(p)=1$ and $\nu(f^\sharp) = \nu^\flat(f)$ for every $f\in C^\flat$. 

By \Cref{lem:expBKF}, the Dieudonn\'{e} module of the $p$-divisible group over $C^+$ at $x$ is given by $(W(C^{\flat+})^{\oplus n}, \mu([p^\flat]-p)\prod_{\alpha\in \Phi_{\mu<0}} i_\alpha([u_\alpha^\flat])w\sigma)$. By \Cref{cor:explevel}, the depth-zero level structure of $x$ is given by $g\in \GL_n(C^{\flat})$ such that $g\sigma(g)^{-1} = \mu(p^\flat)\prod_{\alpha \in \Phi_{\mu<0}} i_{\alpha}(u_\alpha^\flat)w$. Let $e_1,\ldots,e_n$ be the standard basis of $(C^{\flat})^{\oplus n}$ and set $z_i = e_i \cdot g$. For each $2\leq i \leq n$, $\alpha_i \in \Phi_{\mu<0}$ denotes the element with $e_i \cdot (\alpha_i^\vee)(t) = t^{-1}e_i$ and let $u_i^\flat = u_{\alpha_i}^\flat$. Then, the condition on $g$ is written as
\[
    p^\flat z_1 = \sigma(z_n) + \sum_{2\leq i \leq n} u_i^\flat \sigma(z_{i-1}), \quad z_i = \sigma(z_{i-1}) \quad (2\leq i \leq n). 
\]
In particular, $z_i = \sigma^{i-1}(z_1)$ and $p^\flat z_1 = \sigma^n(z_1) + \sum_{1\leq i \leq n-1} u_{i+1}^\flat \sigma^i(z_1)$. Since $p^\flat, u_\alpha^\flat\in C^{\flat, \circ \circ}$, we have $z_1 \in (C^{\flat, \circ \circ})^{\oplus n}$. 

\begin{lem}
    We can take a set of elements $(\varpi_i)_{1\leq i \leq n}, (a_{j,i})_{1\leq i < j \leq n}$ in $C^{\flat,\circ \circ}$ so that 
    \begin{enumerate}
        \item $\nu^\flat(\varpi_1) < \nu^\flat(\varpi_2) < \cdots  < \nu^\flat(\varpi_n)$, 
        \item $z'_i = \varpi_i^{-1} (z_i + \sum_{1\leq j < i} a_{i,j} z_j) \in (C^{\flat +})^{\oplus n}$, and
        \item $(C^{\flat +})^{\oplus n} = \bigoplus_{1\leq i \leq n} C^{\flat +} z'_i$. 
    \end{enumerate}
\end{lem}
\begin{proof}
    We construct $\varpi_i$ and $z'_i$ inductively on $i$. Suppose that we have $\varpi_1,\ldots,\varpi_i$ and $z'_1,\ldots,z'_i$ so that $\nu^\flat(\varpi_1) <\cdots<\nu^\flat(\varpi_i)$ and $\bigoplus_{1\leq j \leq i} C^{\flat +} z'_i$ is a direct summand of $(C^{\flat +})^{\oplus n}$. Since $z_{i+1} = \sigma(z_i)$, $z''_{i+1} = \sigma(\varpi_i)^{-1}(z_{i+1}+\sum_{1\leq j < i} \sigma(a_{i, j})z_{j+1})$ lies in $(C^{\flat +})^{\oplus n}$. Let $M=(C^{\flat +})^{\oplus n}/\bigoplus_{1\leq j \leq i} C^{\flat +} z'_j$. Since $z''_{i+1}$ is linearly independent of $z'_1,\ldots,z'_i$ over $C^\flat$ and $M$ is a free $C^{\flat+}$-module, we can take a nonzero element $c \in C^{\flat+}$ such that $C^{\flat+}\cdot c^{-1}[z''_{i+1}]  = M \cap C^\flat [z''_{i+1}]$ where $[z''_{i+1}]$ is the image of $z''_{i+1}$ in $M$. Then, we can take $c_1,\ldots,c_i \in C^{\flat +}$ so that $z'_{i+1} = c^{-1}(z''_{i+1} + \sum_{1\leq j \leq i} c_j z'_j)$ and $\varpi_{i+1} = c\sigma(\varpi_i)$ satisfy the induction hypothesis for $i+1$. 
\end{proof}

Thus, we have a Cartan decomposition
\[g = \begin{pmatrix} 1 &  &  & \\ a_{2,1} & 1 &  & \\ \vdots & \cdots & \ddots &  \\ a_{n,1} & \cdots & a_{n,n-1} & 1 \end{pmatrix}^{-1} \begin{pmatrix} \varpi_1 & & & \\ & \varpi_2 & & \\ & & \ddots & \\ & & & \varpi_n \end{pmatrix} \begin{pmatrix} z'_1 \\ z'_2 \\ \vdots \\ z'_n \end{pmatrix}. \]
In particular, the specialization of $x$ lies in $x_{\lambda}\cdot (\msf{B}\backslash \msf{G})$. From now on, we will describe the flag given by $z'_1,\ldots,z'_n$ more precisely.

Let $\kappa$ be the residue field of $C^{\flat+}$. The $i$-th term of the flag is determined by the reduction of $z'_1 \wedge \cdots \wedge z'_i$ to $\kappa$. Since $(\prod_{1\leq j \leq i} \varpi_j) z'_1 \wedge \cdots \wedge z'_i = z_1 \wedge \cdots \wedge z_i$, it is enough to calculate $z_1 \wedge \cdots \wedge z_i$ up to scalar. Let $z_1 = (t_1,\ldots,t_n)$. By taking a right translation of $g$ under $\GL_n(\bb{F}_p)$, we may assume that there is a sequence $0 = i_0 < i_1 < \cdots < i_k = n$ such that 
\begin{enumerate}
    \item $\nu^\flat(t_s) = \nu^\flat(t_{i_j})$ for every $i_{j-1} < s \leq i_j$ and the reduction of $\{t_s/t_{i_j}\}_{i_{j-1}<s\leq i_j}$ to $\kappa$ is linearly independent over $\bb{F}_p$ for every $1 \leq j \leq k$, and 
    \item $\nu^\flat(t_{i_j}) < \nu^\flat(t_{i_j+1})$ for every $1 \leq j < k$. 
\end{enumerate}

For each $m\geq 0$, let $\msf{DL}^m \subset \bb{P}^m$ be the Deligne-Lusztig variety as in \cite[Section 2.2]{DL76}. For each $1\leq j \leq k$, let $\ell_j = i_j-i_{j-1}$ and $p_j = [t_{i_{j-1}+1} \colon \cdots \colon t_{i_j-1} \colon t_{i_j}] \in \bb{P}^{\ell_j-1}(C^\flat)$. Let $\ov{p}_j\in \bb{P}^{\ell_j-1}(\kappa)$ be the specialization of $p_j$. Then $\ov{p}_j \in \msf{DL}^{\ell_j-1}$ by condition (1). Moreover, there is a natural inclusion $\msf{DL}^{\ell_1-1} \times \msf{DL}^{\ell_2-1}\times \cdots \times \msf{DL}^{\ell_k-1} \subset \msf{\ov{X}}(w)$ and $(\ov{p}_1,\ov{p}_2,\ldots,\ov{p}_k)$ maps to a flag $(\Lambda_i)_{1\leq i \leq n}$ such that for each $1\leq j \leq k$ and $i_{j-1} < s \leq i_j$, $\Lambda_s$ is generated by $\Lambda_{s-1}$ and $\sigma^{s-i_{j-1}-1}(\ov{p}_j)$. Note that $\Lambda_0=0$ and $\sigma(\Lambda_{i-1}) \subset \Lambda_i$ for every $1\leq i \leq n$. 

\begin{prop} \label{prop:spKGL}
    The specialization of $x$ in $\LT_n^{\KGL_n}$ is $(\ov{p}_1,\ov{p}_2,\ldots,\ov{p}_k) \in \prod_{1\leq j \leq k} \msf{DL}^{\ell_j-1}$. 
\end{prop}

\begin{proof}

For this proof, we need the following two lemmas. 

\begin{lem} \label{lem:nonsing}
    Let $a_1,\ldots,a_n \in \kappa$ be a set of elements that is linearly independent over $\bb{F}_p$. Then, the matrix $(\sigma^{i-1}(a_j))_{1\leq i,j \leq n}$ is nonsingular. 
\end{lem}
\begin{proof}
    Suppose that $(\sigma^{i-1}(a_j))_{1\leq i,j \leq n}$ is singular. Then, there exists $(b_1,\ldots,b_n) \in \kappa^{\oplus n}-\{0\}$ such that $\sum_{1\leq i \leq n} b_i \sigma^{i-1}(a_j)=0$ for every $1\leq j \leq n$. Let $P(x) = \sum_{1\leq i\leq n} b_i x^{p^{i-1}}$. Then, $P(a) = 0$ for every $a=\sum_{1\leq j \leq n} c_j a_j$ with $c_j \in \bb{F}_p$. Then, $P(x)=0$ has $q^n$ different solutions as $a_1,\ldots,a_n$ is linearly independent over $\bb{F}_p$, but it contradicts $\deg P(x) \leq p^{n-1}$. 
\end{proof}

    For each $1\leq j \leq k$ and $i_{j-1} < s \leq i_j$, take $z^s =(t_1^s,\ldots,t_n^s)\in (C^{\flat+})^{\oplus n}$ so that $t_i^s = t_i/t_s$ for $i_{j-1} < i \leq i_j$ and $t_i^s = 0$ otherwise. By \Cref{lem:nonsing}, $\{ \sigma^{s-1}(z^s) \}_{1 \leq s \leq n}$ is a basis of $(C^{\flat+})^{\oplus n}$. 

\begin{lem} \label{lem:orderwedge}
    For each $1\leq i \leq n$, let $\varpi'_i = \prod_{1\leq j \leq i} \sigma^{i-j}(t_j)$. Then, $(\varpi'_i)^{-1} (z_1 \wedge \cdots \wedge z_i) \in \bigwedge^i (C^{\flat+})^{\oplus n}$ and its reduction to $\kappa$ is equal to that of $z^i\wedge \sigma(z^{i-1}) \wedge \cdots \wedge \sigma^{i-1}(z^1)$. 
\end{lem}
\begin{proof}
    For $1\leq j \leq k$ and $1\leq s \leq i$ with $i_{j-1} < s \leq i_j$, $z_{i-s+1} = \sigma^{i-s}(z_1) = \sigma^{i-s}(t_s)\sigma^{i-s}(z^s) +     \sum_{1\leq \ell \leq k, \ell \neq j}  \sigma^{i-s}(t_{i_\ell}) \sigma^{i-s}(z^{i_\ell})$. The claim follows by substituting this into $z_1\wedge \cdots \wedge z_i$. %Here, we use conditions (1) and (2) on $(t_1,\ldots,t_n)$ and \Cref{lem:nonsing}. 
\end{proof}
%\begin{proof}[Proof of \Cref{prop:spKGL}]

    By \Cref{lem:nonsing}, the reduction of $z^i\wedge \sigma(z^{i-1}) \wedge \cdots \wedge \sigma^{i-1}(z^1)$ is nonzero in $\bigwedge^i \kappa^{\oplus n}$. Thus, the specialization of $x$ is given by the reduction of the flag $\{z^i\wedge \sigma(z^{i-1}) \wedge \cdots \wedge \sigma^{i-1}(z^1)\}_{1\leq i \leq n}$. It can be directly checked that it is equal to the flag $(\Lambda_i)$. 
\end{proof}

\begin{cor}
    The underlying space $(\LT_n^{\KGL_n})^\perf_\red$ is equal to $\msf{\ov{X}}(w)^\perf$. 
\end{cor}
\begin{proof}
    Since $\LT_n^{\KGL_n}$ is $p$-torsion free, the specialization map of $\LT_n^{\KGL_n}$ is surjective. Thus, the claim follows from \Cref{lem:inclusion_cpt_DL} and \Cref{prop:spKGL}. 
\end{proof}

\subsection{Specialization in $\LT_n^{1,n-1}$} \label{sssec:spcLT}
In this section, we compute the specialization map of $\LT_n^{1,n-1}$. In \cite{Yos10}, it is studied via formal parameters of $\alpha^\univ$. Instead, we give a Dieudonn\'{e} theoretic description of $\LT_n^1$. 

For a $p$-divisible group $X$ over a formal scheme over $\Spf(\bb{Z}_p)$, let $E(X)$ denote the universal vectorial extension of $X$ (see \cite{Mes72}) and let $M(X) = \Lie(E(X))$. By the theory of \cite{Mes72} and \cite{BBM82}, $M(X)$ extends to a crystal $M_\crys$ on the crystalline site. 

As in \Cref{sssec:spKGLn}, $u_{\alpha_i}$ is denoted by $u_i$ for $2\leq i \leq n$. By \cite{Ito25a} and \cite{Ito25b}, the prismatic Dieudonn\'{e} crystal of $\bb{X}^\univ$ is determined by the evaluation at $(\breve{\bb{Z}}_p \llbracket t, u_2,\ldots, u_n \rrbracket, (t-p))$ and given by $(\breve{\bb{Z}}_p \llbracket t, u_2,\ldots, u_n \rrbracket^{\oplus n}, \mu(t-p)\prod_{2\leq i \leq n}i_{\alpha_i}(u_i)w)$. By applying \cite[Lemma 4.45]{ALB23} to the crystalline prism $(\breve{\bb{Z}}_p \llbracket u_2,\ldots, u_n \rrbracket, (p))$, we obtain a natural isomorphism
\[M(\bb{X}^\univ \otimes \ov{\bb{F}}_p \llbracket u_2,\ldots, u_n \rrbracket) \cong \sigma^* \ov{\bb{F}}_p \llbracket u_2,\ldots, u_n \rrbracket^{\oplus n}.\]
Let $\pr$ denote the projection from $M(\bb{X}^\univ \otimes \ov{\bb{F}}_p \llbracket u_2,\ldots, u_n \rrbracket)$ to the first direct summand on the right-hand side. %Fix a lift $\pr'_i \colon M(\bb{X}^\univ) \to \breve{\bb{Z}}_p \llbracket u_2, \ldots, u_n \rrbracket$ of $\pr_i$. 

\begin{rmk}
    If $(A,I) \in (R_{\bb{X}})_\Prism$ is not crystalline, the relation between the prismatic Dieudonn\'{e} crystal at $(A,I)$ and $M(\bb{X}^\univ\otimes (A/I))$ is not clear (cf.\ \cite[Remark 4.41]{ALB23}). 
\end{rmk}

The Cartier dual of $\alpha_i^\univ\otimes R^1_{\bb{X}}/pR^1_{\bb{X}} \colon \bb{F}_p \to \bb{X}^\univ[p]\otimes R^1_{\bb{X}}/pR^1_{\bb{X}}$ induces
\[\omega_{\mu_{p}} \to \omega_{(\bb{X}^{\univ}[p]\otimes R^1_{\bb{X}}/pR^1_{\bb{X}})^{\vee}} \cong \omega_{(\bb{X}^{\univ})^{\vee}} \otimes R^1_{\bb{X}}/pR^1_{\bb{X}} \to M(\bb{X}^\univ\otimes \ov{\bb{F}}_p \llbracket u_2,\ldots,u_n \rrbracket) \otimes R^1_{\bb{X}}/pR^1_{\bb{X}}.\] 
Here, $\omega_{(-)}$ denotes the dual of $\Lie(-)$. Fix an identification $\omega_{\mu_p} \cong M(\bb{Q}_p/\bb{Z}_p) \cong \bb{F}_p$ and let $\beta_i$ be the image of $1 \in \omega_{\mu_p}$. Let $t'_i\in R^1_{\bb{X}}$ be a lift of $\pr(\beta_i)$. 

\begin{prop} \label{prop:parameter}
    The elements $t'_1,\ldots,t'_n$ generate the maximal ideal of $R^1_{\bb{X}}$. 
\end{prop}

\begin{proof}

For this proof, we take a geometric point $x\in \LT_{n,\eta}^1(C,C^+)$. Since $\LT_{n,\eta}^{\KGL_n} \cong \LT_{n,\eta}^1$, we may keep the notation in \Cref{sssec:spKGLn}. By abuse of notation, the image of $t'_1,\ldots,t'_n$ in $C^+$ is also denoted by $t'_1,\ldots,t'_n$. We will study the relation between $t'_1,\ldots,t'_n$ and $z_1=(t_1,\ldots,t_n)$. Note that we do not assume conditions (1) and (2) on $(t_1,\ldots,t_n)$ at this point. 

\begin{lem} \label{lem:compval}
    We have $t'_i \equiv t_i^p$ under the identification $C^+/(p) \cong C^{\flat+}/(p^\flat)$. 
\end{lem}
\begin{proof}
    Since $C$ is algebraically closed, $\alpha_i^\univ \otimes C^{+}$ can be extended to a homomorphism $\bb{Q}_p/\bb{Z}_p \to \bb{X}^{\univ}$ of $p$-divisible groups over $C^{+}$. By passing it to a map of prismatic Dieudonn\'{e} modules and evaluating it at the crystalline prism $(W(C^{\flat+})/[p^\flat], (p))$, we obtain a commutative diagram from \cite[Lemma 4.45]{ALB23}: 
    \begin{center}
        \begin{tikzcd}
            \omega_{\mu_p} \otimes (C^{+}/p) \ar[r, "\cong"] \ar[d, "(\alpha_i^{\univ,\vee})^*"] & M(\bb{Q}_p/\bb{Z}_p) \otimes (C^{+}/p) \ar[r, "\cong"] \ar[d] & C^{\flat+}/p^\flat \ar[d, "(t_i^p\,\ldots\,t_i^{p^n})"] \\
            \omega_{(\bb{X}^{\univ})^{\vee}} \otimes (C^{+}/p) \ar[r, hook] & M(\bb{X}^\univ) \otimes (C^{+}/p) \ar[r, "\cong"] & (C^{\flat+}/p^\flat)^{\oplus n}. 
        \end{tikzcd}
    \end{center}
    The claim follows by taking a projection to the first direct summand of $(C^{\flat+}/p^\flat)^{\oplus n}$. 
\end{proof}

\begin{lem} \label{lem:minord}
    $\min_{1\leq i \leq n} \nu(t'_i) = \min_{1\leq i \leq n} \nu^\flat(t_i^p) < \tfrac{1}{p-1}$. 
\end{lem}
\begin{proof}
    We may assume that $(t_1,\ldots,t_n)$ satisfies conditions (1) and (2) in \Cref{sssec:spKGLn} and use the notation loc. cit. As $\nu^\flat(\det(g\sigma(g)^{-1}))=\nu^\flat(\det(\mu(p^\flat)))=-1$, $\nu^\flat(\sigma^{n-i}(t_i)) < \tfrac{1}{p-1}$ for $1\leq i \leq n$ by \Cref{lem:orderwedge} as $n\geq 2$. In particular, $\nu^\flat(t_i^p) < \tfrac{1}{p-1}$ for $1\leq i < n$, so the claim follows from \Cref{lem:compval}. 
\end{proof}
%\begin{proof}[Proof of \Cref{prop:parameter}]
    Let $T_i \in R^1_{\bb{X}}$ be a formal parameter of $\alpha_i^\univ$ for $1\leq i \leq n$. By \cite{Dr74}, $T_1,\ldots,T_n$ generate the maximal ideal of $R^1_{\bb{X}}$. For each geometric point $x \in \LT_{n,\eta}^1(C,C^+)$ and $\varpi \in C^+$, $\alpha^\univ\otimes (C^+/\varpi)$ is trivial if and only if the images of $T_1,\ldots,T_n$ are divisible by $\varpi$. 

    On the other hand, if $\varpi$ divides $p$, it follows from the proof of \Cref{lem:compval} that if $\alpha_i^{\univ} \otimes (C^+/\varpi)$ is trivial, then $\varpi$ divides $t'_i$. We will show that the converse also holds. Suppose that $\varpi$ divides $t'_i$. Then, $M(\bb{Q}_p/\bb{Z}_p) \otimes (C^{+}/\varpi)  \to  M(\bb{X}^\univ) \otimes (C^{+}/\varpi)$ is zero. Take $\varpi^\flat \in C^{\flat+}$ so that $(\varpi^\flat)^\sharp = \varpi$. Then, $M(\bb{Q}_p/\bb{Z}_p)(W(C^{\flat+})/[\varpi^{\flat}]) \to M(\bb{X}^\univ)(W(C^{\flat+})/[\varpi^{\flat}])$ factors through the multiplication by $p$. Then, $\bb{Q}_p/\bb{Z}_p \to \bb{X}^{\univ}$ defined over $C^+$ factors through the multiplication by $p$ after the base change to $C^+/\varpi$ by \cite[Theorem 5.7]{Lau18}. It implies that $\alpha_i^{\univ} \otimes (C^+/\varpi)$ is trivial. By \Cref{lem:minord}, it follows that $\alpha^\univ\otimes (C^+/\varpi)$ is trivial if and only if $t'_1,\ldots,t'_n$ are divisible by $\varpi$. 

    In particular, $\min_{1\leq i \leq n} \nu(t'_i) = \min_{1\leq i \leq n} \nu(T_i)$ for every geometric point $x\in \LT_{n,\eta}^1(C,C^+)$. Let $r = \dim_{\ov{\bb{F}}_p } (t'_1,\ldots,t'_n)/(t'_1,\ldots,t'_n) \cap (T_1,\ldots,T_n)^2$ and take a set of elements $f_1,\ldots,f_r \in (t'_1,\ldots,t'_n)$ so that $(t'_1,\ldots,t'_n) = (f_1,\ldots,f_r) +  (t'_1,\ldots,t'_n) \cap (T_1,\ldots,T_n)^2$. Since $\nu(T_i) > 0$, $\min_{1\leq i \leq r} \nu(f_i) = \min_{1\leq i \leq n} \nu(T_i)$ for every $x\in \LT_{n,\eta}^1(C,C^+)$.
    
    It is enough to show that $r=n$. Suppose $r<n$. The dimension of $R_{\bb{X}}^1/(f_1,\ldots,f_r)$ is $n-r>0$, so there is $\vert\cdot \vert \in \Spa(R^1_{\bb{X}})^\an$ such that $\lvert f_i \rvert =0 $ for every $1\leq i\leq r$. Take $1\leq i \leq n$ so that $\lvert T_j \rvert \leq \lvert T_i \rvert \neq 0$ for every $1\leq j \leq n$. Then, the rational localization $R(\tfrac{f_1,\ldots,f_r,T_1^2,\ldots,T_n^2}{T_i^2})\subset \Spa(R^1_{\bb{X}})^\an$ is nonempty, but it cannot meet $\LT_{n, \eta}^1$. Thus, we get a contradiction as $R^1_{\bb{X}}$ is $p$-torsion free. 
\end{proof}

\begin{cor} \label{lem:Yai}
    We have $Y_{(a_i)} = V(\sum_{1\leq i \leq n} [a_i] t'_i) \subset \Spec(R^1_{\bb{X}})$. 
\end{cor}
\begin{proof}
    By \cite[Section 3.2]{Yos10}, $Y_{(a_i)}$ is represented by an $(n-1)$-dimensional regular local ring $R_{(a_i)}$ in characteristic $p$. Then, $\sum_{1\leq i \leq n} a_i \alpha_i^\univ$ induces $\omega_{\mu_p} \to M(\bb{X}^\univ \otimes R_{(a_i)}) \xrightarrow{\pr} R_{(a_i)}$ and it maps $1$ to $\sum_{1\leq i \leq n} a_i t'_i$. In particular, $\sum_{1\leq i \leq n} [a_i] t'_i$ vanishes on $Y_{(a_i)}$. By \Cref{prop:parameter}, $R^1_{\bb{X}}/(\sum_{1\leq i \leq n} [a_i] t'_i) \to R_{(a_i)}$ is a surjection of regular local rings of the same dimension, so it is an isomorphism. 
\end{proof}

%To state some form of compatibility in \Cref{prop:intermlevel}, we will compare $\LT_n^{\KGL_n}$ with $\mfr{Z}_n$. 

In particular, the coordinate of $\LT_{n,\red}^{1,1} \cong \bb{P}^{n-1}$ can be taken as $t'_1,\ldots,t'_n$ and $\LT_{n,\red}^{1,n-1}$ is isomorphic to a successive blowup of $\bb{P}^{n-1}$ along rational linear subspaces with respect to this coordinate. 

Finally, we will compute the specialization map of $\LT_{n}^{1,n-1}$. Again, we take a geometric point $x\in  \LT_{n,\eta}^1(C,C^+)$ and keep the notation in \Cref{sssec:spKGLn}. Here, we will assume conditions (1) and (2) on $(t_1,\ldots,t_n)$ loc. cit. By the proof of \Cref{lem:minord}, $\nu^\flat(t_i^p) < \tfrac{1}{p-1}$ for $1\leq i < n$, so it follows from \Cref{lem:compval} that 
\begin{enumerate}
    \item $\nu(t'_s) = \nu(t'_{i_j})$ for every $i_{j-1} < s \leq i_j$ and the reduction of $\{t'_s/t'_{i_j}\}_{i_{j-1}<s\leq i_j}$ to $\kappa$ is linearly independent over $\bb{F}_p$ for every $1 \leq j \leq k$, and 
    \item $\nu(t'_{i_j}) < \nu(t'_{i_j+1})$ for every $1 \leq j < k$. 
\end{enumerate}
For each $1\leq j \leq k$, let $p'_j = [t'_{i_{j-1}+1} \colon \cdots \colon t'_{i_j-1} \colon t'_{i_j}] \in \bb{P}^{\ell_j-1}(C)$ and let $\ov{p}'_j\in \bb{P}^{\ell_j-1}(\kappa)$ be the specialization of $p'_j$. By condition (1), $\ov{p}'_j \in \msf{DL}^{\ell_j-1}$.

\begin{prop} \label{prop:spYos}
    The specialization of $x$ in $\LT_n^{1, n-1}$ is $(\ov{p}'_1,\ov{p}'_2,\ldots,\ov{p}'_k) \in \prod_{1\leq j \leq k} \msf{DL}^{\ell_j-1}$. 
\end{prop}
\begin{proof}
    The reduction of $x$ in $Z_1$ is $\ov{p}'_1\times\{0\} \in \bb{P}^{n-1}$ and the minimal rational linear subspace $N$ containing $\ov{p}'_1\times\{0\}$ is $\{t'_s = 0 \mid i_1 < s \leq n\}$ by condition (1). By \Cref{lem:Yai}, $Y_N = V(t'_{i_1+1},\ldots,t'_{n})$ and the exceptional divisor of the blowup of $\bb{P}^{n-1}$ along $\bb{P}^{n-1}\cap Y_N$ is isomorphic to $\bb{P}^{i_1-1} \times \bb{P}^{n-i_1-1}$. Moreover, the strict transform of $Y_M$ with $N\subset M$ in the blowup along $Y_N$ intersects the exceptional divisor $\bb{P}^{i_1-1} \times \bb{P}^{n-i_1-1}$ at $\bb{P}^{i_1-1} \times (\bb{P}^{n-i_1-1} \cap M)$. Thus, by repeating this argument, the specialization of $x$ lies in $\prod_{1\leq j \leq k} \msf{DL}^{\ell_j-1} \subset \msf{\ov{X}}(w)$ and the coordinate is $(\ov{p}'_1,\ov{p}'_2,\ldots,\ov{p}'_k)$. 
\end{proof}

\subsection{Summary} \label{ssec:summ}

By \Cref{prop:spKGL} and \Cref{prop:spYos}, we have a commutative diagram
\begin{center}
    \begin{tikzcd}
        \LT_n^{\KGL_n}(C,C^+) \ar[r, "\cong"] \ar[d, "\spc"] & \LT_n^{1,n-1}(C,C^+) \ar[d, "\spc"] \\
        \ov{\msf{X}}(w)(\kappa) \ar[r,"\sigma"] & \ov{\msf{X}}(w)(\kappa)
    \end{tikzcd}
\end{center}
for every algebraically closed perfectoid field $(C,C^+)$ over $\breve{\bb{Q}}_p$ with residue field $\kappa$. 

Let $\mfr{Y}$ be the flat moduli space of depth-zero level structures of type $\KGL_n$ on $\bb{X}^\univ$ and let $\mfr{Y}^\norm\to \mfr{Y}$ be the normalization in the generic fiber. Then, let $\LT_n^{\KGL_n,\norm} = (\mfr{Y}^\norm)^\wedge$ be the $\mfr{m}$-adic completion. By \cite[Theorem 4.2]{Yos10}, $\LT_n^{1,n-1}$ is a normal formal scheme formally of finite type over $\breve{\bb{Z}}_p$, so the maps
\[  
    \LT_{n, \eta}^{\KGL_n,\norm} \cong \LT_{n, \eta}^{1,n-1} ,\quad 
    (\LT_n^{\KGL_n,\norm})^\perf_\red \to (\LT_n^{\KGL_n})^\perf_\red = \ov{\msf{X}}(w) \xrightarrow{\sigma} \ov{\msf{X}}(w) = (\LT_{n}^{1,n-1})^\perf_\red
\]
compatible with specialization maps induce a morphism $\iota \colon \LT_n^{\KGL_n,\norm} \to \LT_n^{1,n-1}$ by \cite[Corollary 6.6]{Lou17}. Since $\iota$ is adic and $\iota_\red$ is finite, $\iota$ is finite itself. Moreover, $\iota_\eta$ is an isomorphism, so $\iota$ is an isomorphism. In particular, we get the following. 

\begin{thm} \label{thm:comparison}
    The normalization of $\LT_n^{\KGL_n}$ is isomorphic to $\LT_n^{1, n-1}$. Moreover, we have $\LT_{n,\eta}^{1,n-1}\cong \LT_{n,\eta}^{\KGL_n}$,  $(\LT_{n}^{1,n-1})^\perf_\red\cong (\LT_n^{\KGL_n})^\perf_\red$ and $\LT_n^{1,n-1}\vert_{\msf{X}(w)}\cong \LT_n^{\KGL_n}\vert_{\msf{X}(w)}$. % as $p$-adic formal schemes. 
\end{thm}
\begin{proof}
    It is enough to prove the last claim. The proof of \Cref{thm:role_of_integralmodel} shows that $\msf{X}(w)^\perf \hookrightarrow \mfr{Y}^\perf_\red$ is an open immersion. In particular, $\LT_n^{\KGL_n}\vert_{\msf{X}(w)} \cong \mfr{Y}\vert_{\msf{X}(w)}$ and $\LT_n^{\KGL_n}\vert_{\msf{X}(w)}$ is $p$-adic. It follows from \cite[Theorem 2.37]{AGLR22} that $(\LT_n^{1,n-1}\vert_{\msf{X}(w)})^\diamond \cong (\LT_n^{\KGL_n}\vert_{\msf{X}(w)})^\diamond$. Since $\mfr{Y}$ admits a maximal good cover by construction, it follows from \Cref{lem:uniqmapmaxlgood} that we have an inverse $\LT_n^{\KGL_n}\vert_{\msf{X}(w)} \to \LT_n^{1,n-1}\vert_{\msf{X}(w)}$. 
\end{proof}

\begin{rmk} \label{rmk:expectisom}
    In fact, we expect $\LT_n^{\KGL_n} \cong \LT_n^{1,n-1}$. We present two possible ways to show this, but we do not pursue the generality here to simplify the computations of specialization maps. 
    \begin{enumerate}
        \item Compare the specialization maps of $\mfr{Y}$ and $Z_n$. 
        \item Work over the entire Shimura variety as in \cite[Section 4.2]{Yos10} instead of $\LT_n$. 
    \end{enumerate}
    In both cases, non-basic points appear in the special fibers and we need to compare their fibers in the specialization maps. The benefit of these calculations is that both formal schemes being compared are $p$-adic, so we may deduce the claim from \cite[Theorem 2.37]{AGLR22} together with \Cref{lem:uniqmapmaxlgood}. 
\end{rmk}

\begin{rmk}\label{prop:intermlevel}
    The construction of the successive blowups $\{Z_i\}$ of $\Spec(R^1_{\bb{X}})$ is parallel to the construction of the successive blowups $\{X^{(i)}\}$ of $\bb{P}(\Mat_n \oplus \triv)$ (see \Cref{ssec:cptGLn}). For $0 \leq i \leq n$, let $\LT_n^{(i)}$ be the $\mfr{m}$-adic completion of the flat moduli space of depth-zero level structures of type $X^{(i)}$ on $\bb{X}^\univ$. For $0 \leq i \leq n-1$, we can verify that the stabilizer of the limit point $x_\lambda$ in $X^{(i)}$ with respect to the right $\GL_n$-action is equal to $\GL_n^{\lambda_i \geq 0}$ with $\lambda_i = (t,t^p,\ldots,t^{p^{i}}, \ldots,t^{p^i})$. Then, $\KGL_n \to X^{(i)}$ induces $x_\lambda \cdot \msf{B} \backslash \msf{G} \to x_\lambda \cdot \GL_n^{\lambda_i \geq 0}\backslash \msf{G}$. For every $g\in \GL_n$, $x_\lambda \cdot g \in \GL_n^{\lambda_i \geq 0}\backslash \GL_n$ corresponds to a flag $(z_1,z_1\wedge z_2,\ldots,z_1\wedge \cdots \wedge z_i)$ where $z_j = e_j \cdot g$ for $1\leq j \leq i$ as in \Cref{sssec:spKGLn}. On the other hand, $\LT_{n,\red}^{1,i}$ parametrizes flags $(\Lambda_j)_{1\leq j \leq i}$ of length $i$ with $\sigma(\Lambda_{j-1}) \subset \Lambda_{j}$ (see the proof of \cite[Section II, Proposition 1.6]{Gen96}). Then, the compatibility of specialization maps of $\LT_n^{(i)}$ and $\LT_n^{1,i}$ follows from \Cref{ssec:summ}, and the normalization of $\LT_n^{(i)}$ is isomorphic to $\LT_n^{1, i}$ as in \Cref{thm:comparison}. 
\end{rmk}

\renewcommand\bibfont{\footnotesize}
\printbibliography

\end{document}